\documentclass[10pt,psamsfonts]{article}

\usepackage{amsmath, amsfonts, amsthm, amssymb, enumerate, url, float}
\usepackage{amscd, enumerate, url, slashed, mathabx}
\usepackage[all,arc,color]{xy}
\usepackage{mathrsfs}
\usepackage{tikz-cd}
\usepackage{mathtools}
\usepackage[utf8]{inputenc}
\usepackage{graphicx}
\graphicspath{ {./images/} }
\usepackage{tikz-cd}

\begingroup
    \makeatletter
    \@for\theoremstyle:=definition,remark,plain\do{%
        \expandafter\g@addto@macro\csname th@\theoremstyle\endcsname{%
            \addtolength\thm@preskip\parskip
            }%
        }
 \endgroup
 
\usepackage[titletoc,toc,title]{appendix}
 
\usepackage{musicography}

\makeatletter
\newcommand*{\itemequation}[3][]{%
  \item
  \begingroup
    \refstepcounter{equation}%
    \ifx\\#1\\%
    \else  
      \label{#1}%
    \fi
    \sbox0{#2}%
    \sbox2{$\displaystyle#3\m@th$}%
    \sbox4{\@eqnnum}%
    \dimen@=.5\dimexpr\linewidth-\wd2\relax
    \ifcase
        \ifdim\wd0>\dimen@
          \z@
        \else
          \ifdim\wd4>\dimen@
            \z@
          \else 
            \@ne
          \fi 
        \fi
      \@latex@warning{Equation is too large}%
    \fi
    \noindent   
    \rlap{\copy0}%
    \rlap{\hbox to \linewidth{\hfill\copy2\hfill}}%
    \hbox to \linewidth{\hfill\copy4}%
    \hspace{0pt}
  \endgroup
  \ignorespaces 
}
\makeatother

   \usepackage[colorlinks = true,
            linkcolor = blue,
            urlcolor  = blue,
            citecolor = blue,
            anchorcolor = blue]{hyperref}
            
   \usepackage[parfill]{parskip}
   \usepackage{anysize}
   \marginsize{1in}{1in}{1in}{1in}
   \usepackage{thmtools}

\declaretheorem[name=Theorem,numberwithin=section]{thm}
\declaretheorem[name=Proposition,numberlike=thm]{prop}
\declaretheorem[name=Lemma,numberlike=thm]{lemma}
\declaretheorem[name=Corollary,numberlike=thm]{cor}

\declaretheorem[name=Definition,style=definition,qed=$\blacktriangle$,numberlike=thm]{defn}

\declaretheorem[name=Remark,style=definition,numberlike=thm]{rmk}

\newcommand{\KN}{\mathbin{\bigcirc\mspace{-16mu}\wedge\mspace{3mu}}}

\newcounter{commentCounter}
\newcommand\tab[1][1cm]{\hspace*{#1}}

\newcommand{\trr}{\operatorname{tr}}

\newcommand{\divv}{\operatorname{Div}}

\newcommand{\su}{\operatorname{SU}}
\newcommand{\uu}{\operatorname{U}}

\newcommand{\ii}{\operatorname{I}}

\DeclareFontFamily{U}{MnSymbolC}{}
\DeclareSymbolFont{MnSyC}{U}{MnSymbolC}{m}{n}
\DeclareFontShape{U}{MnSymbolC}{m}{n}{
    <-6>  MnSymbolC5
   <6-7>  MnSymbolC6
   <7-8>  MnSymbolC7
   <8-9>  MnSymbolC8
   <9-10> MnSymbolC9
  <10-12> MnSymbolC10
  <12->   MnSymbolC12}{}
\DeclareMathSymbol{\intprod}{\mathbin}{MnSyC}{'270}

\newcommand\myarray[1]{%
  \begingroup
  \renewcommand\arraystretch{1.33}
  \left\{ \begin{array}{@{}l@{}} #1 \end{array} \right.
  \endgroup}

\tikzcdset{
arrow style=tikz,
diagrams={>={Straight Barb[scale=0.8]}}
}

\makeatletter
\let\c@equation\c@thm
\makeatother
\numberwithin{equation}{section}

\begin{document}

\title{Betti numbers of nearly $G_2$ and nearly K\"{a}hler $6$-manifolds with Weyl curvature bounds}

\author{Anton Iliashenko \\ {\it Department of Pure Mathematics, University of Waterloo} \\ \tt{ailiashe@uwaterloo.ca}}

\maketitle

\begin{abstract}
In this paper we use the Weitzenböck formulas to get information about the Betti numbers of compact nearly $G_2$ and compact nearly K\"{a}hler $6$-manifolds. First, we establish estimates on two curvature-type self adjoint operators on particular spaces assuming bounds on the sectional curvature. Then using the Weitzenböck formulas on harmonic forms, we get results of the form: if certain lower bounds hold for these curvature operators then certain Betti numbers are zero. Finally, we combine both steps above to get sufficient conditions of vanishing of certain Betti numbers based on the bounds on the sectional curvature.
\end{abstract}

\tableofcontents

\section{Introduction}

\subsection{Motivation}
There is a long history of using Bochner-Weitzenböck technique to conclude vanishing results of Betti numbers of compact Riemannian manifolds assuming curvature bounds. In this paper we establish several results, particulary for compact nearly $G_2$ and compact nearly K\"{a}hler $6$-manifolds. We show that certain bounds on the sectional curvature imply vanishing of the second or the third Betti numbers.\\
Nearly $G_2$ and nearly K\"{a}hler $6$-manifolds are positive Einstein manifolds whose metric cones have $\mathrm{Spin}(7)$ and $\mathrm{G}_2$ holonomy, respectively. These, in turn, are useful from the physics perspective as they provide local models for the simplest type of interesting singularities. Hence, studying the topology of compact nearly $G_2$ and compact nearly K\"{a}hler $6$-manifolds might lead to new insights. For example, Betti numbers of such manifolds are related to Ricci flow. See~\cite{CWMYK} and~\cite{CWMYK2} for results relating Betti numbers and linear stability.
\subsection{Organization of the paper and main results}
Following Bourguignon-Karcher~\cite{BK}, we consider two curvature-type operators $\hat{R} \in \mathcal{S}^2(\Omega^2), \mathring{R} \in \mathcal{S}^2(\mathcal{S}^2)$ and the usual sectional curvature $\bar{R}$ coming from the Riemannian curvature. We prove the following theorems which give us bounds on these operators in terms of the bounds on the sectional curvature.\\
Here is a summary of the main results. Throughout, $[a\pm b]$ means $[a-b,a+b].$\\
First, we reprove the following result from~\cite{BK}.\\
\textbf{Theorem~\ref{hatthm1}}
Assume $\delta \leq \bar{R} \leq \Delta.$ Then the eigenvalues of $\hat{R}$ on $\Omega^2$ lie in the following interval: $$\left[-(\Delta+\delta) \pm \frac{4 \lfloor \frac{n}{2}\rfloor -1}{3} (\Delta - \delta)\right].$$
Then, for nearly $G_2$ or nearly K\"{a}hler $6$-manifolds, we improve the previous result on certain subspaces:
\textbf{Corollary~\ref{improvedaa}}
Assume $\delta \leq \bar{R} \leq \Delta.$ Moreover let $M$ be a nearly $G_2$ or a nearly K\"{a}hler $6$-manifold. Then on $\Omega^2_{14}$ or $\Omega^2_8$, respectively, the eigenvalues of $\hat{R}$ lie in the following interval:
$$\left[-(\Delta+\delta) \pm \frac{7}{3} (\Delta - \delta)\right].$$
Next, we again reprove a theorem from~\cite{BK} for $\mathring{R}$ in the general setting:\\
\textbf{Corollary~\ref{circthm1}}
Assume $\delta \leq \bar{R} \leq \Delta.$ Then all but one of the eigenvalues of $\mathring{R}$ on $\mathcal{S}^2$ lie in the following interval:
$$\left[\frac{1}{2}\bigg((\Delta+\delta) \pm (n-1)(\Delta-\delta)\bigg)\right],$$
and the other one lies in the interval:
$$\left[-(n-1)\Delta, -(n-1) \delta \right].$$ 
Following, we slightly improve the estimates for $\mathring{R}$ in the Einstein case:\\
\textbf{Theorem~\ref{nkmaincor}}
Suppose $M$ is Einstein with Einstein constant $k$. Assume $\delta \leq \bar{R} \leq \Delta.$   Then the eigenvalues of $\mathring{R}$ on $\mathcal{S}^2_0$ lie in the intersection of the following intervals:
$$[-k + n \delta, k - (n-2) \delta], [k - (n-2) \Delta,-k + n \Delta].$$
Next, for the nearly K\"{a}hler $6$-manifolds, we can talk about eigenvalues of $\mathring{R}$ on $\mathcal{S}^2_{+0} \subseteq \mathcal{S}^2$ (see Remarks~\ref{nkdecomofspaces} and~\ref{spacepreserving}). Hence, we are able to get a better estimate in this case:\\
\textbf{Theorem~\ref{nkimprov}}
Assume $\delta \leq \bar{R} \leq \Delta.$ Assume we are in the setting of a nearly K\"{a}hler $6$-manifold. Then the eigenvalues of $\mathring{R}$ on $\mathcal{S}^2_{+0}$ lie in the following interval:
$$\left[\frac{1}{2}\bigg((\Delta+\delta) \pm 3(\Delta-\delta)\bigg)\right]=[2\delta-\Delta, 2 \Delta- \delta].$$
Finally, we will see that again, on a nearly K\"{a}hler $6$-manifold, we have a specific relationship between $\hat{R}$ on $\Omega^2_8$ and $\mathring{R}$ on $\mathcal{S}^2_{+0}$, see Remark~\ref{asfdgdhgfhmg}. This allows us to get estimates for $\hat{R}$ on $\Omega^2_8$ in terms of the ones for $\mathring{R}$ on $\mathcal{S}^2_{+0}$ and vice versa. That is we can combine Corollary~\ref{improvedaa} and Theorem~\ref{nkimprov} to get the following two statements:\\
\textbf{Theorem~\ref{kjlakjkkj}}
Let $M$ be a nearly K\"{a}hler $6$-manifold. Let $\delta \leq \bar{R} \leq \Delta$. Then the eigenvalues of $\hat{R}$ on $\Omega^2_8$ lie in the intersection of the following intervals:
$$\left[-4+(\Delta+\delta) \pm 3(\Delta-\delta)\right],\left[-(\Delta+\delta) \pm \frac{7}{3}(\Delta-\delta)\right].$$
\textbf{Theorem~\ref{sdgfhgffgsdccc}}
Let $M$ be a nearly K\"{a}hler $6$-manifold. Let $\delta \leq \bar{R} \leq \Delta$. Then the eigenvalues of $\mathring{R}$ on $\mathcal{S}^2_{+0}$ lie in the intersection of the following intervals:
$$\left[\frac{1}{2}\bigg((\Delta+\delta) \pm 3(\Delta-\delta)\bigg)\right],\left[2+\frac{1}{2}\bigg(-(\Delta+\delta) \pm \frac{7}{3}(\Delta-\delta)\bigg)\right].$$
In the next section, where we introduce the Weitzenböck formulas (which relates the Laplacian $\Delta$ and the rough Laplacian (or Bochner Laplacian) $\nabla^* \nabla$ in terms of the Riemannian and Ricci curvatures) for $2$-forms and $3$-forms on nearly $G_2$ or nearly K\"{a}hler $6$-manifolds, the results are not new and can be found in the literature. However, we aim to keep the paper as self-contained as possible, so we include all the proofs, but we cite the results when appropriate.\\
The main idea is that for nearly K\"{a}hler and nearly $G_2$ manifolds, harmonic $2$-forms and harmonic $3$-forms are of a special algebraic type. In the case of $2$-forms this means that we need to consider the map $\hat{R}$ (or $\hat{W}$) only on certain subspaces of $\Omega^2$.\\
Moreover, when we apply the Weitzenböck formulas to harmonic forms to obtain sufficient conditions for certain Betti numbers to vanish in terms of lower bounds of $\hat{W}$ and $\mathring{W}$, where $W$ is the Weyl tensor (which is equivalent to some lower bounds on $\bar{R}$ and $\mathring{R}$), we get better estimates by considering the Weitzenböck formulas written in the intermediate forms. For example, consider~\eqref{eq:weiz2forms2}:
$$\Delta \beta =\nabla^{*} \nabla \beta + 8 \beta +  \hat{W}\beta, \text{ for any $\beta \in \Omega^2.$}$$
Assuming $\beta = h \diamond \omega$ (see Section~\ref{diamsection}) is harmonic for some $h\in \mathcal{S}^2_{+0}$, we can rewrite this as:
$$ 0 = \nabla^{*} \nabla \beta + (8h + 2\mathring{W} h) \diamond \omega = (\nabla^* \nabla h - 2h) + (8h + 2\mathring{W} h) \diamond \omega,$$
which is Proposition~\ref{thmnkffff1}. So, even though the last part is a well-known formula, we actually get better sufficient conditions for vanishing of $b_2$ in terms of the lower bound of $\mathring{W}$ by using the intermediate step above. Similar things happen in other cases as well.\\
We summarize the results we obtain in the folowing table:
\begin{center}
\begin{tabular}{ |p{3cm}||p{5cm}|p{5cm}|  }
 \hline
 \multicolumn{3}{|c|}{Sufficient conditions for vanishing of Betti numbers} \\
 \hline
 Manifold type& $b_2 = 0$ & $b_3 = 0$\\
 \hline\hline
Compact nearly $G_2$   & $\mathcal{S}^2(\Omega^2_{14}) \ni \hat{W} \geq - \frac{5 \tau^2_0}{8}$~(\ref{thmm2})    & $\mathcal{S}^2(\mathcal{S}^2_0) \ni \mathring{W}  \geq - \frac{3\tau^2_0}{8}$~(\ref{thmm1}), or\\[0.2cm]
&     & $\mathcal{S}^2(\Omega^2_{14}) \ni \hat{W} \geq -\frac{\tau^2_0}{4}$~(\ref{thmm1})   \\[0.2cm]
\hline
Compact nearly K\"{a}hler&$\mathcal{S}^2(\Omega^2_8) \ni \hat{W}  \geq -8$~(\ref{thmnk1}), or & $\mathcal{S}^2(\mathcal{S}^2_-) \ni \mathring{W}  \geq -\frac{9}{2}$~(\ref{thmnk2}), or\\[0.2cm]
of dim $6$&$\mathcal{S}^2(\mathcal{S}^2_{+0}) \ni \mathring{W} \geq -4$~(\ref{thmnk1}) & $\mathcal{S}^2(\Omega^2_8) \ni \hat{W} \geq -3$~(\ref{thmnk2})\\[0.2cm]
 \hline
\end{tabular}
\end{center}
We also use the fact that there are no parallel non-zero $2$-forms and no parallel non-zero traceless symmetric $2$-tensors. This is true because the holonomy is exactly $\operatorname{SO}(n).$
One can see this by observing that nearly $G_2$ and nearly K\"{a}hler manifolds admit a Killing spinor which implies that they are not locally reducible and nonsymmetric (the arguments can be found in~\cite{HB}), hence the result follows by Berger's classification. As corollaries, we obtain sufficient conditions for vanishing of the Betti numbers from inequalities $\delta \leq \bar{R} \leq \Delta$, which we again summarize in the following table: 
\begin{center}
\begin{tabular}{ |p{3cm}||p{5cm}|p{5cm}|  }
 \hline
 \multicolumn{3}{|c|}{Sufficient conditions for vanishing of Betti numbers} \\
 \hline
 Manifold type& $b_2 = 0$ & $b_3 = 0$\\
 \hline\hline
 Compact nearly $G_2$   & $-(\Delta + \delta) -\frac{7}{3} (\Delta - \delta) \geq -\frac{3\tau_0^2}{4}$~(\ref{aerstdfjhgggggg})    &  $\Delta \leq \frac{11\tau^2_0}{80}$~(\ref{sgfd}), or\\[0.2cm]
&     & $\delta \geq \frac{\tau^2_0}{112}$~(\ref{sgfd})   \\[0.2cm]
\hline
Compact nearly K\"{a}hler&$-(\Delta + \delta) -\frac{7}{3} (\Delta - \delta) \geq -10$~(\ref{cornkk1}), or & $\delta \geq \frac{1}{4}$~(\ref{thmnk3}), or\\[0.2cm]
of dim $6$&$(\Delta + \delta) -3 (\Delta - \delta) \geq -6$~(\ref{cornkk1}) & $\Delta \leq \frac{17}{8}$~(\ref{thmnk3})\\[0.2cm]
 \hline
\end{tabular}
\end{center}
Finally, for both nearly $G_2$ and nearly K\"{a}hler cases, we check our results on one example of a compact normal homogeneous manifold. The corollaries discussed above may not appear to be that useful for the known examples, as calculating the bounds for $\bar{R}$ is a harder process than just getting the bounds for $\hat{R}$ and $\mathring{R}$. However, our theorems are not limited to just these known examples, hence are interesting on their own.
\subsection{Acknowledgements}
I would like to thank my supervisor Spiro Karigiannis for helpful discussions and Uwe Semmelmann who communicated with my supervisor about some of the topics and sent me detailed and useful feedback on an earlier version of this article.

\subsection{Notation}\label{notsec}
Throughout this paper $(M^n,g)$ is a compact connected Riemannian manifold. We always work in a local orthonormal frame  $e_1, \dots ,e_n$. We often use the Riemannian metric (via the musical isomorphism) to identify vector fields and 1-forms. By $\mathcal{T}^k$ we denote $k$-tensors, by $\mathcal{S}^k$ the symmetric $k$-tensors, by $\mathcal{S}^2_0$ the traceless symmetric $2$-tensors, and by $\Omega^k$ the $k$-forms.\\
We also define the wedge product without any constants, meaning that for $\alpha, \beta \in \Omega^1$ we set $$\alpha \wedge \beta \coloneqq  \alpha \otimes \beta - \beta \otimes \alpha,$$ and extend to the higher order forms to preserve associativity.\\
The inner product on $k$-forms we define as follows. For $\alpha,\beta \in \Omega^k$: $$\langle \alpha, \beta \rangle = \frac{1}{k!} \alpha_{i_1 \dots i_k} \beta_{i_1 \dots i_k}.$$
We define the \textit{Kulkarni–Nomizu} product as follows. For $s,t \in \mathcal{T}^2$ we define $s \KN t \in \mathcal{T}^4$ to be:
\begin{equation}\label{eq:NKprod}
(s \KN t)_{ijkl} \coloneqq s_{il}t_{jk} +  s_{jk}t_{il} - s_{ik}t_{jl}- s_{jl}t_{ik}.
\end{equation}
We define the \textit{divergence operator} $\divv: \mathcal{S}^2 \rightarrow \Omega^1$ as $$(\divv h)_k = \nabla_i h_{ik}.$$
Next, for $\sigma \in \Omega^k$ and $h \in \mathcal{T}^2$, we define:
$$(h\diamond \sigma)_{i_1 \cdots i_k} \coloneqq h_{i_1 p}\sigma_{p i_2 \cdots i_k}+ h_{i_2 p}\sigma_{i_1 p i_3 \cdots i_k} + \cdots + h_{i_k p}\sigma_{i_1 \cdots i_{k-1} p}.$$
To simplify the interval notation, by $[a \pm b]$, we mean $[a-b,a+b],$ for $a,b \in \mathbb{R}, b>0.$\\
Finally, whenever we reffer to $\delta \leq \Delta$, these are any real numbers.

\begin{rmk}\label{rmkcurvature}
Let $\operatorname{Riem}$ be the Riemann curvature tensor. We will write $R_{ijkl}$ for $\operatorname{Riem}_{ijkl}$, $R_{ij}$ for $\operatorname{Ric}_{ij} \coloneqq R_{kijl} g^{kl}$ and $R \coloneqq \operatorname{Ric}_{ij} g^{ij} $ for the scalar curvature when there is no confusion.\\
Also we define the traceless Ricci tensor: $$\operatorname{Ric}^{0} \coloneqq \operatorname{Ric} - \frac{1}{n} R g.$$
Then on a general Riemannian manifold of dimension $n \geq 3$ we have the following orthogonal decomposition of $\operatorname{Riem}$ (see~\cite{AB}). Define
\begin{align*}
\text{traceless Ricci part: } E &\coloneqq  \frac{1}{n-2} \operatorname{Ric}^{0} \KN g,\\
\text{scalar part: } S &\coloneqq  \frac{R}{2n(n-1)} g \KN g,\\
\text{Weyl part: } W &\coloneqq \operatorname{Riem} - E -S.
\end{align*}
Then we have: $$\operatorname{Riem}=S+E+W.$$
Also, we say that $(M^n,g)$ is \textit{Einstein} with \textit{Einstein constant} $k$ if $\operatorname{Ric} = kg.$ In this case the scalar curvature is $R = nk$ and $\operatorname{Ric}^0 = 0$, thus $E = 0$. So, $\operatorname{Riem} = S + W,$ for an Einstein metric.\\
We also have, by construction, that $W_{kijl} g_{kl} = 0$.
\end{rmk}

\section{Curvature estimates}

Throughout this section, we let $(M,g)$ be a Riemannian manifold. First, we define a notion of a curvature tensor $A$. Then following Bourguignon-Karcher~\cite{BK} we introduce two self-adjoint operators $\hat{A}$ and $\mathring{A}$ and in Sections~\ref{estrhatggg} and~\ref{boundsmathringadsad} we obtain multiple results for bounds of $\hat{W}$ and $\mathring{W}$ in terms of bounds on the sectional curvature $\hat{R}$, where $W$ is the Weyl tensor. In particular, we strengthen some of the results from~\cite{BK} in the nearly $G_2$ and nearly K\"{a}hler of dimension $6$ settings.

\begin{defn}
We say an element $A \in \mathcal{T}^4$ is a \textit{curvature tensor}, if the following properties hold:
\begin{itemize}
\item $A_{ijkl}=-A_{jikl}=-A_{ijlk}=A_{klij}.$
\item $A_{ijkl}+A_{kijl}+A_{jkil} = 0\text{ (Bianchi identity)}.$
\end{itemize}
Let $\mathcal{R}$ be the set of curvature tensors. Note that $\mathcal{R}$ is a module over $C^{\infty}.$
\end{defn}

\begin{rmk}
If $s,t \in \mathcal{S}^2$, then $s \KN t \in \mathcal{R}$. This follows directly from the definition of $\KN$ in \eqref{eq:NKprod}. Hence, it follows from Remark~\ref{rmkcurvature} that $W$ is also a curvature tensor.
\end{rmk}

\begin{defn}\label{def1}
Let $A \in \mathcal{R}$. Following Bourguignon-Karcher \cite{BK}, we define $$\hat{A} \in \mathcal{S}^2(\Omega^2) \text{ as } (\hat{A}\beta)_{ij} = A_{ijkl}\beta_{kl}, \text{ for } \beta \in \Omega^2.$$
$$\mathring{A} \in \mathcal{S}^2(\mathcal{S}^2) \text{ by } (\mathring{A}h)_{ij} = A_{kilj}h_{kl}, \text{ by } h \in \mathcal{S}^2,$$
$$\bar{A} \text{ by } \bar{A}(X \wedge Y) = \frac{A(X,Y,Y,X)}{\|X \wedge Y\|^2}, \text{ for linearly independent } X,Y \in \Gamma(TM).$$
In particular, in an orthonormal frame: $\bar{A}(e_i \wedge e_j) = A_{ijji}$, for $i\neq j$ (with no sum over indices). We also call $\bar{A}$ the \textit{sectional curvature} of $A$, it is a smooth function on the space of $2$-planes on $M$.
\end{defn}
For the sake of completeness, we show that indeed, $\hat{A} \in \mathcal{S}^2(\Omega^2)$ and $\mathring{A} \in \mathcal{S}^2(\mathcal{S}^2)$.\\
Let $\beta, \gamma \in \Omega^2$. Then:
$$(\hat{A}\beta)_{ij} = A_{ijkl} \beta_{kl} = - A_{jikl} \beta_{kl}= -(\hat{A}\beta)_{ji}.$$
$$\langle \hat{A}\beta, \gamma\rangle = \frac{1}{2} A_{ijkl} \beta_{kl} \gamma_{ij} = \frac{1}{2} \beta_{kl}A_{klij} \gamma_{ij} = \langle\beta, \hat{A}\gamma\rangle. $$
Now, let $h,s \in \mathcal{S}^2.$ Then:
$$(\mathring{A}h)_{ij} = A_{kilj} h_{kl} = A_{ljki} h_{lk} = (\mathring{A}h)_{ji}.$$
$$\langle \mathring{A}h, s\rangle = A_{kilj} h_{kl} s_{ij} = h_{kl}A_{ljki} s_{ij} =   h_{lk} A_{jlik} s_{ji} = \langle h, \mathring{A} s\rangle.$$

\begin{rmk}\label{extendhat}
We can extend the map $\hat{A}$ to any $k$-form for $k > 2$ as follows: for $\beta \in \Omega^k$ we define $\hat{A} \beta \in \Omega^2 \bigotimes \Omega^{k-2}$ as $$(\hat{A} \beta)_{i_1 \dots i_{k-2}} \coloneqq A_{i_1 i_2 ab} \beta_{abi_3 \dots i_{k-2}},$$
that is, we just fix the last $k-2$ indices and think of $\beta$ as a $2$-form in the first two indices.\\
Also, note that $W_{kijl} g_{kl} = 0$ implies that for any $h \in \mathcal{S}^2$, we have $\mathring{W} h \in \mathcal{S}^2_0.$\\
From now on we will also use $R,\hat{R},\mathring{R},\bar{R}$ instead of $\operatorname{Riem},\hat{\operatorname{Riem}}$, etc., which should be clear from the context, and similarly for $W$.
\end{rmk}

\begin{lemma}\label{rmk1} The following identities hold:
\begin{itemize}
\item $g \bar{\KN} g = 2.$
\item $g \hat{\KN} g = -4 \operatorname{Id}.$
\item $g \mathring{\KN} g = 2 \operatorname{Id}$ on $\mathcal{S}^2_0$,
\end{itemize}
where by $g \bar{\KN} g$ we mean that we apply the $\bar{\phantom{}}$ operator to $g \KN g \in \mathcal{R}.$ Similarly, for the $g \hat{\KN} g$ and $g \mathring{\KN} g$.
\begin{proof}
For any $X,Y \in \Gamma(TM)$, we have:
\begin{align*}
(g \bar{\KN} g)(X \wedge Y) &= \frac{(g {\KN} g)(X,Y,Y,X)}{\|X\wedge Y\|^2}\\
&= \frac{2\|X\|^2\|Y\|^2 -2 \langle X,Y \rangle^2 }{\|X\wedge Y\|^2}\\
&=2.
\end{align*}
Next, let $\beta \in \Omega^2$ and $h \in \mathcal{S}^2_0$. Then in an orthonormal frame:
\begin{align*}
(\hat{(g \KN g)}\beta)_{ij} &= (g \KN g)_{ijkl} \beta_{kl}\\
= & 2 (g_{il}g_{jk} - g_{ik}g_{jl}) \beta_{kl}\\
= & 2 (\beta_{ji} - \beta_{ij}) \\
= & -4 \beta_{ij},
\end{align*}
and 
\begin{align*}
(\mathring{(g \KN g)}h )_{jl} &= (g \KN g)_{ijkl} h_{ik}\\
= & 2 (g_{il}g_{jk} - g_{ik}g_{jl}) h_{ik}\\
= & 2 ( h_{lj} - \trr(h) g_{jl}) \\
= & 2 h_{jl}.
\end{align*}
giving us the required results.
\end{proof}
\end{lemma}

In order to simplify the proofs of the following theorems we make the following definition:
\begin{defn}\label{defrest}
Assume that: $$\delta \leq \bar{R} \leq \Delta,$$
where $\delta, \Delta$ are any real constants. This means that for all $X,Y \in \Gamma(TM)$ with $\|X \wedge Y\|^2 = 1$, we have $\delta \leq \bar{R}(X \wedge Y) \leq \Delta$.\\
Define $$R_{0} \coloneqq R - \frac{\delta + \Delta}{4} g \KN g.$$
Now, by Lemma~\ref{rmk1}, $\bar{R}_{0} = \bar{R} - \frac{\delta + \Delta}{2}$, so that 
\begin{equation}\label{eq:RBARineq}
|\bar{R}_{0}| \leq \frac{\Delta - \delta}{2}.
\end{equation}
Note that $R_0 \in \mathcal{R}$, because both $R, g \KN g \in \mathcal{R}$.
\end{defn}

Next, we note that in the Einstein case, $\hat{W}$ and $\hat{R}$ differ by a constant multiple of the identity. The same holds for $\mathring{W}$ and $\mathring{R}$ on $\mathcal{S}^2_0$ (the constant is not the same though). 

\begin{lemma}\label{hatcirclemma}
Assume $M$ is Einstein with Einstein constant $k$. Then
\begin{align*}
\hat{W} &= \hat{R} + \frac{2k}{n-1} \operatorname{Id},\\
\mathring{W} &= \mathring{R} - \frac{k}{n-1} \operatorname{Id}, \text{ on $\mathcal{S}^2_0.$}
\end{align*}
\begin{proof}
By Remark~\ref{rmkcurvature}, $W = R - S$. Using Lemma~\ref{rmk1}, we have
$$\hat{S} =   \frac{R}{2n(n-1)} g \hat{\KN} g =  -\frac{nk}{2n(n-1)} 4 \operatorname{Id} =  -\frac{2k}{n-1} \operatorname{Id}.$$
Similarly on $\mathcal{S}^2_0$ we have
$$\mathring{S} =  \frac{R}{2n(n-1)} g \mathring{\KN} g = \frac{nk}{2n(n-1)} 2 \operatorname{Id} = \frac{k}{n-1} \operatorname{Id},$$
hence, the results follow.
\end{proof}
\end{lemma}

Finally, we have an observation about the a priori values of $\delta, \Delta$ in the Einstein case. 
\begin{rmk}\label{uselessfact}
Assume $(M^n,g)$ is Einstein with Einstein constant $k$. Let $\delta \leq \bar{R} \leq \Delta.$ Then: $$ \delta \leq \frac{k}{n-1} \leq \Delta.$$
\begin{proof} We compute
$$nk = R = \sum_{i=1}^n R_{ii} = \sum_{i,j=1}^n R_{ijji} = \sum_{i\neq j} \bar{R}(e_i \wedge e_j) \leq n(n-1) \Delta,$$
as when $i=j$, $R_{ijji} = 0$. So, $k \leq (n-1) \Delta$. The other inequality is done similarly.
\end{proof}
\end{rmk}

\subsection{Estimates for $\hat{R}$}\label{estrhatggg}
In this section we investiate what sectional curvature bounds tell us about the bounds of $\hat{R}$. Since in the Einstein case, $\hat{R}$ and $\hat{W}$ differ by a constant multiple of the identity map, one can use the result above to get bounds for $\hat{W}$.\\
First, we prove a lemma which gives us bounds for $R_0$ in terms of bounds of $\bar{R}$. Note that one can similarly obtain bounds for $R$ itself, but we do not need this.

\begin{lemma}\label{lemm1}
Assume $\delta \leq \bar{R} \leq \Delta.$ Let $X,Y,Z,W \in TM$ be unit length. Then $|R_0(X,Y,Z,W)| \leq \frac{2}{3} (\Delta - \delta).$
\begin{proof}
This result is Lemma 3.7 in~\cite{BK}, but we provide all the details.\\
Without loss of generality, assume $X \neq \pm W$ and $Y \neq \pm Z$. Otherwise, swap $Z$ and $W$. If even after swapping, that is not achieved, it means, $Z$ and $W$ are multiples of each other, so $R_0(X,Y,Z,W)=0.$\\
We claim that
\begin{equation}\label{eq:rhates1}
\begin{split}
6 R_0(X,Y,Z,W) =& R_0(X,Y+Z,Y+Z,W)-R_0(Y,X+Z,X+Z,W)\\
&-R_0(X,Y-Z,Y-Z,W)+R_0(Y,X-Z,X-Z,W).
\end{split}
\end{equation}
Expanding the RHS we get:
\begin{align*}
&R_0(X,Y,Y,W)+ R_0(X,Z,Z,W)+R_0(X,Y,Z,W)+R_0(X,Z,Y,W)\\
&-R_0(Y,X,X,W)- R_0(Y,Z,Z,W)-R_0(Y,X,Z,W)-R_0(Y,Z,X,W)\\
&-R_0(X,Y,Y,W)- R_0(X,Z,Z,W)+R_0(X,Y,Z,W)+R_0(X,Z,Y,W)\\
&+R_0(Y,X,X,W)+ R_0(Y,Z,Z,W)-R_0(Y,X,Z,W)-R_0(Y,Z,X,W)\\
=& 4 R_0(X,Y,Z,W) - 2 \Big(R_0(Z,X,Y,W) +R_0(Y,Z,X,W)\Big)\\
=& 6 R_0(X,Y,Z,W),
\end{align*}
as claimed. Now, consider one of the terms $R_0(X,Y+Z,Y+Z,W)$:
\begin{align*}
R_0(X,Y+Z,Y+Z,W) =& \frac{1}{4} \Big( R_0(X+W,Y+Z,Y+Z,X+W) - R_0(X-W,Y+Z,Y+Z,X-W) \Big)\\
 =& \frac{\|Y+Z\|^2}{4} \Big( \|X+W\|^2  R_0(\frac{X+W}{\|X+W\|},\frac{Y+Z}{\|Y+Z\|},\frac{Y+Z}{\|Y+Z\|},\frac{X+W}{\|X+W\|})\\
 &- \|X-W\|^2 R_0(\frac{X-W}{\|X-W\|},\frac{Y+Z}{\|Y+Z\|},\frac{Y+Z}{\|Y+Z\|},\frac{X-W}{\|X-W\|}) \Big).
\end{align*}
Now, note that for unit length vectors $S,T$ we have:
$$|R_0(S,T,T,S)| = |\bar{R}_{0}(S,T)| (\|S\|^2\|T\|^2 - \langle S,T\rangle^2) \leq |\bar{R}_{0}(S \wedge T)|. $$
Thus:
\begin{align*}
|R_0(X,Y+Z,Y+Z,W)| & \leq \frac{\|Y+Z\|^2}{4} \Big( \|X+W\|^2 |\bar{R}_{0}(\frac{X+W}{\|X+W\|} \wedge \frac{Y+Z}{\|Y+Z\|})|\\
&+ \|X-W\|^2 |\bar{R}_{0}(\frac{X-W}{\|X-W\|} \wedge \frac{Y+Z}{\|Y+Z\|})| \Big)\\
& \leq \frac{\|Y+Z\|^2}{4} (\|X+W\|^2 + \|X-W\|^2) \frac{\Delta - \delta}{2}\text{ (by \eqref{eq:RBARineq})}\\
& = \|Y+Z\|^2 \frac{\Delta - \delta}{2}.
\end{align*}
Hence, applying the same inequalities for the other terms, equation~\eqref{eq:rhates1} becomes:
\begin{align*}
6|R_0(X,Y,Z,W)| & \leq (\|Y+Z\|^2 + \|X+Z\|^2+\|Y-Z\|^2+ \|X-Z\|^2)\frac{\Delta - \delta}{2}\\
&=4(\Delta - \delta),
\end{align*}
which yields the desired result.
\end{proof}
\end{lemma}

We are ready to get to the main theorem of this section. The first part applies to any manifold, however on certain subspaces of manifolds with $G_2$ or $SU(3)$-structure, we can improve the result.
\begin{thm}\label{hatthm1}
Assume $\delta \leq \bar{R} \leq \Delta.$ Then the eigenvalues of $\hat{R}$ lie in the following interval: $$\left[-(\Delta+\delta) \pm \frac{4 \lfloor \frac{n}{2}\rfloor -1}{3} (\Delta - \delta)\right].$$
\begin{proof}
Assume $\hat{r}$ is an eigenvalue of $\hat{R}$ with $0 \neq \beta \in \Omega^2$ the corresponding unit eigenvector. Note that $\hat{R}_0 = \hat{R} + (\Delta+\delta) \operatorname{Id}$, by Remark~\ref{rmk1} and Definition~\ref{defrest}. So, $\beta$ is also an eigenvector for $\hat{R}_0$ with the eigenvalue $\hat{r}_0 = \hat{r} + (\delta + \Delta)$.\\
Assume $\beta$ is of rank $2p$, so there exists an orthonormal basis $\{e_1, \dots , e_n\}$ such that $\displaystyle \beta = \sum_{i=1}^p \beta_i e_i \wedge e_{\bar{i}}$, where $\bar{i} = i + p$. Then we have
\begin{align*}
\hat{r}_0 \beta_j &= (\hat{r}_0 \beta)_{j \bar{j}} \\
&= (\hat{R}_0 \beta)_{j \bar{j}}\nonumber \\
&= \sum_{i=1}^p \beta_i (\hat{R}_0 (e_i \wedge e_{\bar{i}}) )_{j \bar{j}} \\
&= \sum_{i=1}^p \beta_i (R_0)_{pl j \bar{j}} (e_i \wedge e_{\bar{i}}) _{pl} \\
&= \sum_{i=1}^p \beta_i (R_0)_{pl j \bar{j}} (\delta_{ip} \delta_{\bar{i} l}- \delta_{il} \delta_{\bar{i} p}) \\
&= 2 \sum_{i=1}^p \beta_i (R_0)_{i \bar{i} j \bar{j}}.
\end{align*}
Now, take $|\beta_j| \neq 0$ maximal to obtain from the above that
\begin{align}
|\hat{r}_0| &\leq 2 \sum_{i=1}^p \Big|\frac{\beta_i}{\beta_j}\Big| |(R_0)_{i \bar{i} j \bar{j}}|\nonumber\\
&= 2 \sum_{i \neq j} \big|\frac{\beta_i}{\beta_j}\big| |(R_0)_{i \bar{i} j \bar{j}}| + 2 |(\bar{R}_0)(e_j \wedge e_{\bar{j}})|\label{randomequation12}\\
&\leq 2 (p-1) \frac{2}{3} (\Delta - \delta) + 2 \frac{\Delta - \delta}{2} \text{ (by Lemma~\ref{lemm1} and~\eqref{eq:RBARineq})}\nonumber\\
&= \frac{4p-1}{3} (\Delta - \delta)\nonumber\\
&\leq \frac{4 \lfloor \frac{n}{2}\rfloor -1}{3} (\Delta - \delta).\nonumber
\end{align}
Recalling that $\hat{r}_0 = \hat{r} + (\delta + \Delta)$, we get the required result.
\end{proof}
\end{thm}

Adding onto the work of Bourguignon-Karcher \cite{BK}, the previous theorem can be improved for nearly $G_2$ or nearly K\"{a}hler $6$-manifolds on certain subspaces.
\begin{cor}\label{improvedaa}
In the nearly $G_2$ case on $\Omega^2_{14}$ or in the nearly K\"{a}hler case on $\Omega^2_8$ the eigenvalues of $\hat{R}$ lie in the following interval:
$$\left[-(\Delta+\delta) \pm \frac{7}{3} (\Delta - \delta)\right].$$
See Sections~\ref{g2prelim} and~\ref{nkprelim} for the descriptions of these manifolds and subspaces. Note that just the presence of a $G_2$ or an $SU(3)$ structure is not enough, as we need $\hat{R}$ to preserve those subspaces.\\
Note that the previous Theorem~\ref{hatthm1}, only would have given us $\frac{11}{3}$ instead of $\frac{7}{3}$.
\begin{proof}
For both the $G_2$-structure case on $\Omega^2_{14}$ or for the $SU(3)$-structure case on $\Omega^2_8$, if we assume $\beta$ is of rank $2p = 2,4,6$, then there exist canonical forms $\displaystyle \beta = \sum_{i=1}^p \beta_i e_i \wedge e_{\bar{i}}$, where $\bar{i} = i + p$, such that $\displaystyle \sum_{i=1}^k \beta_i = 0$, for some orthonormal basis $\{e_1, \dots e_n\}$ (in the case of $G_2$-structures, see~\cite{Sp3}, and in the case of $SU(3)$-structure, this follows because $\Lambda^2_{8} \cong \mathfrak{su}(3)$). Taking $|\beta_j| \neq 0$ maximal forces the other $\beta_i$'s, of which there are at most two, to be of the same sign, meaning that $\displaystyle |\beta_j| = \sum_{i \neq j} |\beta_i|.$ Thus, continuing from~\eqref{randomequation12}, we can improve the previous estimate to:
\begin{align*}
|\hat{r}_0| &\leq 2 \sum_{i \neq j} \big|\frac{\beta_i}{\beta_j}\big| |(R_0)_{i \bar{i} j \bar{j}}| + 2 |(\bar{R}_0)(e_j \wedge e_{\bar{j}})|\\
&\leq 2 \frac{2}{3} (\Delta - \delta) + 2 \frac{\Delta - \delta}{2}\text{ (by Lemma~\ref{lemm1} and~\eqref{eq:RBARineq})}\\
&= \frac{7}{3} (\Delta - \delta).
\end{align*}
which is enough to conclude the result.
\end{proof}
\end{cor}

We will see that in the nearly K\"{a}hler case, the operators $\hat{W}$ and $\mathring{W}$ are closely related on certain subspaces. See Remark~\ref{asfdgdhgfhmg}. Hence, we summarize the estimates for $\hat{R}$ on $\Omega^2_8$ in the following Corollary:
\begin{cor}\label{kjlakjkkj}
Let $M$ be a nearly K\"{a}hler $6$-manifold. Let $\delta \leq \bar{R} \leq \Delta$. Then the eigenvalues of $\hat{R}$ on $\Omega^2_8$ lie in the intersection of the following intervals:
$$\left[-4+(\Delta+\delta) \pm 3(\Delta-\delta)\right],\left[-(\Delta+\delta) \pm \frac{7}{3}(\Delta-\delta)\right].$$
\begin{proof}
This follows from Remark~\ref{asfdgdhgfhmg}.
\end{proof}
\end{cor}

\subsection{Estimates for $\mathring{R}$}\label{boundsmathringadsad}

Note that when $M$ is Einstein, $\mathring{R},\mathring{W}$ preserve $\mathcal{S}^2_0$. This is because $\mathring{W} h \in \mathcal{S}^2_0$ for any $h \in \mathcal{S}$, by the properties of the Weyl tensor, and since  $\mathring{R}$ and  $\mathring{W}$ differ by a constant on $\mathcal{S}^2_0$, we get the required observation.

First, we prove a therorem that gives us bounds for $\mathring{R}$ on $\mathcal{S}^2$ in terms of bounds of $\bar{R}$. Next, we assume that $M$ is Einstein which allows us to improve the result on $\mathcal{S}^2_0$.
\begin{thm}\label{circthm1}
Assume $\delta \leq \bar{R} \leq \Delta.$ Then all but one of the eigenvalues of $\mathring{R}$ on $\mathcal{S}^2$ lie in the following interval:
$$\left[\frac{1}{2}\bigg((\Delta+\delta) \pm (n-1)(\Delta-\delta)\bigg)\right],$$
and the other one lies in the interval:
$$\left[-(n-1)\Delta, -(n-1) \delta \right].$$ 
\begin{proof}
On $\mathcal{S}^2_0$, $\mathring{R} = \mathring{R}_0 + \frac{\delta + \Delta}{4} g \mathring{\KN} g = \mathring{R}_0 +\frac{\delta + \Delta}{2} \operatorname{Id}$, by Lemma~\ref{rmk1} and Definition~\ref{defrest}.\\
Recall that by Definition~\ref{defrest} we have that $\mathring{R} = \mathring{R}_0 + \frac{\delta + \Delta}{4} g \mathring{\KN} g.$\\
First, we show that $|\mathring{R}_0 |\leq \frac{n-1}{2} (\Delta -\delta)$:
Let $0 \neq h \in \mathcal{S}^2_0$ be a unit eigenvector of  $\mathring{R}_{0}$ with the eigenvalue $\mathring{r}_0$. Assume $h$ is of rank $p$ for some $1 \leq p \leq n$. Then there exists an orthonormal basis $\{e_1, \dots, e_n\}$ such that $\displaystyle h = \sum_{i=1}^{p} h_i e_i \otimes e_i$. Thus:
\begin{align*}
\mathring{r}_{0} h_j &= (\mathring{R}_{0} h)_{jj}\\
&= (R_0)_{mjlj} h_{ml}\\
&= (R_0)_{mjlj} h_m \delta_{ml}\\
&= \sum_{m} (R_0)_{mjmj} h_m. \\
\end{align*}
Take $|h_j| \neq 0$ maximal. We obtain from the above that:
\begin{align*}
|\mathring{r}_{0}| &\leq \sum_{m} \Big|\frac{h_m}{h_j} \Big| |(R_0)_{mjmj}|\\
&\leq (p-1)  |\bar{R_{0}}|\\
&\leq (n-1)  \frac{\Delta - \delta}{2},\\
\end{align*}
yielding the required result.\\
Next, we investigate the eigenvalues of $\mathring{R}- \mathring{R}_0 = \frac{\delta + \Delta}{4} g \mathring{\KN} g.$ It is easy to check that $(g \mathring{\KN} g) g = 2(1-n)g$, and we know that $g \mathring{\KN} g  = 2 \operatorname{Id}$ on  $\mathcal{S}^2_0,$ by Lemma~\ref{rmk1}. \\Hence, the result follows from the Weyl's inequality for eigenvalues applied to $\mathring{R} = \mathring{R}_0+ ( \mathring{R}- \mathring{R}_0)$.
\end{proof}
\end{thm}

\begin{thm}\label{nkmaincor}
Suppose $M$ is Einstein with Einstein constant $k$. Assume $\delta \leq \bar{R} \leq \Delta.$   Then the eigenvalues of $\mathring{R}$ on $\mathcal{S}^2_0$ lie in the intersection of the following intervals:
$$[-k + n \delta, k - (n-2) \delta], [k - (n-2) \Delta,-k + n \Delta].$$
\begin{proof}
For simplicity, introduce $R^{'} \coloneqq R - \frac{\delta}{2} g \KN g$. Then $\bar{R}^{'} = \bar{R} - \delta$ and $\mathring{R}^{'} = \mathring{R} - \delta \operatorname{Id}$, by Remark~\ref{rmk1}. Hence, $\bar{R}^{'} \geq 0.$ By $\mathring{R}{'}$ we will mean $\mathring{\phantom{}}$ applied to $R{'}$, and similarly for $\bar{R}{'}.$\\
Let $0 \neq h \in S^2_0$ be a unit eigenvector of  $\mathring{R}$ with the eigenvalue $\mathring{r}$. Note that $h$ is also an eigenvector of  $\mathring{R}^{'}$ with the eigenvalue $\mathring{r}' = \mathring{r} - \delta$.  Assume $h$ is of rank $p$ for some $1 \leq p \leq n$. Then there exists an orthonormal basis $\{e_1, \dots, e_n\}$ such that $\displaystyle h = \sum_{i=1}^{p} h_i e_i \otimes e_i$. Thus:
\begin{align*}
\mathring{r}{'} h_j &= (\mathring{R}{'} h)_{jj}\\
&= (R{'})_{mjlj} h_{ml}\\
&= (R{'})_{mjlj} h_m \delta_{ml}\\
&= \sum_{m=1}^{p} (R{'})_{mjmj} h_m\\
&= -\sum_{m=1}^{p} (\bar{R}{'})_{mj} h_m. \\
\end{align*}
Take $|h_j| \neq 0$ maximal. By replacing $h$ by $-h$, if necessary, assume that $h_j > 0$. Note that now for all $m$, $-1 \leq \frac{h_m}{h_j} \leq 1$. Then since $\bar{R}^{'} \geq 0$, we have:
\begin{align*}
-\mathring{r}' &= \sum_{m=1}^{p} \frac{h_m}{h_j} \bar{R}{'}(e_m \wedge e_j)\\
& \leq \sum_{m=1}^{n} \bar{R}{'}(e_m \wedge e_j)\\
& = \sum_{m=1}^{n} (\bar{R} - \delta)(e_m \wedge e_j) \\
& = \sum_{m=1}^{n} \bar{R}(e_m \wedge e_j)  - (n-1) \delta.\\
\end{align*}
Finally, note that $\displaystyle \sum_{m=1}^{n} (\bar{R})_{mj} =\sum_{m=1}^{n} R_{mjjm} = R_{jj} = k g_{jj} = k$ (where the $j$ was fixed.)
Hence, $$-(\mathring{r} - \delta) = -\mathring{r}' \leq k - (n-1) \delta,$$
which gives the required $$\mathring{r} \geq -k + n \delta.$$
However, (this was not present in the Bourguignon-Karcher paper) since $-1 \leq \frac{h_m}{h_j}$, we can also do the following:
\begin{align*}
-\mathring{r}' &= \sum_{m=1}^{p} \frac{h_m}{h_j} \bar{R}{'}(e_m \wedge e_j)\\
&\geq -\sum_{m=1}^{p} \bar{R}{'}(e_m \wedge e_j)\\
&\geq -\sum_{m=1}^{n} \bar{R}{'}(e_m \wedge e_j)\\
&=-(k - (n-1)\delta).
\end{align*}
Hence, we also get
$$-(\mathring{r} - \delta) = -\mathring{r}{'} \geq -k + (n-1) \delta,$$
which is just
$$ \mathring{r} \leq k - (n-2) \delta.$$
Thus, we have
$$-k + n \delta \leq \mathring{r} \leq k - (n-2) \delta.$$
The other inequality $$k - (n-2) \Delta \leq \mathring{r} \leq -k + n \Delta$$ is proven in the similar way by introducing $R{''} \coloneqq R - \frac{\Delta}{2} g \KN g$, so $\bar{R}{''}\leq 0.$
\end{proof}
\end{thm}

\begin{rmk}
In~\cite{BK}, the authors proved the estimate $$-k + n \delta \leq \mathring{R} \leq -k + n \Delta, \tab \text{on $\mathcal{S}^2_0.$}$$
which is weaker than Theorem~\ref{nkmaincor}.
\end{rmk}

\begin{prop}\label{nkimprov}
Assume $\delta \leq \bar{R} \leq \Delta.$ Assume we are in the setting of a nearly K\"{a}hler $6$-manifold. Then by Remark~\ref{spacepreserving}, $\mathring{R}$ preserves $\mathcal{S}^2_{0+}$ and we claim that the eigenvalues of $\mathring{R}$ on $\mathcal{S}^2_{+0}$ lie in the following interval:
$$\left[\frac{1}{2}\bigg((\Delta+\delta) \pm 3(\Delta-\delta)\bigg)\right]=[2\delta-\Delta, 2 \Delta- \delta].$$
\begin{proof}
As in Corollary~\ref{improvedaa}, we use the fact that for an element $ h \in \mathcal{S}^2_{0+}$ (which by Proposition~\ref{h2forms} is isomorphic to $\Lambda^2_{8} \cong \mathfrak{su}(3)$) we can find a canonical form $h = \sum_{i=1}^{6} h_i e_i \otimes e_i$ with $e_1, \dots, e_6$ an orthonormal frame and $h_1=h_2, h_3=h_4, h_5=h_6$ with $h_1+h_2+h_3 = 0.$\\
We proceed in the same way as in the proof of Theorem~\ref{circthm1}. Let $0 \neq h \in \mathcal{S}^2_{0+}$ be a unit eigenvector of $\mathring{R_0}$ and let $\mathring{r}_0$ be the corresponding eigenvalue. Put $h$ in the canonical form as above. By replacing $h$ by $-h$ and by swapping $h_i$'s if necessary, we can assume $|h_1| > 0$ is maximal, $h_1 = h_2 > 0$ and as before, $h_3=h_4, h_5 = h_6$. This forces $h_3,h_5$ to be non-positive with
\begin{equation}\label{eq:afdsddddd}
|h_3|+|h_5| = |h_4|+|h_6| = h_1.
\end{equation}
As before, we get:
$$\mathring{r}_{0} h_1 =  \sum_{m} (R_0)_{m1m1} h_m.$$
Dividing through by $h_1 > 0$, and using~\eqref{eq:RBARineq} with~\eqref{eq:afdsddddd} we get:
\begin{align*}
|\mathring{r}_{0}| &\leq \sum_{m} \frac{h_m}{h_1} |(R_0)_{m1m1}|\\
&= |(R_0)_{2121}| + \frac{|h_3|}{h_1} |(R_0)_{3131}|+ \frac{|h_4|}{h_1} |(R_0)_{4141}|+ \frac{|h_5|}{h_1} |(R_0)_{5151}|+ \frac{|h_6|}{h_1} |(R_0)_{6161}|\\
&\leq 3 \operatorname{max}|\bar{R}_0|\\
&\leq \frac{3(\Delta - \delta)}{2},
\end{align*}
which along with the same details as in Theorem~\ref{circthm1} concludes the proof.\qedhere
\end{proof}
\end{prop}
\begin{rmk}\label{dfgdsfgas}
In the Einstein setting, Theorem~\ref{circthm1} tells us that the eigenvalues of $\mathring{R}$ on $\mathcal{S}^2_0$ lie in $$\left[\frac{1}{2}\big((\Delta+\delta) \pm (n-1)(\Delta-\delta)\big)\right],$$
which is a weaker result than the one in Theorem~\ref{nkmaincor}. It is immediately clear that on a nearly K\"{a}hler $6$-manifold, on $\mathcal{S}^2_{0+}$, the interval from Theorem~\ref{nkimprov} is a better result than Theorem~\ref{circthm1}.\\
One can also show that Theorem~\ref{nkimprov} is also stronger than Theorem~\ref{nkmaincor}. For example, let us show that $2\delta - \Delta \geq -5 + 6\delta.$
So, we need $5 \geq 4\delta + \Delta$. This is clearly true, as we can pick an orthonormal frame where $R_{1212} = \Delta$. Then we know that $R_{1i1i} \geq \delta$ for $i = 3,4,5$ and $\sum_{i=2}^5 R_{1i1i}  = 5$, which is the Einstein constant for a nearly K\"{a}hler $6$-manifold. Hence the result follows. All the other inequalities are similar.
\end{rmk}
\begin{rmk}\label{asfdgdhgfhmg}
In the proof of Theorem~\ref{thmnkffff1}, we will show that on a nearly K\"{a}hler $6$-manifold, for $\beta \in \Omega^2_8$, which must equal $h \diamond \omega$ for some unique $h \in \mathcal{S}^2_{0+}$, we have that $\hat{W}\beta = (2 \mathring Wh) \diamond \omega.$ Hence, if $\beta = h \diamond \omega$ is an eigenvector of $\hat{W}$ with the eigenvalue $\lambda$, then $h$ is an eigenvector $\mathring{W}$ with the eigenvalue $\frac{\lambda}{2}.$ This clearly means that $\operatorname{range}(\hat{W}) =2\operatorname{range}(\mathring{W}),$ where by $\operatorname{range}$ of a self-adjoint operator we mean the closed interval from the smallest eigenvalue to the largest one.\\
By Lemma~\ref{hatcirclemma} we have that $\hat{W} = \hat{R} + 2 \operatorname{Id},\mathring{W} = \mathring{R} - \operatorname{Id}$ on $\mathcal{S}^2_0.$\\
Hence, assume that $\delta \leq \bar{R} \leq \Delta$. Then by Corollary~\ref{improvedaa}, the range of $\hat{W}$ on $\Omega^2_8$ lies in
\begin{equation}\label{eq:interval1}
\left[2-(\Delta+\delta) \pm \frac{7}{3} (\Delta - \delta)\right].
\end{equation}
Similarly, by Theorem~\ref{nkimprov}, the range of $\mathring{W}$ on $\mathcal{S}^2_{0+}$ lies in
\begin{equation}\label{eq:interval2}
\left[-1+\frac{1}{2}\bigg((\Delta+\delta) \pm 3(\Delta-\delta)\bigg)\right].
\end{equation}
However, since $\operatorname{range}(\hat{W}) =2\operatorname{range}(\mathring{W}),$ $\hat{W}$ on $\Omega^2_8$ also lies in
$$\left[-2+(\Delta+\delta) \pm 3(\Delta-\delta)\right],$$
which is clearly not the same interval as in~\eqref{eq:interval1}. Since we cannot say that one of the intervals is always better than the other one, we will use them both by taking their intersection. Similarly, we can also obtain a second interval for $\mathring{W}$ on $\mathcal{S}^2_{0+}$.
\end{rmk}
We summarize the estimates in the case of a nearly K\"{a}hler $6$-manifold for $\mathring{R}$ on $\mathcal{S}^2_0.$
\begin{cor}\label{sdgfhgffgsdccc}
Let $M$ be a nearly K\"{a}hler $6$-manifold. Let $\delta \leq \bar{R} \leq \Delta$. Then the eigenvalues of $\mathring{R}$ on $\mathcal{S}^2_{+0}$ lie in the intersection of the following intervals:
$$\left[\frac{1}{2}\bigg((\Delta+\delta) \pm 3(\Delta-\delta)\bigg)\right],\left[2+\frac{1}{2}\bigg(-(\Delta+\delta) \pm \frac{7}{3}(\Delta-\delta)\bigg)\right].$$
\begin{proof}
This follows from Remark~\ref{asfdgdhgfhmg}.
\end{proof}
\end{cor}

\section{General Weitzenböck formulas}
In this section we rederive the well-known general Weitzenböck formula and then simplify it in the case that the manifold is Einstein. More information can be found in~\cite{USGW},~\cite{PPMW},~\cite{SUU}.\\
Let $(M^n,g)$ be a Riemannian manifold.
For $\alpha \in \Omega^k$ and $T \in \mathcal{T}^{k}$ we have:
$$(d\alpha)_{i_1 \dots i_{k+1}} = \sum_{j=1}^{k+1} (-1)^{j-1} \nabla_{i_j} \alpha_{i_1 \dots \hat{i_j} \dots i_{k+1}},$$
$$(d^*\alpha)_{i_1 \dots i_{k-1}} = -\nabla_p \alpha_{p i_1 \dots i_{k-1}},$$
$$(\nabla^* \nabla \alpha)_{i_1 \dots i_{k-1}} = - \nabla_p\nabla_p  \alpha_{i_1 \dots i_{k-1}},$$
$$(\nabla^* T)_{i_2 \dots i_k} = - \nabla_{p} T_{p i_2 \dots i_k}.$$

\begin{prop}\label{weiz} \textbf{General Weitzenböck formula}\\
For $\alpha \in \Omega^k$ we have:
\begin{align*}
(\Delta \alpha)_{i_1 \dots \i_k} =  (\nabla^* \nabla \alpha)_{i_1 \dots \i_k} +& \sum_{j=1}^{k} \alpha_{i_1 \dots u \dots i_k} R_{i_j u} \text{    ($u$ is at $j^{\text{th}}$ position)}\\
+& \sum_{1\leq l < j \leq k}  \alpha_{i_1 \dots u \dots p \dots i_k} R_{i_j i_l p u} \text{    ($u$ and $p$ are at $l^{\text{th}}$ and $j^{\text{th}}$ positions respectively)}
\end{align*}
\end{prop}

From Proposition~\ref{weiz} we obtain the \textbf{Weitzenböck formula for $2$-forms}:\\
Let $\beta \in \Omega^2$. Then:
\begin{equation}\label{eq:weiz2forms}
(\Delta \beta)_{ab} =(\nabla^{*} \nabla \beta)_{ab} + R_{ap} \beta_{pb} +  R_{bp} \beta_{ap} + R_{abpq} \beta_{pq}.
\end{equation}

\begin{cor}\label{weiz2simpl}
Assume $M$ is Einstein with Einstein constant $k$. Then the Weitzenböck formula for $2$-forms simplifies to:
$$\Delta \beta =\nabla^{*} \nabla \beta + 2k \frac{n-2}{n-1} \beta + \hat{W} \beta,$$
where the $\hat{W}$ notation is defined in Section~\ref{def1}.
\begin{proof}
In the Einstein case we have $Ric = kg.$ Hence, each of the Ricci terms in~\eqref{eq:weiz2forms} is equal to $k \beta_{ab}$.\\
The last term is: $$R_{abpq} \beta_{pq} = (\hat{R}\beta)_{ab} =  ((\hat{W}-\frac{2k}{n-1} \operatorname{Id})\beta)_{ab} =(\hat{W}\beta)_{ab} - \frac{2k}{n-1}\beta_{ab},$$
by Lemma~\ref{hatcirclemma}.
Thus, putting everything together, we get: $$\Delta \beta = \nabla^{*} \nabla \beta + 2(k\beta)+(\hat{W}\beta -\frac{2k}{n-1}\beta) = \nabla^{*} \nabla \beta +  2k \frac{n-2}{n-1} \beta +  \hat{W}\beta,$$
as required.
\end{proof}
\end{cor}

From Proposition~\ref{weiz} we obtain the \textbf{Weitzenböck formula for $3$-forms}:\\
Let $\beta \in \Omega^3$. Then:
\begin{equation}\label{eq:weiz3forms}
(\Delta \beta)_{abc} =(\nabla^{*} \nabla \beta)_{abc} + R_{au} \beta_{ubc} +  R_{bu} \beta_{auc} +  R_{cu} \beta_{abu} +  R_{abpu} \beta_{puc}+ R_{acpu} \beta_{pbu}+ R_{bcpu} \beta_{apu}.
\end{equation}

\begin{cor}\label{weiz3simpl}
Assume $M$ is Einstein with Einstein constant $k$. Then the Weitzenböck formula for $3$-forms simplifies to:
$$(\Delta \beta)_{abc} =(\nabla^{*} \nabla \beta)_{abc} + 3k\frac{n-3}{n-1} \beta_{abc} +  W_{abpu} \beta_{puc}+ W_{acpu} \beta_{pbu}+ W_{bcpu} \beta_{apu}.$$
Note that the last three terms can be written as $(\hat{W} \beta)_{abc}+(\hat{W} \beta)_{cab}+(\hat{W} \beta)_{bca},$ by using notation of Remark~\ref{extendhat}.
\begin{proof}
In the Einstein case we have $Ric = kg.$ Hence, each of the Ricci terms in~\eqref{eq:weiz3forms} is equal to $k \beta_{abc}$.\\
Now, consider the term $R_{abpu} \beta_{puc}$ of~\eqref{eq:weiz3forms}. This is $(\hat{R} \beta)_{abc}$, which is equal to $((\hat{W}-\frac{2k}{n-1}\operatorname{Id}) \beta)_{abc} = W_{abpu} \beta_{puc} - \frac{2k}{n-1}\beta_{abc},$ by Lemma~\ref{hatcirclemma}. Similarly, we can do the same for the other terms to get the required result:
\begin{align*}
(\Delta \beta)_{abc} &= (\nabla^{*} \nabla \beta)_{abc} + 3(k \beta_{abc})+W_{abpu} \beta_{puc}+ W_{acpu} \beta_{pbu}+ W_{bcpu} \beta_{apu} - 3(\frac{2k}{n-1}\beta_{abc})\\
&= (\nabla^{*} \nabla \beta)_{abc} + 3k\frac{n-3}{n-1} \beta_{abc} +  W_{abpu} \beta_{puc}+ W_{acpu} \beta_{pbu}+ W_{bcpu} \beta_{apu}.\qedhere
\end{align*}
\end{proof}
\end{cor}

\section{Nearly $G_2$ manifolds}

 First, in Section~\ref{g2prelim} we recall some facts about $G_2$ structures and nearly $G_2$ manifolds. In Section~\ref{ng2jkcurg} we observe some properties about the curvature of nearly $G_2$ manifolds. Finally, in Section~\ref{ng2weindgfgf} we simplify the Weizenbock formulas for harmonic $2$ and $3$-forms and using the assumption of compactness of our manifolds, we get the necessary conditions of vanishing of $b_2$ and $b_3$ in terms of bounds on $\bar{R}, \mathring{R},$ and $\bar{R}$.

\subsection{Preliminaries}\label{g2prelim}
Throughout this section $M^7$ is a manifold with a $G_2$ structure. That means that $M$ admits a non-degenerate $3$-form $\varphi$ (see~\cite{Sp1} for more details).  Note that $\varphi$ determines a metric $g$ and orientation, hence also the Hodge-star $\star$. Then we also have that $\psi \coloneqq \star \varphi$ is a non-degenerate $4$-form. First, we list some results for manifolds with a $G_2$-structure.

\begin{prop}\label{g2idenflows} We use the following identities from~\cite{Sp1}:
\begin{itemize}
\item $\varphi_{ijk}\varphi_{abk} = \delta_{ia} \delta_{jb}-\delta_{ib} \delta_{ja} - \psi_{ijab}.$
\item $\varphi_{ijk} \psi_{abck} = \delta_{ia} \varphi_{jbc}+\delta_{ib} \varphi_{ajc}+\delta_{ic} \varphi_{abj}-\delta_{aj} \varphi_{ibc}-\delta_{bj} \varphi_{aic}-\delta_{cj} \varphi_{abi}.$
\item $\varphi_{ijk} \psi_{abjk} =-4 \varphi_{iab}.$
\item $\psi_{ijkl} \psi_{abkl} =4 \delta_{ia} \delta_{jb}-4\delta_{ib} \delta_{ja}-2\psi_{ijab}.$
\item $\psi_{ijkl} \psi_{ajkl} =24\delta_{ia}.$
\end{itemize}
\end{prop}

\begin{rmk}\label{g2description}
We have the following descriptions of the orthogonal decompositions of $\Omega^2$ and $\Omega^3$ into irreducible subspaces (see~\cite{Sp1}):\\
\begin{align*}
\Omega^2 &= \Omega^2_7 \oplus \Omega^2_{14}\\
\Omega^3 &= \Omega^3_1 \oplus \Omega^3_{7} \oplus \Omega^3_{27},
\end{align*}
where the indices denote the corresponding dimensions. In particular:
\begin{itemize}
\item $\Omega^2_7 = \{X \intprod \varphi: X \in \Gamma(TM) \} = \{\beta \in \Omega^2: \star(\varphi \wedge \beta) = -2 \beta \},$
\end{itemize}
\tab or equivalently, $\beta \in \Omega^2_7$ iff $\beta_{ij} \psi_{ijkl} = -4 \beta_{kl} \Leftrightarrow  \beta_{ij} = X_k \varphi_{ijk} \text{ for } X_k = \frac{1}{6} \beta_{ij} \varphi_{ijk}$
\begin{itemize}
\item $\Omega^2_{14} = \{\beta \in \Omega^2: \beta \wedge \psi =  0 \} = \{\beta \in \Omega^2: \star(\varphi \wedge \beta) =  \beta \},$
\end{itemize}
\tab or equivalently, $ \beta \in \Omega^2_{14}$ iff $\beta_{ij} \psi_{ijkl} = 2 \beta_{kl} \Leftrightarrow  \beta_{ij}\varphi_{ijk}= 0$.
\begin{itemize}
\item $\Omega^3_{1} = \{f \varphi: f \in C^{\infty}(M)\} \cong \mathbb{R}g$,
\item $\Omega^3_{7} = \{X \intprod \psi: X \in \Gamma(TM) \} \cong \Omega^2_7$,
\item $\Omega^3_{27} \cong  \mathcal{S}^2_0$,
\end{itemize}
where the isomorphisms are obtained using the $\diamond$ map with $\varphi$ (although different notation is used in~\cite{Sp1}, instead of $\diamond$ there).\\
Similarly, we have isomorphisms between the irreducible subpspaces of $\Omega^4$ and $\mathcal{S}^2 \oplus \Omega^2_7$ via $\diamond$ with $ \psi.$  In particular, let $\gamma \in \Omega^3, \zeta \in \Omega^4$. Then $\gamma = A \diamond \varphi$ and $\zeta = B \diamond \psi$ for some unique $A = \frac{1}{7} (\trr A) g + A_0 + A_7, B = \frac{1}{7}(\trr B) g + B_0 + B_7 \in \mathcal{S}^2 \oplus \Omega^2_7$. Define $$\hat{\gamma}_{ia} \coloneqq \gamma_{ijk} \varphi_{ajk} \text{\tab and\tab}\hat{\zeta}_{ia} \coloneqq \zeta_{ijkl} \psi_{ajkl}.$$ Then:
\begin{align*}
\trr{A} &= \frac{1}{18} \trr(\hat{\gamma}),\\
(A_0)_{ia} &= \frac{1}{8} (\hat{\gamma}_{ia} + \hat{\gamma}_{ai}) - \frac{1}{28}\trr(\hat{\gamma})g_{ia},\\
(A_7)_{ia} &= \frac{1}{24} (\hat{\gamma}_{ia} - \hat{\gamma}_{ai}).
\end{align*}
and we have similar formulas for $B$, but they will not be used here.
\end{rmk}

\begin{defn}
A manifold $M$ with a $G_2$ structure $\varphi$ has four independent torsion forms corresponding to a $G_2$ structure $\varphi$:
$$\tau_0 \in C^{\infty}(M),\tab \tau_1 \in \Omega^1_7, \tab\tau_2 \in \Omega^2_{14}, \tab\tau_3 \in \Omega^{3}_{27},$$
defined by the equations:
\begin{align*}
d \varphi &= \tau_0 \psi + 3 \tau_1 \wedge \varphi + \star \tau_3\\
d \psi &= 4 \tau_1 \wedge \psi + \star \tau_2.
\end{align*}
We say that $M$ is \textbf{nearly $G_2$} if $d \varphi = \tau_0 \psi$ for some constant $\tau_0 \neq 0$ and $d\psi = 0.$ A priori, $\tau_0$ could be a function on $M$ but it can be deduced that it has to be constant. Note that the condition of being nearly $G_2$ is also equivalent to the fact that the only non-zero component of the torsion tensor is $\tau_0$. These manifolds are positive Einstein and one might want to scale the metric so that $\tau_0 = 4$ (in this case we will also have $\operatorname{Ric} = 6g$), as we do for the nearly  K\"{a}hler case. However, we keep it more general.
\end{defn}

\begin{prop}\label{g2f}
For a nearly $G_2$ manifold we have the following formulas:
\begin{align*}
\nabla_p \varphi_{ijk} =& \frac{\tau_0}{4} \psi_{pijk},\\
\nabla_p \psi_{ijkl} =&  \frac{\tau_0}{4}(\delta_{lp} \varphi_{ijk}+\delta_{jp} \varphi_{ikl}-\delta_{kp} \varphi_{ijl}-\delta_{ip} \varphi_{ljk}),\\
\sum_p \nabla_p \nabla_p \varphi_{ijk} =& - \frac{\tau_0^2}{4} \varphi_{ijk}.
\end{align*}
\begin{proof}
The first two formulas are in~\cite{Sp1}. The third formula is demonstrated in Proposition 2.4 of~\cite{AS}.
\end{proof}
\end{prop}

\subsection{Curvature identities}\label{ng2jkcurg}
On a nearly $G_2$ manifold we have: $Ric = \frac{3\tau_0^2}{8}g$ (also see~\cite{Sp1}), so the Einstein constant $k=\frac{3\tau^2_0}{8}$ and $R= \frac{21\tau_0^2}{8}$.
Applying the result from Lemma~\ref{hatcirclemma} we get:
\begin{equation}\label{hatlemmag2}
\begin{split}
\hat{W} &= \hat{R} + \frac{\tau_0^2}{8} \operatorname{Id},\\
\mathring{W} &= \mathring{R} - \frac{\tau_0^2}{16} \operatorname{Id}, \text{ on $\mathcal{S}^2_0.$}
\end{split}
\end{equation}

Also, note that the Weyl tensor $W_{ijkl}$ lies in $\Omega^2_{14}$ in the first two or the last two indices.
This is because from Theorem $4.2$ in~\cite{Sp1}, we have $R_{ijkl}\varphi_{klm} = -\frac{\tau^2_0}{8} \varphi_{ijm}.$ By Remark~\ref{extendhat}, we can write this as $\hat{R} \varphi = -\frac{\tau^2_0}{8} \varphi.$ By~\eqref{hatlemmag2}, we get $\hat{W} \varphi= \hat{R} \varphi +  \frac{\tau_0^2}{8} \varphi = 0.$ By Remark~\ref{g2description}, this is equivalent to the fact that $W$ lies in $\Omega^2_{14}$ in the first two indices. Because of its symmetries, the same holds in the last two indices.\\
Hence, we can also conclude that $\hat{W}$, $\hat{R}$ preserve the space $\Omega^2_{14}.$ Consider $\hat{W}$ first. The $2$-form $(\hat{W} \beta)_{ab} = W_{abij} \beta_{ij}$ will always lie in $\Omega^2_{14}$. So, vacuously, it preserves $\Omega^2_{14}.$ Next, since $\hat{R}$ and $\hat{W}$ differ by a constant multiple of the identity, $\hat{R}$ also preserves $\Omega^2_{14}$. This fact means that we can consider $\hat{W}$ (and $\hat{R}$) as a self-adjoint operator only on $\Omega^2_{14}$ which will provide better estimates when we apply the Bochner-Weitzenböck techniques.

\subsection{Weitzenböck formulas}\label{ng2weindgfgf}
In this section we establish sufficient conditions for the Betti numbers $b_2$ or $b_3$ to vanish, in terms of bounds on $\hat{W}$ and $\mathring{W}$ respectively. As a corollary, we can get those sufficient conditions in terms of bounds on $\hat{R}$.\\
The simplified Weitzenböck formulas obtained in this section can be found in the literature, but possibly in different forms (see~\cite{AS}). As we mentioned in the introduction, we again reprove all the results in a simple, direct way.\\
We will use Theorems 3.8 and 3.9 from~\cite{SD}, which state the every harmonic $2$-form lies in $ \Omega^2_{14}$, and every harmonic $3$-form lies in $ \Omega^3_{27}.$

\subsubsection{$2$-forms}

We apply Corollary~\ref{weiz2simpl} to the nearly $G_2$ setting to get:
\begin{equation}\label{eq:weiz2forms2g2}
\Delta \beta =\nabla^{*} \nabla \beta + \frac{5 \tau^2_0}{8} \beta +  \hat{W}\beta \text{, for any $\beta\in \Omega^2.$}
\end{equation}

\begin{thm}\label{thmm2}
Let $M$ be a compact nearly $G_2$ manifold. If $\mathcal{S}^2(\Omega^2_{14}) \ni \hat{W} \geq - \frac{5 \tau^2_0}{8}$, then $b_2 = 0.$
\begin{proof}
Let $\beta \in \Omega^2$ be harmonic. Then $\beta \in \Omega^2_{14}.$ Substituting it in~\eqref{eq:weiz2forms2g2}, and using the assumption that $\hat{W} \geq -\frac{5 \tau^2_0}{8}$, we get by integration that $\nabla \beta = 0.$ Hence $\beta = 0,$ as there are no parallel non-zero $2$-forms.
\end{proof}
\end{thm}

\begin{thm}\label{aerstdfjhgggggg}
Let $M$ be a compact nearly $G_2$ manifold. Let $\delta \leq \bar{R} \leq \Delta$ with $-(\Delta + \delta) -\frac{7}{3} (\Delta - \delta) \geq -\frac{3\tau_0^2}{4}$. Then $b_2=0$.
\begin{proof} 
If $-(\Delta + \delta) -\frac{7}{3} (\Delta - \delta) \geq -\frac{3\tau_0^2}{4}$, then by Corollary~\ref{improvedaa} we have that $\hat{R} \geq -\frac{3\tau_0^2}{4}$. So, we use equation~\eqref{hatlemmag2} to get that $\hat{W} \geq - \frac{5 \tau^2_0}{8}$ and hence $b_2 = 0$ by Theorem~\ref{thmm2}.
\end{proof}
\end{thm}

\subsubsection{$3$-forms}
We apply Corollary~\ref{weiz3simpl} to the nearly $G_2$ setting to get:
\begin{equation}\label{eq:W3}
(\Delta \beta)_{abc} =(\nabla^{*} \nabla \beta)_{abc} + \frac{3}{4}\tau_0^2 \beta_{abc} +  W_{abpu} \beta_{puc}+ W_{acpu} \beta_{pbu}+ W_{bcpu} \beta_{apu} \text{, for any $\beta\in \Omega^3.$}
\end{equation}

The aim now is to simplify this formula for harmonic $\beta$. First, we do this more generally, we will just assume $\beta \in \Omega^3_1 \oplus \Omega^3_{27}$, which includes all the harmonic forms. Then we will see what we can get from the assumption of $\beta$ being harmonic and then we use all these steps to get a simpler formula.\\
Recall definitions of $\divv$ and $\diamond$ from Section~\ref{notsec}, and also for $h \in \mathcal{S}^2$, let $\tilde{h} \in \mathcal{T}^2$ be defined as
\begin{equation}\label{eq:tildedefinition}
\tilde{h}_{kc}=(\nabla_i h_{jk}) \varphi_{ijc}.
\end{equation}

\begin{prop}\label{W3h}
Let $M$ be nearly $G_2$. Let $\beta \in \Omega^3_1 \oplus\Omega^3_{27}$, so that $\beta = h \diamond \varphi$ for some $h\in \mathcal{S}^2$. Then:
$$\Delta (h \diamond \varphi) = (\nabla^* \nabla h + \tau_0^2 h + \frac{\tau_0}{2} \tilde{h}_{symm} - \frac{\tau_0}{12} ((\divv h-\nabla \trr h) \intprod \varphi) +2\mathring{W}h ) \diamond \varphi.$$
\begin{proof}
First, consider the term $W_{abpu} \beta_{puc}$ from~\eqref{eq:W3}:
\begin{align*}
W_{abpu} \beta_{puc} &= W_{abpu} (h_{ps} \varphi_{suc} +h_{us} \varphi_{psc}+h_{cs} \varphi_{pus})\\
&= 2 W_{abpu} h_{ps} \varphi_{suc},
\end{align*}
as the first two terms in the brackets are skew in $p,u$ and the last term vanishes because $W \in \mathcal{S}^2(\Omega^2_{14})$.
Hence, 
$$W_{abpu} \beta_{puc}+ W_{acpu} \beta_{pbu}+ W_{bcpu} \beta_{apu} = 2h_{ps} (W_{abpu} \varphi_{suc}+W_{acpu} \varphi_{sbu}+W_{bcpu} \varphi_{asu}), $$
and define $\gamma_{abc}$ to be equal to this expression. Since $\gamma \in \Omega^3, \gamma = A \diamond \varphi$ for some $A \in \mathcal{S}^2 \oplus \Omega^2_7.$ We have:
\begin{align*}
\hat{\gamma}_{at} =& \gamma_{abc}\varphi_{tbc}\\
=& 2h_{ps} (W_{abpu} \varphi_{suc}+ W_{acpu} \varphi_{sbu}+W_{bcpu} \varphi_{asu}) \varphi_{tbc}\\
=& 4h_{ps} (W_{abpu} \varphi_{suc} \varphi_{tbc}) \text{ (since $W \in \mathcal{S}^2(\Omega^2_{14})$ and we use  skew-symmetry in $b,c$ on the first two terms)}\\
=& 4  h_{ps} W_{abpu} (\delta_{st}\delta_{ub}-\delta_{sb}\delta_{ut} - \psi_{sutb})\text{ (by Proposition~\ref{g2idenflows})}\\
=&-4 (h_{ps} W_{aspt} + h_{ps} W_{abpu} \psi_{sutb}). \text{ (because the Ricci tensor of $W$ is zero, i.e. $W_{abbu} = 0$)}
\end{align*}
Now, note that:
\begin{align*}
W_{abpu} \psi_{sutb} &= - (W_{apub} + W_{aubp}) \psi_{sutb}\\
&= W_{apub}\psi_{stub} - W_{aubp}\psi_{sutb}\\
&= 2 W_{apst} - W_{abup}\psi_{sbtu} \text{ (here we swap the indices $b,u$ and use that $W \in S^2(\Omega^2_{14})$ with Remark~\ref{g2description})}\\
&= 2 W_{apst} - W_{abpu}\psi_{sutb}.
\end{align*}
Hence, we have
$$W_{abpu} \psi_{sutb} = W_{apst},$$
and thus
\begin{align*}
\hat{\gamma}_{at} &= - 4 ( h_{ps} W_{aspt} + h_{ps} W_{apst})\\
&= 8  (\mathring{W}h)_{at}.
\end{align*}
Next, we have that $\operatorname{tr}(\hat{\gamma}) = 0$ because $W_{ipqi} = 0$. Also, we see that $\hat{\gamma}$ is symmetric. Hence, by Remark~\ref{g2description}, $$A = A_0 = \frac{1}{4} \hat{\gamma} = 2\mathring{W}h.$$ Thus, the term with Weyl tensors in the Weitzenböck formula is equal to $ \gamma =  A \diamond \varphi = 2(\mathring{W}h) \diamond \varphi$.\\
Next, we compute $\nabla^* \nabla \beta = \nabla^* \nabla (h \diamond \phi)$ as follows:
\begin{align*}
\nabla^* \nabla (h \diamond \phi)_{abc} =&- \nabla_s \nabla_s (h_{ap} \varphi_{pbc} +h_{bp} \varphi_{apc}+h_{cp} \varphi_{abp})\\
=& -  (\nabla_s (\nabla_s h_{ap}) \varphi_{pbc} + \nabla_s (\nabla_s h_{bp}) \varphi_{apc}+ \nabla_s (\nabla_s h_{cp}) \varphi_{abp}))\\
&-2 (\nabla_s (h_{ap}) \nabla_s( \varphi_{pbc}) + \nabla_s(h_{bp}) \nabla_s(\varphi_{apc})+\nabla_s (h_{cp}) \nabla_s (\varphi_{abp}))\\
&-(h_{ap} \nabla_s \nabla_s(\varphi_{pbc}) +h_{bp} \nabla_s \nabla_s(\varphi_{apc})+h_{cp} \nabla_s \nabla_s(\varphi_{abp}))\\
=& ((\nabla^* \nabla h + \frac{\tau_0^2}{4} h) \diamond \varphi)_{abc} \text{ (by Proposition~\ref{g2f})}\\
&-\frac{\tau_0}{2} (\nabla_s (h_{ap}) \psi_{spbc}+ \nabla_s(h_{bp}) \psi_{sapc}+\nabla_s (h_{cp}) \psi_{sabp}).
\end{align*}
Let $\sigma_{abc} \coloneqq \nabla_s (h_{ap}) \psi_{spbc}+ \nabla_s(h_{bp}) \psi_{sapc}+\nabla_s (h_{cp}) \psi_{sabp} $. As $\sigma \in \Omega^3$, $\sigma = B \diamond \varphi$ for some unique $B \in  \mathcal{S}^2 \oplus \Omega^2_7.$ Then:
\begin{align*}
\hat{\sigma}_{at} &= \sigma_{abc} \varphi_{tbc}\\
&= \nabla_s (h_{ap}) \psi_{spbc} \varphi_{tbc} + 2 \nabla_s(h_{bp}) \psi_{sapc} \varphi_{tbc}\\
&= \nabla_s (h_{ap}) (-4 \varphi_{tsp}) + 2 \nabla_s(h_{bp}) (\delta_{ts}\varphi_{bap}+\delta_{ta}\varphi_{sbp}+\delta_{tp}\varphi_{sab}-\delta_{sb}\varphi_{tap}-\delta_{ab}\varphi_{stp}-\delta_{pb}\varphi_{sat})\text{ (by Proposition~\ref{g2idenflows})}\\
&=-4 \nabla_s (h_{pa}) \varphi_{spt}-2 \nabla_s(h_{bt}) \varphi_{sba} + 2 \nabla_s(h_{sp}) \varphi_{pat} +2 \nabla_s(h_{pa})\varphi_{spt}-2 \nabla_s(h_{pp}) \varphi_{sat}\\
&=-2\tilde{h}_{at}-2\tilde{h}_{ta}+2 (\divv h \intprod \varphi)_{at}-2  (\nabla \trr h \intprod \varphi)_{at}\\
&=-4(\tilde{h}_{symm})_{at}+2 ((\divv h-\nabla \trr h) \intprod \varphi)_{at}.
\end{align*}
Note that $\trr{\hat{\sigma}} = -4(\tilde{h}_{symm})_{aa}= -4\tilde{h}_{aa}=-4 \nabla_i (h_{ja}) \varphi_{ija} = 0$.
Thus, by Remark~\ref{g2description} we have:
$$(B_0)_{ia} = \frac{1}{8} (\hat{\sigma}_{ia}+\hat{\sigma}_{ai}) = - (\tilde{h}_{symm})_{ia},$$
$$(B_7)_{ia} = \frac{1}{24} (\hat{\sigma}_{ia}-\hat{\sigma}_{ai}) =  \frac{1}{6} ((\divv h-\nabla \trr h) \intprod \varphi)_{ia}.$$
We conclude that:
$$\sigma = (-\tilde{h}_{symm} + \frac{1}{6} ((\divv h-\nabla \trr h) \intprod \varphi)) \diamond \varphi.$$
Putting everything together we get:
\begin{align*}
\Delta (h \diamond \varphi) =& (\nabla^* \nabla h + \frac{\tau_0^2}{4} h - \frac{\tau_0}{2}  (-\tilde{h}_{symm} + \frac{1}{6} ((\divv h-\nabla \trr h) \intprod \varphi)) + \frac{3}{4}\tau_0^2 h + 2\mathring{W}h ) \diamond \varphi\\
=&(\nabla^* \nabla h + \tau_0^2 h + \frac{\tau_0}{2} \tilde{h}_{symm} - \frac{\tau_0}{12} ((\divv h-\nabla \trr h) \intprod \varphi) +2\mathring{W}h ) \diamond \varphi,
\end{align*}
as claimed.
\end{proof}
\end{prop}

\begin{prop}\label{harm-h}
Let $M$ be a compact nearly $G_2$ manifold. Let $\beta \in \Omega^3_{27} \supseteq \mathcal{H}^3$, so that $\beta = h \diamond \varphi$ for some $h\in \mathcal{S}^2_0$. Then $\beta$ is harmonic if and only if:
$$\tilde{h} = -\frac{3\tau_0}{4}h \in \mathcal{S}^2\text{\tab and \tab} \divv h= 0.$$
\begin{proof}
Since $M$ is compact, $\beta$ is harmonic if and only if $d^* \beta = 0$ and $d\beta = 0.$ \\
First, we calculate the $\Omega^2_7$ component of $\tilde{h}$. This is done by contracting it with $\varphi$. That is, by Remark~\ref{g2description}, $\pi_7(\tilde{h}) = X \intprod \varphi$, where $X_k = \frac{1}{6} \tilde{h}_{ij} \varphi_{kij}$. So:
\begin{align*}
X_k &= \frac{1}{6} \tilde{h}_{ij} \varphi_{kij}\\
&= \frac{1}{6} (\nabla_a h_{bi}) \varphi_{abj} \varphi_{kij}\\
&= \frac{1}{6} (\nabla_a h_{bi}) (\delta_{ak} \delta_{bi}-\delta_{ai} \delta_{bk} - \psi_{abki})\\
&= \frac{1}{6} (\nabla_k h_{bb} - \nabla_i h_{ki})\\
&= -\frac{1}{6} (\divv h)_k.
\end{align*}
Thus:
$$\pi_7(\tilde{h}) = X \intprod \varphi \text{\tab where\tab} X = - \frac{1}{6} \divv h.$$
Now, consider the condition $d^* \beta = 0$:
\begin{align*}
-(d^* \beta)_{kl} =& \nabla_j \beta_{jkl}\\
=& \nabla_j (h_{jp} \varphi_{pkl}+h_{kp} \varphi_{jpl}+h_{lp} \varphi_{jkp})\\
=& \nabla_j (h_{jp}) \varphi_{pkl} + h_{jp} \nabla_j (\varphi_{pkl})+\nabla_j (h_{kp}) \varphi_{jpl} + h_{kp} \nabla_j (\varphi_{jpl})+\nabla_j (h_{lp}) \varphi_{jkp} + h_{lp} \nabla_j (\varphi_{jkp})\\
=& (\divv h \intprod \varphi)_{kl} + h_{jp} \frac{\tau_0}{4} \psi_{jpkl}+\tilde{h}_{kl} + h_{kp} \frac{\tau_0}{4} \psi_{jjpl}-\tilde{h}_{lk} + h_{lp}\frac{\tau_0}{4} \psi_{jjkp}\\
=& (\divv h \intprod \varphi)_{kl}+ 2(\tilde{h}_{skew})_{kl}.
\end{align*}
Hence, using the formula for $\pi_7(\tilde{h})$, we get:
\begin{align*}
-d^* \beta =& \divv h \intprod \varphi + 2 (\pi_7(\tilde{h}) + \pi_{14}(\tilde{h}))\\
=&  \divv h \intprod \varphi - \frac{1}{3} \divv h \intprod \varphi + 2\pi_{14}(\tilde{h})\\
=& 2\pi_{14}(\tilde{h}) +  \frac{2}{3} \divv h \intprod \varphi.
\end{align*}
Thus, we have that:
\[d^* \beta = 0 \text{\tab if and only if\tab} \myarray{\pi_{14}(\tilde{h}) = 0, \\ \divv h = 0.}\]
Next, we consider the second condition $d\beta = 0$. Since $d\beta \in \Omega^4$, by Remark~\ref{g2description}, $d\beta = B \diamond \psi$ for some unique $B \in \mathcal{S}^2 \oplus \Omega^2_7$. Then $d \beta = 0$ iff $B =0$ iff $\widehat{d\beta} = 0$. So:
\begin{align*}
(\widehat{d\beta})_{ia} =& (d\beta)_{ijkl} \psi_{ajkl}\\
=&(\nabla_i \beta_{jkl}-\nabla_j \beta_{ikl}+\nabla_k \beta_{ijl}-\nabla_l \beta_{ijk}) \psi_{ajkl}\\
=&(\nabla_i \beta_{jkl}-3 \nabla_j \beta_{ikl}) \psi_{ajkl}\\
=& (\nabla_i(h_{jp} \varphi_{pkl}+h_{kp} \varphi_{jpl}+h_{lp} \varphi_{jkp})) \psi_{ajkl}\\
&-3  (\nabla_j(h_{ip} \varphi_{pkl}+h_{kp} \varphi_{ipl}+h_{lp} \varphi_{ikp})) \psi_{ajkl}\\
=& 3 \nabla_i(h_{jp} \varphi_{pkl}) \psi_{ajkl}-3 \nabla_j(h_{ip} \varphi_{pkl}) \psi_{ajkl}- 6\nabla_j(h_{kp} \varphi_{ipl}) \psi_{ajkl}.
\end{align*}
We calculate each of these terms separately:
\begin{align*}
3 \nabla_i(h_{jp} \varphi_{pkl}) \psi_{ajkl} &= 3 \nabla_i(h_{jp}) \varphi_{pkl} \psi_{ajkl} + 3 h_{jp} \nabla_i(\varphi_{pkl}) \psi_{ajkl}\\
&= -12 \nabla_i(h_{jp}) \varphi_{paj} + \frac{3 \tau_0}{4} h_{jp} \psi_{ipkl} \psi_{ajkl}\\
&= \frac{3 \tau_0}{4} h_{jp} (4 \delta_{ia} \delta_{pj}-4 \delta_{ij} \delta_{pa} - 2 \psi_{ipaj})\text{ (by Proposition~\ref{g2idenflows})}\\
&= 3 \tau_0 ((\trr{h}) \delta_{ia} - h_{ia})\\
&= -3 \tau_0  h_{ia}.
\end{align*}
Next:
\begin{align*}
-3 \nabla_j(h_{ip} \varphi_{pkl}) \psi_{ajkl} &=-3 \nabla_j(h_{ip}) \varphi_{pkl} \psi_{ajkl}-3 h_{ip} \nabla_j(\varphi_{pkl}) \psi_{ajkl}\\
&=12 \nabla_j(h_{ip}) \varphi_{paj} -\frac{3 \tau_0}{4} h_{ip} \psi_{jpkl} \psi_{ajkl}\\
&=12 \nabla_j(h_{ip}) \varphi_{paj} +\frac{3 \tau_0}{4} h_{ip} 24 \delta_{pa}\text{ (by Proposition~\ref{g2idenflows})}\\
&=12 \tilde{h}_{ia} +18 \tau_0 h_{ia}.
\end{align*}
Finally:
\begin{align*}
- 6\nabla_j(h_{kp} \varphi_{ipl}) \psi_{ajkl} =& - 6\nabla_j(h_{kp}) \varphi_{ipl} \psi_{ajkl}- 6h_{kp} \frac{\tau_0}{4} \psi_{jipl} \psi_{ajkl}\\
=& - 6\nabla_j(h_{kp}) \varphi_{ipl} \psi_{ajkl}- \frac{3\tau_0}{2}h_{kp} \psi_{pijl} \psi_{akjl}\\
=& - 6\nabla_j(h_{kp}) (\delta_{ia} \varphi_{pjk}+\delta_{ij} \varphi_{apk}+\delta_{ik} \varphi_{ajp}-\delta_{ap} \varphi_{ijk}-\delta_{jp} \varphi_{aik}-\delta_{kp} \varphi_{aji})\\
&- \frac{3\tau_0}{2}h_{kp} (4 \delta_{pa}\delta_{ik}-4 \delta_{pk}\delta_{ia}-2\psi_{piak})\text{ (by Proposition~\ref{g2idenflows})}\\
=&-6(\nabla_j(h_{ip})\varphi_{ajp}-\nabla_j(h_{ka})\varphi_{ijk}-\nabla_j(h_{kj})\varphi_{aik}-\nabla_j(\trr h)\varphi_{aji})-6\tau_0 (h_{ia}-(\trr h)\delta_{ia})\\
=& -6 \tilde{h}_{ia}+6 \tilde{h}_{ai}-6(\divv h \intprod \varphi)_{ia} -6\tau_0 h_{ia}.
\end{align*}
Combining all the results we get:
\begin{align*}
(\widehat{d\beta})_{ia} =&~6 (\tilde{h}_{ia} + \tilde{h}_{ai}) + 9\tau_0 h_{ia}  -6 (\divv h \intprod \varphi)_{ia}\\
=&~12 (\tilde{h}_{symm})_{ia} + 9\tau_0 h_{ia} - 6 (\divv h  \intprod \varphi)_{ia}.
\end{align*}
Thus, we have that:
\[d \beta = 0 \text{\tab if and only if\tab} \myarray{\tilde{h}_{symm} = -\frac{3\tau_0}{4} h, \\ \divv h = 0.}\]
Summarizing, $\beta$ is harmonic if and only if $\tilde{h}_{symm} = -\frac{3\tau_0}{4}h, \pi_{14}(\tilde{h}) = 0, \divv h =0$. But we know that $\divv h$ vanishes if and only if $\pi_{7}(\tilde{h})$ vanishes, so $\tilde{h}_{skew}=0$ and $\tilde{h} = \tilde{h}_{symm}$. Hence we get the required result.
\end{proof}
\end{prop}

\begin{cor}\label{corh}
Let $M$ be a compact nearly $G_2$ manifold. Assume $\beta$ is harmonic. Then $\beta \in \Omega^3_{27}$, so that $\beta = h \diamond \varphi$ for some $h\in \mathcal{S}^2_0$. Then:
$$ \nabla^* \nabla h + \frac{5\tau_0^2}{8} h +2\mathring{W}h = 0.$$
\begin{proof}
By assumption, $\beta$ is harmonic, so the left hand side of the Weitzenböck formula in Proposition~\ref{W3h} vanishes. Next, by Proposition~\ref{harm-h}, $\tilde{h}_{symm} = \tilde{h} = -\frac{3\tau_0}{4}h$ and $\divv h = 0$. Also, we know $\trr h = 0$. Substituting all these terms into Proposition~\ref{W3h}, we get the required result.
\end{proof}
\end{cor}
\begin{thm}\label{thmm1}
Let $M$ be a compact nearly $G_2$ manifold. If $\mathcal{S}^2(\mathcal{S}^2_0) \ni \mathring{W}  \geq - \frac{3\tau^2_0}{8}$\textbf{or} $\mathcal{S}^2(\Omega^2_{14}) \ni \hat{W} \geq -\frac{\tau^2_0}{4}$,  then $b_3 = 0.$
\begin{proof}
For the first part, we use~\eqref{eq:W3} and the proof of Proposition~\ref{W3h} to get that if $\beta = h \diamond \varphi \in \Omega^3_{27}$ is harmonic for some $h \in \mathcal{S}^2_0,$ then:
$$\nabla^* \nabla \beta + (\frac{3\tau^2_0}{4} h + 2\mathring{W} h)\diamond \varphi = 0.$$
Hence, the result follows by intergration and the fact that there are no parallel non-zero $3$-forms. Note that using Corollary~\ref{corh} in order to get a similar result would have been worse, as we would have been able to only conclude that if $\mathcal{S}^2(\mathcal{S}^2_{0}) \ni \mathring{W} \geq - \frac{5\tau^2_0}{16}$ then $b_3 = 0.$ This is because  $\nabla^* \nabla \beta = (\nabla^* \nabla h - \frac{\tau^2_0}{8} h) \diamond \varphi$, so we can see that even though the left hand side is obviously non-negative, we cannot conclude that from the right hand side.\\
The second part trivially follows from~\eqref{eq:W3}.
\end{proof}
\end{thm}
\begin{thm}\label{sgfd}
Let $M$ be a compact nearly $G_2$ manifold. Let $\delta \leq \bar{R} \leq \Delta$ with $\Delta \leq \frac{11\tau^2_0}{80}$ \textbf{or} $\delta \geq \frac{\tau^2_0}{112}$. Then $b_3=0$.
\begin{proof}
Recall that the Einstein constant $k = \frac{3 \tau^2_0}{8}.$ Then by Theorem~\ref{nkmaincor}, on $\mathcal{S}^2_0$, $\mathring{R} \geq -\frac{3\tau^2_0}{8} + 7\delta$ and $\mathring{R} \geq \frac{3\tau^2_0}{8} - 5\Delta$. Hence, by~\eqref{hatlemmag2}, $\mathring{W} \geq -\frac{7\tau^2_0}{16} + 7\delta$ and $\mathring{W} \geq \frac{5\tau^2_0}{16} - 5\Delta$.\\
In order for $b_3 = 0$, by Theorem~\ref{thmm1} we want $\mathring{W}  \geq - \frac{3\tau^2_0}{8}$. We have $-\frac{7\tau^2_0}{16} + 7\delta \geq - \frac{3\tau^2_0}{8}$  iff $\delta \geq \frac{\tau^2_0}{112}$; and $\frac{5\tau^2_0}{16} -5 \Delta \geq - \frac{3\tau^2_0}{8}$  iff $\Delta \leq \frac{11\tau^2_0}{80}$. Hence, the result follows from Theorem~\ref{thmm1}. Recall that, a priori, by Remark~\ref{uselessfact}, we have that $\delta \leq \frac{\tau^2_0}{16} \leq \Delta$.\\Also, note that we do not use Corollary~\ref{improvedaa} along with the statement that $\mathcal{S}^2(\Omega^2_{14}) \ni \hat{W} \geq -\frac{\tau^2_0}{4}$ implies that $b_3 = 0.$ This is because the sufficient conditions in terms of the bounds on the sectional curvature we would have obtained imply that $\Delta \leq \frac{11\tau^2_0}{80}$ or $\delta \geq \frac{\tau^2_0}{112}$.
\end{proof}
\end{thm}
\section{Nearly K\"{a}hler $6$-manifolds}
First, in Sections~\ref{nkprelim} and~\ref{diamsection} we establish some preliminaries for $6$-manifolds with an $SU(3)$-structure. In particular, we give various descriptions of irreducible subspaces of $\Omega^2$ using the $\diamond$ map. Next, in Section~\ref{nkmanifoldsintro} we introduce nearly K\"{a}hler $6$-manifolds and in Sections~\ref{curvvv} and \ref{harmshuidfh} we establish several identities involving curvature and harmonic $2$ and $3$-forms. Finally, in Section~\ref{weifojahsd} we simplify the Weizenbock formulas for harmonic $2$ and $3$-forms  and using the assumption of compactness of our manifolds, we get the necessary conditions of vanishing of $b_2$ and $b_3$ in terms of bounds on $\bar{R}, \mathring{R},$ and $\bar{R}$.

\subsection{Preliminaries}\label{nkprelim}
For this section as well, most of the results can be found in~\cite{LF},~\cite{GR},~\cite{CWMYK},~\cite{AMPANUW},~\cite{AMPANUW2} and other sources on nearly K\"{a}hler manifolds. Nevertheless, we include as many details as possible.\\
First we consider a general $SU(3)$-structure on a Riemannian $6$-manifold $(M,g)$. This means that $M^6$ has an almost complex structure $J$ compatible with the metric $g$, and a complex $3$-form $\Omega = \psi^{+} + i\psi^{-}$ satisfying $$\psi^{+} \wedge \psi^{-} = \frac{2}{3} \omega^3 = 4 \operatorname{vol_M}.$$
Also in this case, $g_C$, $\varphi$ and $\psi$ defined as:
\begin{align*}
g_C & \coloneqq r^2 g + dr^2,\\
\varphi &\coloneqq -r^2 dr \wedge \omega + r^3 \psi^{+},\\
\psi &\coloneqq - r^3 dr \wedge \psi^{-} - r^4 \frac{w^2}{2}.
\end{align*}
give a metric cone $G_2$-structure on $\mathbb{R}^{+}\times M.$\\
Hence, in a local orthonormal frame we get:
\begin{align}\label{phipsiri}
\varphi_{0ij} &= -\omega_{ij} \text{\tab} \psi_{0ijk} = - \psi^{-}_{ijk}\\
\varphi_{ijk} &= \psi^{+}_{ijk} \text{\tab}\psi_{ijkl} =- (\star \omega)_{ijkl}\nonumber.
\end{align}

We also list the following identities that hold, without proof:
\begin{equation}\label{su3iden}
\omega_{ik} \omega_{il} = \delta_{kl}, \star \psi^{+} = \psi^{-}, \star \psi^{-} = -\psi^{+}, \star \omega = \frac{1}{2} \omega^2.
\end{equation}
Looking at the last identity in coordinates gives us:
$$(\star \omega)_{ijkl} = \omega_{ij} \omega_{kl} + \omega_{jk} \omega_{il} +\omega_{lj} \omega_{ik}.$$

\begin{prop}\label{nkup}
The following identities hold:
\begin{itemize}
\itemequation[eq:pi1]{$\psi^{+}_{ijk} \omega_{ak} = -\psi^{-}_{ija}.$}{}
\itemequation[eq:pi2] {$\psi^{-}_{ijk} \omega_{ak} = \psi^{+}_{ija}.$}{}
\itemequation[eq:pinew] {$\psi^{+}_{ijk} \omega_{jk} = 0 =\psi^{-}_{ijk} \omega_{jk}.$}{}
\itemequation[eq:pi3] {$\psi^{+}_{ijk} \psi^{+}_{abk} = \delta_{ia} \delta_{jb}-\delta_{ib} \delta_{ja} - \omega_{ia} \omega_{jb} +\omega_{ib} \omega_{ja}.$}{}
\itemequation[eq:pi4] {$\psi^{+}_{ijk} \psi^{+}_{ajk} = 4 \delta_{ik}.$ (contraction of the previous one)}{}
\itemequation[eq:pi5] {$\psi^{+}_{ijk} \psi^{-}_{abk} = \delta_{ia} \omega_{jb}+\delta_{jb} \omega_{ia}-\delta_{ib} \omega_{ja}-\delta_{ja} \omega_{ib}.$}{}
\itemequation[eq:pi6] {$\psi^{+}_{ijk} \psi^{-}_{ajk} = 4 \omega_{ia}.$ (contraction of the previous one)}{}
\itemequation[eq:pi7] {$\psi^{-}_{ijk} \psi^{-}_{abk} = \delta_{ia} \delta_{jb}-\delta_{ib} \delta_{ja} - \omega_{ia} \omega_{jb} +\omega_{ib} \omega_{ja}.$ (same as $\psi^{+}_{ijk} \psi^{+}_{abk}$)}{}
\itemequation[eq:pi8] {$\psi^{-}_{ijk} \psi^{-}_{ajk} = 4 \delta_{ik}.$ (contraction of the previous one)}{}
\itemequation[eq:pi9] {$\omega_{ik} (\star \omega)_{abck} =  \delta_{ia} \omega_{bc}+\delta_{ib} \omega_{ca}+\delta_{ic} \omega_{ab}.$}{}
\itemequation[eq:pi10] {$\omega_{ik} (\star \omega)_{abik} =  4 \omega_{ab}.$ (contraction of the previous one)}{}
\itemequation[eq:pi11] {$\psi^{+}_{ijk} (\star \omega)_{abck} =- \delta_{ia} \psi^{+}_{jbc}-\delta_{ib} \psi^{+}_{ajc}-\delta_{ic} \psi^{+}_{abj}+\delta_{aj} \psi^{+}_{ibc}+\delta_{bj} \psi^{+}_{aic}+\delta_{cj} \psi^{+}_{abi} - \omega_{ij} \psi^{-}_{abc}.$}{}
\itemequation[eq:pi12] {$\psi^{+}_{ijk} (\star \omega)_{abck} = -\psi^{-}_{ija} \omega_{bc}-\psi^{-}_{ijb} \omega_{ca}-\psi^{-}_{ijc} \omega_{ab}.$ (alternative expression to the previous one)}{}
\itemequation[eq:pi13] {$\psi^{+}_{ijk} (\star \omega)_{abjk} =2\psi^{+}_{iab}.$ (contraction of the previous one)}{}
\itemequation[eq:pi14] {$\psi^{-}_{ijk} (\star \omega)_{abck}=\psi^{+}_{ija} \omega_{bc}+\psi^{+}_{ijb} \omega_{ca}+\psi^{+}_{ijc} \omega_{ab}.$}{}
\itemequation[eq:pi15] {$\psi^{-}_{ijk} (\star \omega)_{abjk} = 2\psi^{-}_{iab}.$ (contraction of the previous one)}{}
\itemequation[eq:pi16] {$(\star \omega)_{ijkl}  (\star \omega)_{abkl} = 2 \delta_{ia}\delta_{jb} - 2 \delta_{ib}\delta_{ja}+2 \omega_{ij} \omega_{ab}.$}{}
\itemequation[eq:pi17] {$(\star \omega)_{ijkl}  (\star \omega)_{ajkl} = 12 \delta_{ia}.$ (contraction of the previous one)}{}
\end{itemize}
\begin{proof}
We repeatedly use~\eqref{phipsiri} along with $G_2$-contraction identities. For~\eqref{eq:pi1}:
\begin{align*}
\psi^{+}_{ijk} \omega_{ak} &= \sum_{k=1}^{6} \varphi_{ijk}(-\varphi_{0ak})\\
&= -\sum_{k=0}^{6} \varphi_{ijk}\varphi_{0ak}\\
&= -(\delta_{i0} \delta_{ja}-\delta_{ia} \delta_{j0} - \psi_{ij0a})\\
&= \psi_{0ija}\\
&= -\psi^{-}_{ija}.
\end{align*}
Contracting both sides of~\eqref{eq:pi1} with $w_{au}$ we get:
\begin{align*}
\psi^{+}_{ijk} \omega_{ak} \omega_{au} &= -\psi^{-}_{ija} \omega_{au}\\
\psi^{+}_{ijk} \delta_{ku} &= \psi^{-}_{ija} \omega_{ua}\\
\psi^{+}_{iju} &=\psi^{-}_{ija} \omega_{ua},
\end{align*}
which gives us~\eqref{eq:pi2}.\\
Contracting on~\eqref{eq:pi1} and~\eqref{eq:pi2} on $j,a$ immediately gives~\eqref{eq:pinew}.\\
Next, for~\eqref{eq:pi3}:
\begin{align*}
\psi^{+}_{ijk} \psi^{+}_{abk} &= \sum_{k=1}^{6} \varphi_{ijk}\varphi_{abk}\\
&= \sum_{p=0}^{6} \varphi_{ijk}\varphi_{abk}- \varphi_{ij0}\varphi_{ab0}\\
&=(\delta_{ia} \delta_{jb}-\delta_{ib} \delta_{ja} - \psi_{ijab}) - \omega_{ij}  \omega_{ab}\\
&=(\delta_{ia} \delta_{jb}-\delta_{ib} \delta_{ja}) +( (\star \omega)_{ijab}- \omega_{ij}  \omega_{ab}).\\
&=\delta_{ia} \delta_{jb}-\delta_{ib} \delta_{ja} + \omega_{ja} \omega_{ib} +\omega_{bj} \omega_{ia}\text{ (by Lemma~\ref{staromegalemma})}\\
&=\delta_{ia} \delta_{jb}-\delta_{ib} \delta_{ja} - \omega_{ia} \omega_{jb} +\omega_{ib} \omega_{ja}.
\end{align*}
Contracting~\eqref{eq:pi3} on $j,b$ gives us~\eqref{eq:pi4}:
\begin{align*}
\psi^{+}_{ijk} \psi^{+}_{ajk} &= 6\delta_{ia} -\delta_{ia} + \omega_{ij} \omega_{ja}\\
&= 6\delta_{ia} -\delta_{ia}  -\delta_{ia}\\
&= 4\delta_{ia}.
\end{align*}
Next, for~\eqref{eq:pi5}:
\begin{align*}
\psi^{+}_{ijk} \psi^{-}_{abk} &= \sum_{k=1}^{6} \varphi_{ijk} (-\psi_{0abk})\\ 
&= -\sum_{k=0}^{6} \varphi_{ijk} \psi_{0abk}\\ 
&= -(\delta_{i0} \varphi_{jab}+\delta_{ia} \varphi_{0jb}+\delta_{ib} \varphi_{0aj}-\delta_{0j} \varphi_{iab}-\delta_{aj} \varphi_{0ib}-\delta_{bj} \varphi_{0ai})\\
&= \delta_{ia} \omega_{jb}+\delta_{ib} \omega_{aj}-\delta_{aj} \omega_{ib}-\delta_{bj} \omega_{ai}\\
&= \delta_{ia} \omega_{jb}+\delta_{jb} \omega_{ia}-\delta_{ib} \omega_{ja}-\delta_{ja} \omega_{ib}.
\end{align*}
Contracting~\eqref{eq:pi5} on $j,b$ will give us~\eqref{eq:pi6}:
\begin{align*}
\psi^{+}_{ijk} \psi^{-}_{ajk} &= 6 \omega_{ia} - \omega_{ia} - \omega_{ia}\\
&= 4\omega_{ia}.
\end{align*}
Next, we show that $\psi^{-}_{ijk} \psi^{-}_{abk} = \psi^{+}_{ijk} \psi^{+}_{abk}$, which means that~\eqref{eq:pi7} and~\eqref{eq:pi3} are the same. Using~\eqref{eq:pi1}, we have:
\begin{align*}
\psi^{-}_{ijk} \psi^{-}_{abk} &=\psi^{+}_{ijs} \omega_{ks} \psi^{+}_{abt} \omega_{kt}\\
&= \psi^{+}_{ijs} \psi^{+}_{abt} \delta_{st}\\
&=\psi^{+}_{ijs} \psi^{+}_{abs}.
\end{align*}
Thus, we also get~\eqref{eq:pi8}: $$\psi^{-}_{ijk} \psi^{-}_{ajk} = \psi^{+}_{ijk} \psi^{+}_{ajk} = 4 \delta_{ia}.$$
Next, for~\eqref{eq:pi9}:
\begin{align*}
\omega_{ik} (\star \omega)_{abck} &= \sum_{k=1}^6 (-\varphi_{0ik}) (-\psi_{abck})\\
&= \sum_{k=0}^6 \varphi_{0ik} \psi_{abck}\\
&= \delta_{0a} \varphi_{ibc}+\delta_{0b} \varphi_{aic}+\delta_{0c} \varphi_{abi}-\delta_{ai} \varphi_{0bc}-\delta_{bi} \varphi_{a0c}-\delta_{ci} \varphi_{ab0}\\
&= \delta_{ia} \omega_{bc}+\delta_{ib} \omega_{ca}+\delta_{ic} \omega_{ab}.
\end{align*}
Contracting~\eqref{eq:pi9} on $i,c$ yields~\eqref{eq:pi10}:
\begin{align*}
\omega_{ik} (\star \omega)_{abik} &=\omega_{ba}+\omega_{ba}+6 \omega_{ab}\\
&= 4 \omega_{ab}.
\end{align*}
For the next identity, there are two ways of computing the desired contractions yielding two different expressions~\eqref{eq:pi11} and~\eqref{eq:pi12}. First, we use the usual way:
\begin{align*}
\psi^{+}_{ijk} (\star \omega)_{abck} &= \sum_{k=1}^{6} \varphi_{ijk}(-\psi_{abck})\\
&= -\sum_{k=0}^{6} \varphi_{ijk}\psi_{abck} + \varphi_{ij0}\psi_{abc0}\\
&= -(\delta_{ia} \varphi_{jbc}+\delta_{ib} \varphi_{ajc}+\delta_{ic} \varphi_{abj}-\delta_{aj} \varphi_{ibc}-\delta_{bj} \varphi_{aic}-\delta_{cj} \varphi_{abi}) + (-\omega_{ij}) \psi^{-}_{abc}\\
&= - \delta_{ia} \psi^{+}_{jbc}-\delta_{ib} \psi^{+}_{ajc}-\delta_{ic} \psi^{+}_{abj}+\delta_{aj} \psi^{+}_{ibc}+\delta_{bj}\psi^{+}_{aic}+\delta_{cj} \psi^{+}_{abi} - \omega_{ij} \psi^{-}_{abc}.
\end{align*}
Second, we can also use the previous results to get:
\begin{align*}
\psi^{+}_{ijk} (\star \omega)_{abck} &= \psi^{-}_{iju}  \omega_{ku} (\star \omega)_{abck}\\
&= -\psi^{-}_{iju}  \omega_{uk} (\star \omega)_{abck}\\
&= -\psi^{-}_{iju} (\delta_{au} \omega_{bc}+\delta_{bu} \omega_{ca}+\delta_{cu} \omega_{ab})\\
&= -\psi^{-}_{ija} \omega_{bc}-\psi^{-}_{ijb} \omega_{ca}-\psi^{-}_{ijc} \omega_{ab}.
\end{align*}
Note that both contractions of~\eqref{eq:pi11} and~\eqref{eq:pi12} on $j,c$ yiled the same result~\eqref{eq:pi13}:
\begin{align*}
\psi^{+}_{ijk} (\star \omega)_{abjk} &= -\psi^{+}_{abi}+ \psi^{+}_{iba}+ \psi^{+}_{aib}+6 \psi^{+}_{abi} - \omega_{ij} \psi^{-}_{abj}\\
&= - \psi^{+}_{iab}- \psi^{+}_{iab}- \psi^{+}_{iab}+6 \psi^{+}_{iab} -  \psi^{+}_{abi}\\
&=2 \psi^{+}_{iab},
\end{align*}
and
\begin{align*}
\psi^{+}_{ijk} (\star \omega)_{abjk} &= -\psi^{-}_{ija} \omega_{bj}-\psi^{-}_{ijb} \omega_{ja}\\
&= \psi^{-}_{iaj} \omega_{bj}-\psi^{-}_{ibj} \omega_{aj}\\
&= \psi^{+}_{iab} -\psi^{+}_{iba}\\
&=2 \psi^{+}_{iab}.
\end{align*}
Since the second way of computing the contraction of $\psi^{+}$ and $\star \omega$ gave us a nicer expression, we use it again for~\eqref{eq:pi14}:
\begin{align*}
\psi^{-}_{ijk} (\star \omega)_{abck} &= -\psi^{+}_{iju}  \omega_{ku} (\star \omega)_{abck}\\
&= \psi^{+}_{iju}  \omega_{uk} (\star \omega)_{abck}\\
&= \psi^{+}_{iju} (\delta_{au} \omega_{bc}+\delta_{bu} \omega_{ca}+\delta_{cu} \omega_{ab})\\
&= \psi^{+}_{ija} \omega_{bc}+\psi^{+}_{ijb} \omega_{ca}+\psi^{+}_{ijc} \omega_{ab}.
\end{align*}
Contracting~\eqref{eq:pi14} on $j,c$ yields~\eqref{eq:pi15}:
\begin{align*}
\psi^{-}_{ijk} (\star \omega)_{abjk} &=  \psi^{+}_{ija} \omega_{bj}+\psi^{+}_{ijb} \omega_{ja}\\
&=  -\psi^{+}_{iaj} \omega_{bj}+\psi^{+}_{ibj} \omega_{aj}\\
&= \psi^{-}_{iab} - \psi^{-}_{iba}\\
&=2 \psi^{-}_{iab}.
\end{align*}
Finally, we compute~\eqref{eq:pi16}:
\begin{align*}
(\star \omega)_{ijkl}  (\star \omega)_{abkl} =& \sum_{k,l=1}^6 \psi_{ijkl}\psi_{abkl}\\
=& \sum_{k,l=0}^6 \psi_{ijkl}\psi_{abkl}-\sum_{k=0}^6 \psi_{ijk0}\psi_{abk0}-\sum_{l=0}^6 \psi_{ij0l}\psi_{ab0l}\\
=& \sum_{k,l=0}^6 \psi_{ijkl}\psi_{abkl}-2\sum_{k=1}^6 \psi^{-}_{ijk}\psi^{-}_{abk}\\
=& (4 \delta_{ia}\delta_{jb} - 4 \delta_{ib}\delta_{ja} - 2\psi_{ijab})-2( \delta_{ia} \delta_{jb}-\delta_{ib} \delta_{ja} - \omega_{ia} \omega_{jb} +\omega_{ib} \omega_{ja})\\
=& 2 \delta_{ia}\delta_{jb} - 2 \delta_{ib}\delta_{ja}+2 ((\star \omega)_{ijab} + \omega_{ia}\omega_{jb} + \omega_{aj}\omega_{ib})\\
=& 2 \delta_{ia}\delta_{jb} - 2 \delta_{ib}\delta_{ja}+2 ((\omega_{ij} \omega_{ab} + \omega_{ja} \omega_{ib} +\omega_{bj} \omega_{ia}) + \omega_{ia}\omega_{jb} + \omega_{aj}\omega_{ib})\\
=& 2 \delta_{ia}\delta_{jb} - 2 \delta_{ib}\delta_{ja}+2 \omega_{ij} \omega_{ab}.
\end{align*}
Contracting~\eqref{eq:pi16} on $b,j$ gives us~\eqref{eq:pi17}:
\begin{align*}
(\star \omega)_{ijkl}  (\star \omega)_{ajkl} = 12 \delta_{ia} - 2 \delta_{ia}+2 \delta_{ia}= 12 \delta_{ia}.\tab\qedhere
\end{align*}
\end{proof}
\end{prop}

\begin{rmk}\label{nkdecomp}
We have the following descriptions of the orthogonal decompositions of $\Omega^2$ and $\Omega^3$ into irreducible subspaces (see~\cite{LF}):\\
\begin{align*}
\Omega^2 &= \Omega^2_1 \oplus \Omega^2_6 \oplus \Omega^2_8,\\
 \Omega^3 &= \Omega^3_{1\oplus 1} \oplus \Omega^3_6 \oplus \Omega^3_{12},
\end{align*}
where the indices denote the corresponding dimensions. In particular:
\begin{itemize}
\item $\Omega^2_1 = \{\beta \in \Omega^2: \star(\beta \wedge \omega) = 2\beta \} =\mathbb{R} \omega$,
\item $\Omega^2_6 =  \{\beta \in \Omega^2: \star(\beta \wedge \omega) = \beta \} = \{ X \intprod \psi^{+}: X \in \Gamma(TM) \}$,
\item $\Omega^2_8= \{\beta \in \Omega^2: \star(\beta \wedge \omega) = -\beta \} \text{ is the space of primitive forms of type $(1,1)$},$
\end{itemize}
\tab or equivalently, $\beta \in \Omega^2_8$ iff $\beta_{ij}\psi^{+}_{ijk} = 0$ and $\beta_{ij} w_{ij} = 0,$
\begin{itemize}
\item $\Omega^3_{1\oplus 1} = \mathbb{R} \psi^{+} \bigoplus \mathbb{R} \psi^{-}$,
\item $\Omega^3_6 = \{X \wedge \omega: X \in \Gamma(TM)\},$
\item $\Omega^3_{12} \text{ is the space of primitive forms of type $(1,2)+(2,1)$, or equivalently, $\Omega^3_{12} = \mathcal{S}^2_{-} \diamond \psi^{+}$,}$
\end{itemize}
where the $\mathcal{S}^2_{-}$ is defined in Section~\ref{diamsection}.

\end{rmk}

\begin{rmk}\label{mapP}
Consider the map $\mathcal{P}:\Omega^2 \rightarrow \Omega^2$ given by $\mathcal{P}(\beta) = \star(\beta \wedge \omega),$ for $\beta \in \Omega^2.$\\
In a local orthonormal frame: $$(\mathcal{P} \beta)_{ab} = \frac{1}{2} \beta_{ij} (\frac{\omega^2}{2})_{ijab}=\frac{1}{2} \beta_{ij} (\star \omega)_{ijab}.$$
We can extend the map $\mathcal{P}$ to all of $\mathcal{T}^2$ via the formula above. Then we have $\mathcal{S}^2 = \operatorname{ker}(\mathcal{P})$ and for $\beta \in \Omega^2$, Remark~\ref{nkdecomp} says that:
\begin{equation}\label{pmapprop}
\mathcal{P}(\beta) = 2\pi_1(\beta)+\pi_6(\beta)-\pi_8(\beta).
\end{equation}
\end{rmk}

\begin{prop}\label{types}
Let $\beta = \beta_0 + \lambda \omega + X \intprod \psi^{+}$, where $\beta_0 \in \Omega^2_8$. Then:
\begin{itemize}
\item $\lambda = \frac{1}{6} \beta_{ij}\omega_{ij}.$
\item $X_k = \frac{1}{4} \beta_{ij} \psi^{+}_{kij}.$
\end{itemize}
\begin{proof}
Recall that $(\beta_0)_{ij}\psi^{+}_{ijk} = 0,(\beta_0)_{ij} \omega_{ij} =0$, and $\omega_{ij} \psi^{+}_{ijk} = 0.$\\
For the first identity, contract 
\begin{equation}\label{eq:req1}
\beta_{ij} =(\beta_0)_{ij} + \lambda \omega_{ij} +  X_a \psi^{+}_{aij}
\end{equation} with $\omega_{ij}$ to get:
\begin{align*}
\beta_{ij} \omega_{ij} &= \lambda \omega_{ij} \omega_{ij}= 6 \lambda.
\end{align*}
Similarly, contracting~\eqref{eq:req1} with $\psi^{+}_{kij}$ gives us:
\begin{align*}
\beta_{ij} \psi^{+}_{kij} &=  X_a \psi^{+}_{aij} \psi^{+}_{kij}= 4 X_a \delta_{ak}= 4 X_k.
\end{align*}
as claimed.
\end{proof}
\end{prop}

\begin{lemma}\label{28omegacommute}
Let $\beta \in \Omega^2_8$. Then $\beta \omega = \omega \beta,$ where by $\beta \omega \in \mathcal{T}^2$ we mean $(\beta \omega)_{ij} = \beta_{ik} \omega_{kl},$ and similarly for $\omega \beta$.
\begin{proof}
Since, $\beta \in \Omega^2_8$, by Remark~\ref{mapP}, we have that $\mathcal{P} \beta = -\beta$, which in a local orthonormal frame is $\frac{1}{2} \beta_{ij} (\star \omega)_{ijab} = -\beta_{ab}.$ Also, recall that $\beta_{ij} \omega_{ij} = 0$. Using Proposition~\ref{nkup}, we have:
\begin{align*}
(\beta \omega)_{st} &= \beta_{su} \omega_{ut}\\
& = -\frac{1}{2} \beta_{ij} (\star \omega)_{ijsu}\omega_{ut}\\
&=  \frac{1}{2} \beta_{ij} \omega_{tu}(\star \omega)_{ijsu}\\
&=  \frac{1}{2} \beta_{ij} (\delta_{it} \omega_{js}+\delta_{jt} \omega_{si}+\delta_{st} \omega_{ij})\\
&=  \beta_{ij} \delta_{it} \omega_{js}+0\\
&=  \beta_{tj} \omega_{js}\\
&=  \omega_{sj} \beta_{jt}\\
&= (\omega \beta)_{st},
\end{align*}
as claimed.
\end{proof}
\end{lemma}

\subsection{The $\diamond$ operator}\label{diamsection}
The results in this section are very similar to the ones in Remark~\ref{g2description}. We describe the isomorphisms coming from the $\diamond$ map. This time, however, we give most of the details.\\
Recall the definition of the $\diamond$ map: let $\sigma \in \Omega^k.$ For $h \in \mathcal{T}^2$, we define:
$$(h\diamond \sigma)_{i_1 \cdots i_k} \coloneqq h_{i_1 p}\sigma_{p i_2 \cdots i_k}+ h_{i_2 p}\sigma_{i_1 p i_3 \cdots i_k} + \cdots + h_{i_k p}\sigma_{i_1 \cdots i_{k-1} p}.$$

\begin{defn}\label{defhatnk}
Let $\beta$ be a $2,3$, or $4$-form. Then we define $\hat{\beta} \in \mathcal{T}^2$ as follows:
\begin{center}
\begin{align*}
\text{for } \beta \in \Omega^2,~~\hat{\beta}_{ia} &\coloneqq \beta_{ik} \omega_{ak},\\
\text{for } \beta \in \Omega^3,~~\hat{\beta}_{ia} &\coloneqq \beta_{ijk} \psi^{+}_{ajk},\\
\text{for } \beta \in \Omega^4,~~\hat{\beta}_{ia} &\coloneqq \beta_{ijkl} (\star \omega)_{ajkl}.\qedhere
\end{align*}
\end{center}
\end{defn}

\begin{rmk}\label{nkdecomofspaces}
Let $$\mathcal{S}^2_{+} \coloneqq \{h \in \mathcal{S}^2| h\omega - \omega h = 0\}$$
$$\mathcal{S}^2_{-} \coloneqq \{h \in \mathcal{S}^2| h\omega + \omega h = 0\}$$
which are the spaces of symmetric $2$-tensors which commute and anti-commute with $\omega$ (or equivalently with $J$), respectively.\\
Note that $\mathcal{S}^2_{-} \subset \mathcal{S}_0$. This can be easily seen by recalling that $\omega^2 = -\operatorname{Id}$, and so $$\trr(h) = - \trr((h \omega) \omega) = -\trr ((-\omega h) \omega) =\trr (\omega h \omega)= \trr (h) \omega^2 = -\trr(h).$$\\
Hence, we can further decompose $$\mathcal{S}^2_{+} = \mathbb{R} g \oplus \mathcal{S}^2_{+0}$$
where $\mathcal{S}^2_{+0}$ are the traceless elements of $\mathcal{S}^2_{+}.$\\
Concluding, we have the orthogonal decomposition: $$\mathcal{S}^2 =\mathbb{R} g \oplus \mathcal{S}^2_{+0}  \oplus  \mathcal{S}^2_{-}.$$
It is easy to check that $\mathcal{S}^2_{+0}$ has dimension $8$ and $\mathcal{S}^2_{-}$ has dimension $12$.
\end{rmk}

\subsubsection{$2$-forms}

\begin{prop}\label{h2forms}
In the case of $2$-forms, the $\cdot \diamond \omega$ map is an isomorphism of the following spaces:
\begin{center}
\begin{align*}
\mathbb{R} g \cong& \Omega^2_1,\\
\Omega^2_6 \cong& \Omega^2_6,\\
\mathcal{S}^2_{+0} \cong& \Omega^2_8.\\
\end{align*}
\end{center}
\begin{proof}
For the first two maps it will be clear that they are isomorphisms. For the last one, we just check that the image under the $\cdot \diamond \omega$ map lies in the required subspace. Then by the next Proposition~\ref{diam2}, which shows that the map is invertible, we conclude that it is also an isomorphism.
So, we have: $$(g \diamond \omega)_{ij} = g_{ip} \omega_{pj} + g_{jp} \omega_{ip} = 2 \omega_{ij}.$$
Next, take any $\beta = X \intprod \psi^{+} \in \Omega^2_6.$ By Proposition~\ref{nkup} and~\ref{types}, we have:
\begin{align*}
(\beta \diamond \omega )_{ij} &= \beta_{ip} \omega_{pj} + \beta_{jp} \omega_{ip}\\
&= X_a \psi^{+}_{aip} \omega_{pj} + X_a \psi^{+}_{ajp} \omega_{ip}\\
&= -X_a \psi^{+}_{aip} \omega_{jp} + X_a \psi^{+}_{ajp} \omega_{ip}\\
&= X_a \psi^{-}_{aij} - X_a \psi^{-}_{aji}\\
&= 2X_a \psi^{-}_{aij}\\
&= 2X_a \psi^{-}_{ija}\\
&= -2X_a \omega_{ap} \psi^{+}_{ijp}\\
&= -(2J(X) \intprod \psi^{+})_{ij},
\end{align*}
where we have used that $$(J(X))_p = \langle J(X), e_p\rangle = \omega(X,e_p) = X_a \omega_{ap}.$$
Finally, take any $h \in S^2_{+0}.$ Then:
\begin{align*}
(h \diamond \omega)_{ij} \omega_{ij} &= (h_{ip} \omega_{pj} + h_{jp} \omega_{ip}) \omega_{ij}\\
&= h_{ip} \delta_{pi} + h_{jp} \delta_{pj}\\
&= 2 \trr{h}\\
&= 0,
\end{align*}
and
\begin{align*}
(h \diamond \omega)_{ij} \psi^{+}_{ijk} &= (h_{ip} \omega_{pj} + h_{jp} \omega_{ip}) \psi^{+}_{ijk}\\
&=  (( h \omega)_{ij} + (\omega h)_{ij}) \psi^{+}_{ijk}\\
&=   2 ( h \omega)_{ij} \psi^{+}_{ijk} \text{ (since $h\omega = \omega h$ for $h \in \mathcal{S}^2_+$)}\\
&=   2h_{ip} \omega_{pj} \psi^{+}_{ijk}\\
&=   2h_{ip} \omega_{pj} \psi^{+}_{kij}\\
&=   -2h_{ip} \psi^{-}_{kip}\\
&=  0,
\end{align*}
because $h \in \mathcal{S}^2.$ This shows that $h \diamond \omega \in \Omega^2_8$. Hence, the result follows.
\end{proof}
\end{prop}

\begin{prop}\label{diam2}
Let $\beta \in \Omega^2$. Then $\beta = h \diamond \omega$ for some unique $h = \frac{1}{6} \trr(h) g + X \intprod \psi^{+} + h_{+0}$, where $X \in \Gamma(TM),$ $h_{+0} \in \mathcal{S}^2_{+0}$. Then $h = \frac{1}{2} \hat{\beta},$ where $\hat{\beta}$ is as in Definition~\ref{defhatnk}. This implies that $$\trr(h) = \frac{1}{2}\trr(\hat{\beta}),\tab X_k = \frac{1}{8} \hat{\beta}_{ia} \psi^{+}_{kia},\tab h_{+0} = \frac{1}{2} \hat{\beta}_{symm} - \frac{1}{12} \trr(\hat{\beta}) g.$$
Also, clearly $\beta = 0$ iff $h =0$ iff $\hat{\beta} = 0.$
\begin{proof}
We compute $\hat{\beta}$ as follows:
\begin{align*}
\hat{\beta}_{ia} &= \beta_{ik}\omega_{ak}\\
&= (h \diamond \omega)_{ik}\omega_{ak}\\
&= (h_{ip} \omega_{pk}+ h_{kp} \omega_{ip})\omega_{ak}\\
&= h_{ip} \delta_{ap} + h_{kp} \omega_{ip}\omega_{ak}\\
&= h_{ia} + (\frac{1}{6} \trr(h) \delta_{kp}+ X_u \psi^{+}_{ukp}+(h_{+0})_{kp}) \omega_{ip} \omega_{ak}\\
&= h_{ia} + \frac{1}{6} \trr(h) \delta_{ia} +X_u  \psi^{+}_{ukp} \omega_{ip} \omega_{ak} - \omega_{ip}(h_{+0})_{pk}\omega_{ka} \\
&= h_{ia} + \frac{1}{6} \trr(h) \delta_{ia} -X_u  \psi^{-}_{uki} \omega_{ak} - (\omega h_{+0} \omega)_{ia}\\
&= h_{ia} + \frac{1}{6} \trr(h) \delta_{ia} +X_u  \psi^{-}_{uik} \omega_{ak} - (h_{+0} \omega^2)_{ia}\\
&= h_{ia} + \frac{1}{6} \trr(h) \delta_{ia} +X_u  \psi^{+}_{uia} + (h_{+0})_{ia}\\
&= 2h_{ia}.
\end{align*}
as claimed.\\
So, we get $2 \trr(h) = \trr(\hat{\beta})$ along with $$\hat{\beta}_{symm} = 2 (\frac{1}{6} \trr(h) g + h_{+0})$$ which means that
$$ h_{+0} =\frac{1}{2} \hat{\beta}_{symm} - \frac{1}{6} \trr(h) g =\frac{1}{2} \hat{\beta}_{symm} - \frac{1}{12} \trr(\hat{\beta}) g.$$
Finally, by Proposition~\ref{types}, we get that $$2X_k = \frac{1}{4} (\hat{\beta}_{skew})_{ia} \psi^{+}_{kia} = \frac{1}{4} \hat{\beta}_{ia} \psi^{+}_{kia}.\qedhere$$
\end{proof}
\end{prop}

\subsubsection{$4$-forms}

\begin{prop}
Let $\beta \in \Omega^4$. Then $\beta = h \diamond (\star \omega)$ for some unique $h \in \Omega^2_6 \oplus \mathcal{S}^2_{+}.$
\begin{proof}
It is easy to check that: $$\star(h \diamond \omega) = (\frac{1}{4} \trr(h) g - h^{T}) \diamond (\star \omega).$$
Now, since $\beta \in \Omega^4$, we have $\star \beta \in \Omega^2$. Then by Proposition~\ref{h2forms}, $\star \beta = h \diamond \omega$, for some unique $h \in\Omega^2_6 \oplus \mathcal{S}^2_{+}.$\\
Hence, $$\beta = \star(\star \beta) =  (\frac{1}{4} \trr(h) g - h^{T}) \diamond (\star \omega).$$
Note that the map $h \mapsto  (\frac{1}{4} \trr(h) g - h^{T})$ is an automorphism of $\Lambda^2_6 \oplus \mathcal{S}^2_{+}$. This is because it can be seen that under this map, $\Omega^2_6$ is mapped to itself and for $h \in \mathcal{S}^2_{+}$ we have:
$$h \mapsto \frac{1}{4}\trr(h) g - h,$$
which is injective, as $\frac{1}{4}\trr(h) g - h = 0$ iff $h = cg$, for some $c\in \mathbb{R}$, but then $\frac{1}{4} 6c g = cg$, hence $c=0$.
Also, since $g,h$ commute with $\omega$, $\frac{1}{4}\trr(h) g - h$ also commutes with $\omega$, so is in $\mathcal{S}^2_{+}.$ Thus, we get the required result.
\end{proof}
\end{prop}

\begin{prop}\label{prop4forms1}
Let $\beta \in \Omega^4$, so $\beta = h \diamond (\star \omega)$ for some unique $h \in \Omega^2_6 \oplus \mathcal{S}^2_{+}$. Then $$\hat{\beta} = 8\trr(h)g + 12h_6 + 12h_{+0},$$
where $\hat{\beta}$ is as in Definition~\ref{defhatnk}.
\begin{proof}
We compute $\hat{\beta}$:
\begin{align*}
\hat{\beta}_{ia} &= \beta_{ijkl} (\star \omega)_{ajkl}\\
&= (h_{ip}  (\star \omega)_{pjkl}+h_{jp}  (\star \omega)_{ipkl}+h_{kp}  (\star \omega)_{ijpl}+h_{lp}  (\star \omega)_{ijkp}) (\star \omega)_{ajkl}\\
&= h_{ip}  (\star \omega)_{pjkl}  (\star \omega)_{ajkl} +3h_{jp}  (\star \omega)_{ipkl} (\star \omega)_{ajkl}\\
&=h_{ip} 12 \delta_{pa} + 3 h_{jp} (2 \delta_{ia} \delta_{pj}-2 \delta_{ij} \delta_{pa} + 2 \omega_{ip} \omega_{aj})\\
&=12h_{ia}+ 6 \trr(h)\delta_{ia}-6 h_{ia} + 6 h_{jp} \omega_{ip} \omega_{aj}
\end{align*}
In the proof of Proposition~\ref{diam2}, we computed that for $h \in \Omega^2_6 \oplus \mathcal{S}^2_{+}$, $h_{kp} \omega_{ip}\omega_{ak} = h_{ia}$. Hence,
\begin{align*}
\hat{\beta}_{ia} &= 6h_{ia}+6 \trr(h)\delta_{ia}+6 h_{ia}\\
&= 12h_{ia}+6 \trr(h)\delta_{ia}\\
 &= 12(\frac{1}{6} \trr(h) g + h_6 + h_{+0})_{ia})+6 \trr(h)\delta_{ia}\\
 &=8 \trr(h)\delta_{ia} + 12 (h_6)_{ia} + 12(h_{+0})_{ia}.\qedhere
\end{align*}
\end{proof}
\end{prop}

\begin{cor}\label{cor4diam}
Let $\beta \in \Omega^4$, so $\beta = h \diamond (\star \omega)$ for some unique $h = \frac{1}{6} \trr(h)g+ h_6 + h_{+0}$. Then:
\begin{center}
\begin{align*}
\trr(h) =& \frac{1}{48} \trr(\hat{\beta}),\\
(h_{+0})_{ia} =& \frac{1}{12} (\hat{\beta}_{symm})_{ia} - \frac{1}{72} \trr(\hat{\beta}) \delta_{ia} = \frac{1}{24} (\hat{\beta}_{ia}+\hat{\beta}_{ai}) - \frac{1}{72} \trr(\hat{\beta}) \delta_{ia},\\
(h_6)_{ia} =& \frac{1}{12}(\hat{\beta}_{skew})_{ia}=\frac{1}{24} (\hat{\beta}_{ia}-\hat{\beta}_{ai}).
\end{align*}
\end{center}
Also, clearly $\beta = 0$ iff $h =0$ iff $\hat{\beta} = 0.$
\begin{proof}
In Proposition~\ref{prop4forms1} we proved that
\begin{equation}\label{eq:cor4diam1}
\hat{\beta} = 8\trr(h)g + 12h_6 + 12h_{+0}.
\end{equation}
Taking traces of both sides, we get $$\trr(\hat{\beta}) = 48\trr(h).$$
Hence, $$h_{+0} = \frac{1}{12} (\hat{\beta}_{symm} - 8 \trr(h)g) = \frac{1}{12} \hat{\beta}_{symm} - \frac{2}{3} \cdot \frac{1}{48}\trr(\hat{\beta})g = \frac{1}{12}\hat{\beta}_{symm} - \frac{1}{72}\trr(\hat{\beta})g.$$
Taking skew-symmetric parts of ~\eqref{eq:cor4diam1} gives us the required $$12h_6 = \hat{\beta}_{skew}.\qedhere$$
\end{proof}
\end{cor}

\subsubsection{$3$-forms}

\begin{prop}
In the case of the $3$-forms, the $\cdot \diamond \psi^{+}$ map is an isomorphism of the following spaces:
\begin{center}
\begin{align*}
\mathbb{R} g \oplus \mathbb{R} \omega \cong&~ \Omega^3_{1\oplus 1},\\
\Omega^2_6 \cong&~ \Omega^3_6,\\
\mathcal{S}^2_{-} \cong&~ \Omega^3_{12}.\\
\end{align*}
\end{center}
\begin{proof}
Computing $g \diamond \psi^{+}$ and $\omega \diamond \psi^{+}$ gives us:
\begin{align*}
g \diamond \psi^{+} &= 3 \psi^{+}.\\
(\omega \diamond \psi^{+})_{ijk} &= \omega_{ip} \psi^{+}_{pjk}+\omega_{jp} \psi^{+}_{ipk}+\omega_{kp} \psi^{+}_{ijp}\\
&= \omega_{ip} \psi^{+}_{jkp}+\omega_{jp} \psi^{+}_{kip}+\omega_{kp} \psi^{+}_{ijp}\\
&= - \psi^{-}_{jki}- \psi^{-}_{kij}- \psi^{-}_{ijk}\\
&= -3 \psi^{-}_{ijk},
\end{align*}
which is enough to conclude that $\mathbb{R} g \oplus \mathbb{R} \omega \cong \Omega^3_{1\oplus 1}$.\\
Next, take any $X \intprod \psi^{+},$ with $X \in \Gamma(TM)$. Then:
\begin{align*}
((X \intprod \psi^{+}) \diamond \psi^{+})_{ijk} =& (X \intprod \psi^{+})_{ip} \psi^{+}_{pjk} + (X \intprod \psi^{+})_{jp} \psi^{+}_{ipk} + (X \intprod \psi^{+})_{kp} \psi^{+}_{ijp}\\
=& X_u \psi^{+}_{uip} \psi^{+}_{pjk} + X_u \psi^{+}_{ujp} \psi^{+}_{ipk} + X_u \psi^{+}_{ukp} \psi^{+}_{ijp}\\
=& X_u (\psi^{+}_{uip} \psi^{+}_{jkp} + \psi^{+}_{ujp} \psi^{+}_{kip} + \psi^{+}_{ukp} \psi^{+}_{ijp})\\
=& X_u (\delta_{uj} \delta_{ik} - \delta_{uk} \delta_{ij} - \omega_{uj} \omega_{ik} + \omega_{uk} \omega_{ij}\\
&+ \delta_{uk} \delta_{ji} - \delta_{ui} \delta_{jk} - \omega_{uk} \omega_{ji} + \omega_{ui} \omega_{jk}\\
&+ \delta_{ui} \delta_{kj} - \delta_{uj} \delta_{ki} - \omega_{ui} \omega_{kj} + \omega_{uj} \omega_{ki})\\
=& 2 X_u (\omega_{ui} \omega_{jk} + \omega_{uj} \omega_{ki} + \omega_{uk} \omega_{ij})\\
=& 2 (JX \wedge \omega)_{ijk},
\end{align*}
which again is enough to see that $\Omega^2_6 \cong \Omega^3_6$.\\
For the last isomorphism, we avoid the details, because this is how we defined $\Omega^3_{12}$ in Remark~\ref{nkdecomp}.
\end{proof}
\end{prop}

\begin{prop}\label{cor3diam}
Let $\beta \in \Omega^3$, so $\beta = h \diamond \psi^{+}$ for some unique $h \in \mathbb{R} g \oplus \mathbb{R} \omega \oplus \Omega^2_6 \oplus \mathcal{S}^2_{-}$. Then $$\hat{\beta} =  2 \trr(h) g + 12 \lambda \omega + 4 h_6 + 4 h_{-},$$
where $\hat{\beta}$ is as in Definition~\ref{defhatnk}, and $\lambda$ is the coefficient of $\omega$ in $h$, meaning that the unique part of $h$ in $\mathbb{R} \omega$ is $\lambda \omega.$
\begin{proof}
Let $h_6 = X \intprod \psi^{+}$, for some unique $X \in \Gamma(TM).$ Now, we just compute $\hat{\beta}:$
\begin{align*}
\hat{\beta}_{ia} &= \beta_{ijk} \psi^{+}_{ajk}\\
 &= (h \diamond \psi^{+}) \psi^{+}_{ajk}\\
 &= (h_{ip} \psi^{+}_{pjk} + h_{jp} \psi^{+}_{ipk}+h_{kp} \psi^{+}_{ijp}) \psi^{+}_{ajk}\\
 &= h_{ip} \psi^{+}_{pjk} \psi^{+}_{ajk}+2  h_{jp} \psi^{+}_{ipk} \psi^{+}_{ajk}\\
 &= h_{ip} 4 \delta_{pa} +2  h_{jp} (\delta_{ia} \delta_{pj} - \delta_{ij} \delta_{pa} - \omega_{ia}\omega_{pj} + \omega_{ij}\omega_{pa})\\
 &= 4 h_{ia} + 2 \trr(h) \delta_{ia} -2 h_{ia} + 2 h_{jp} \omega_{jp}\omega_{ia} + 2 h_{jp} \omega_{pa}\omega_{ij}.
\end{align*}
We compute the last two terms separately:
\begin{align*}
2 h_{jp} \omega_{pa}\omega_{ij} &= 2 ( \frac{1}{6} \trr(h) \delta_{jp} + \lambda \omega_{jp} + X_u \psi^{+}_{ujp} + (h_{-})_{jp})\omega_{pa}\omega_{ij}\\
&=  \frac{1}{3} \trr(h) \omega_{pa}\omega_{ip} + 2\lambda  \omega_{jp} \omega_{pa}\omega_{ij} - 2X_u \psi^{+}_{ujp}\omega_{ap} \omega_{ij} + 2\omega_{ij}  (h_{-})_{jp}\omega_{pa}\\
&=  -\frac{1}{3} \trr(h) \delta_{ia} - 2\lambda \delta_{aj} \omega_{ij} +2 X_u \psi^{-}_{uja} \omega_{ij} + 2(\omega h_{-} \omega)_{ia}\\
&=  -\frac{1}{3} \trr(h) \delta_{ia} - 2\lambda \omega_{ia} - 2X_u \psi^{-}_{uaj} \omega_{ij} - 2(h_{-} \omega^2)_{ia}\\
&=  -\frac{1}{3} \trr(h) \delta_{ia} - 2\lambda \omega_{ia} - 2X_u \psi^{+}_{uai} + 2(h_{-})_{ia}\\
&=  -\frac{1}{3} \trr(h) \delta_{ia} - 2\lambda \omega_{ia} + 2(h_6)_{ia} + 2(h_{-})_{ia}.
\end{align*}
Next,
\begin{align*}
2 h_{jp} \omega_{jp}\omega_{ia} &= 2 ( \frac{1}{6} \trr(h) \delta_{jp} + \lambda \omega_{jp} + X_u \psi^{+}_{ujp} + (h_{-})_{jp})\omega_{jp}\omega_{ia}\\
&= 0 + 12 \lambda \omega_{ia} +0 + 0\\
&= 12 \lambda \omega_{ia}.\\
\end{align*}
Hence, combining these parts we get:
\begin{align*}
\hat{\beta}_{ia} &= 2 h_{ia} + 2 \trr(h) \delta_{ia} +(-\frac{1}{3} \trr(h) \delta_{ia} - 2\lambda \omega_{ia} + 2(h_6)_{ia} + 2(h_{-})_{ia}) + 12 \lambda \omega_{ia}\\
&= 2(\frac{1}{6} \trr(h) \delta_{ia} + \lambda \omega_{ia} + (h_6)_{ia} + (h_{-})_{ia})+ \frac{5}{3} \trr(h) \delta_{ia} + 10\lambda \omega_{ia} + 2(h_6)_{ia} + 2(h_{-})_{ia}\\
&= 2 \trr(h) \delta_{ia} + 12\lambda \omega_{ia} + 4(h_6)_{ia} +4 (h_{-})_{ia},
\end{align*}
as claimed.
\end{proof}
\end{prop}

\begin{cor}\label{diamcorol3}
Let $\beta \in \Omega^3$, so $\beta = h \diamond \psi^{+}$ for some unique $h = \frac{1}{6} \trr(h) g + \lambda \omega+ X \intprod \psi^{+} + h_{-}$, where $X \in \Gamma(TM)$. Then:
\begin{center}
\begin{align*}
\trr(h) =& \frac{1}{12} \trr(\hat{\beta}),\\
(h_{-})_{ia} =& \frac{1}{4} (\hat{\beta}_{symm})_{ia} - \frac{1}{24} \trr(\hat{\beta}) \delta_{ia} = \frac{1}{8} (\hat{\beta}_{ia}+\hat{\beta}_{ai}) - \frac{1}{24} \trr(\hat{\beta}) \delta_{ia},\\
\lambda =& \frac{1}{72} \hat{\beta}_{ia} \omega_{ia},\\
X_k =& \frac{1}{16}  \hat{\beta}_{ia} \psi^{+}_{kia}.
\end{align*}
\end{center}
Also, clearly $\beta = 0$ iff $h =0$ iff $\hat{\beta} = 0.$
\begin{proof}
In Proposition~\ref{cor3diam} we proved that:
\begin{equation}\label{eq:cor3diam1}
\hat{\beta} =  2 \trr(h) g + 12 \lambda \omega + 4 h_6 + 4 h_{-}.
\end{equation}
Taking traces of both sides yields $$\trr(\hat{\beta}) = 12 \trr(h).$$
Next, taking symmetric parts of both sides of~\eqref{eq:cor3diam1} gives us: $$\hat{\beta}_{symm} =  2 \trr(h) g+ 4 h_{-}.$$ Thus, $$h_{-} = \frac{1}{4} (\hat{\beta}_{symm} - 2 \trr(h) g) = \frac{1}{4} \hat{\beta}_{symm} - \frac{1}{24} \trr(\hat{\beta}) g.$$
On the other hand, comparing skew-symmetric parts of both sides of~\eqref{eq:cor3diam1} gives us:
$$\hat{\beta}_{skew} = 12 \lambda \omega+ 4 h_6$$
We recall Proposition~\ref{types} to get
$$12 \lambda = \frac{1}{6}(\hat{\beta}_{skew})_{ia}  \omega_{ia} =   \frac{1}{6}\hat{\beta}_{ia}  \omega_{ia}$$
and
$$4 X_k= \frac{1}{4}(\hat{\beta}_{skew})_{ia}  \psi^{+}_{kia} =   \frac{1}{4}\hat{\beta}_{ia}  \psi^{+}_{kia},$$
which concludes the proof.
\end{proof}
\end{cor}

\subsection{Nearly K\"{a}hler $6$-manifolds}\label{nkmanifoldsintro}
 
Let $(M^{6},g,J,\Omega)$ be a compact connected $6$-manifold with an $SU(3)$-structure. We say it is \textbf{nearly K\"{a}hler} if:
\begin{equation}\label{defnk1}
\nabla_X \omega = -X \intprod \psi^{+} \text{ and } \nabla_X \psi^{+} = X \wedge \omega.
\end{equation}
In dimension $6$ it is equivalent to $(\nabla_X J)(X) = 0$, for all $X \in \Gamma(TM)$, but $\nabla J \neq 0.$ 
Also, by~\cite{GR} it is also equivalent to $d\omega = 3 \nabla \omega$ or that $d \omega = - 3\psi^{+}$ and $d\psi^{-} = 4 \frac{\omega^2}{2}.$ Moreover, one can check that in this case the conical $G_2$ structure on $M \times \mathbb{R}$ is torsion-free. Finally, it is a fact that all nearly K\"{a}hler manifolds in dimension $6$ are positive Einstein. With our choice of normalization, the Einstein constant is $5$.\\
In a local orthonormal frame we can write~\eqref{defnk1} as:
\begin{equation}\label{defnk2}
\nabla_i \omega_{jk} = -\psi^{+}_{ijk}\text{ and } \nabla_i \psi^{+}_{jkl} = \delta_{ij} \omega_{kl} + \delta_{ik} \omega_{lj}+ \delta_{il} \omega_{jk}
\end{equation}
Note that contracting the second identity on $i,j$ gives us
\begin{equation}\label{defnk3}
\nabla_i \psi^{+}_{ikl} =6 \omega_{kl} + \omega_{lk}+ \omega_{lk} = 4 \omega_{kl}.
\end{equation}

\subsection{Curvature identities}\label{curvvv}
On a nearly K\"{a}hler manifold we have the Einstein constant $k=5$. Applying the result from Lemma~\ref{hatcirclemma} we get:
\begin{equation}\label{curidnkah}
\begin{split}
\hat{W} &= \hat{R} + 2 \operatorname{Id},\\
\mathring{W} &= \mathring{R} - \operatorname{Id}, \text{ on $\mathcal{S}^2_0.$}
\end{split}
\end{equation}

\begin{prop}\label{curiden}
The following identities hold:
\begin{itemize}
\itemequation[eq:cur1] {$R_{pqiu} \psi^{+}_{liu} = -2 \psi^{+}_{pql}.$}{}
\itemequation[eq:cur2] {$R_{pqiu} \psi^{-}_{viu} = -2 \psi^{-}_{pqv}.$}{}
\itemequation[eq:cur3] {$R_{pqju} \omega_{ju} = -2 \omega_{pq}.$}{}
\end{itemize}
\begin{proof}
For the first identity~\eqref{eq:cur1}, we show that computing contraction of $\psi^{-}$ and $\nabla \nabla \omega$ yields the required result. Explicitly,
\begin{align*}
\nabla_p \nabla_q \omega_{ij} &= \nabla_p (-\psi^{+}_{qij})\\
&= - (\delta_{pq} \omega_{ij}+\delta_{pi} \omega_{jq}+\delta_{pj} \omega_{qi}) \text{ (by~\eqref{defnk2})}
\end{align*}
Now, we use the Ricci identity to get:
\begin{align*}
-R_{pqiu} \omega_{uj} - R_{pqju} \omega_{iu} &= (\nabla_p \nabla_q - \nabla_q \nabla_p) \omega_{ij}\\
&= -(\delta_{pi} \omega_{jq}+\delta_{pj} \omega_{qi}) + (\delta_{qi} \omega_{jp}+\delta_{qj} \omega_{pi}).
\end{align*}
Contracting both sides with $\psi^{-}_{ijl}$ and using skew-symmetry of both sides in $i,j$ we get:
\begin{align*}
-2R_{pqiu} \omega_{uj} \psi^{-}_{ijl} &= -2 \delta_{pi} \omega_{jq} \psi^{-}_{ijl} + 2 \delta_{qi} \omega_{jp} \psi^{-}_{ijl}\\
2R_{pqiu} \psi^{-}_{ilj} \omega_{uj} &=  -2 \psi^{-}_{pjl} \omega_{jq} + 2 \psi^{-}_{qjl} \omega_{jp}\\
2R_{pqiu} \psi^{+}_{ilu} &=  -2\psi^{-}_{plj} \omega_{qj} +  2\psi^{-}_{qlj} \omega_{pj}\\
-2R_{pqiu} \psi^{+}_{liu} &=  - 2\psi^{+}_{plq} + 2 \psi^{+}_{qlp}\\
&=4 \psi^{+}_{pql},
\end{align*}
which yields~\eqref{eq:cur1}.\\
For~\eqref{eq:cur2}, contract~\eqref{eq:cur1} with $\omega_{vl}$ to get:
\begin{align*}
R_{pqiu} \psi^{+}_{liu} \omega_{vl}&= -2 \psi^{+}_{pql} \omega_{vl}\\
R_{pqiu} \psi^{+}_{iul} \omega_{vl}&= -2 (-\psi^{-}_{pqv})\\
R_{pqiu} (-\psi^{-}_{iuv})&= -2 (-\psi^{-}_{pqv})\\
R_{pqiu} \psi^{-}_{viu}&= -2\psi^{-}_{pqv},
\end{align*}
as desired.\\
Finally, as for the first identity, we first compute $\nabla \nabla \psi^{+}$ and then contract it with $\psi^{-}$. Explicitly,
\begin{align*}
\nabla_p \nabla_q \psi^{+}_{jkl} &= \nabla_p (\delta_{qj} \omega_{kl} + \delta_{qk} \omega_{lj}+\delta_{ql} \omega_{jk})\\
&= -(\delta_{qj} \psi^{+}_{pkl} + \delta_{qk} \psi^{+}_{plj}+\delta_{ql} \psi^{+}_{pjk})  \text{ (by~\eqref{defnk2}).}
\end{align*}
Now we use the Ricci identity to get:
\begin{align*}
-R_{pqju}\psi^{+}_{ukl}-R_{pqku}\psi^{+}_{jul}-R_{pqlu}\psi^{+}_{jku} &= (\nabla_p \nabla_q-\nabla_q \nabla_p) \psi^{+}_{jkl}\\
&= -(\delta_{qj} \psi^{+}_{pkl} + \delta_{qk} \psi^{+}_{plj}+\delta_{ql}\psi^{+}_{pjk})+(\delta_{pj} \psi^{+}_{qkl} + \delta_{pk} \psi^{+}_{qlj}+\delta_{pl} \psi^{+}_{qjk}).
\end{align*}
Contracting both sides with $\psi^{-}_{jkl}$ and using the skew-symmetry in $j,k,l$ we get~\eqref{eq:cur3}:
\begin{align*}
-3R_{pqju}\psi^{+}_{ukl}\psi^{-}_{jkl} &= -3\delta_{qj} \psi^{+}_{pkl} \psi^{-}_{jkl}+ 3\delta_{pj} \psi^{+}_{qkl} \psi^{-}_{jkl}\\
-3R_{pqju} (4\omega_{uj}) &= -3 \psi^{+}_{pkl} \psi^{-}_{qkl}+ 3 \psi^{+}_{qkl} \psi^{-}_{pkl}\\
12R_{pqju} \omega_{ju} &= -12 \omega_{pq}+ 12 \omega_{qp}\\
R_{pqju} \omega_{ju} &= -2 \omega_{pq}. \qedhere
\end{align*}
\end{proof}
\end{prop}
\begin{rmk}\label{spacepreserving}
Proposition~\ref{curiden} says that $\hat{R} = -2 \operatorname{Id}$ on $\psi^{+}_{ijk}, \psi^{+}_{ijk}, \omega_{ij}$.
Recall that by~\eqref{curidnkah}, we have $\hat{W} = \hat{R} + 2 \operatorname{Id}.$ Hence, $\hat{W} \psi^{+}, \hat{W} \psi^{-}, \hat{W} \omega$ are all equal to $0$, which is exactly what is needed in order for $W$ to be in $\Omega^2_8$ (in the first two or the last two indices), by Remark~\ref{nkdecomp}. Hence, we have that $(W \beta)_{ab} = W_{abij} \beta_{ij}$ will always lie in $\Omega^2_8$. Thus, since $\hat{R}$ and $\hat{W}$ dffer by a constant, we can conclude that both $\hat{W}$ and $\hat{R}$ preserve $\Omega^2_8$.\\
We claim that $\mathring{W}$ preserves both $\mathcal{S}^2_{-}$ and $\mathcal{S}^2_{+0}$. For the first subspace, let $h \in \mathcal{S}^2_{-}$. Then we compute:
\begin{align*}
((\mathring{W}h) \omega)_{ab} &= (\mathring{W}h)_{au} \omega_{ub} = W_{kalu} h_{kl} \omega_{ub}\\
&= -(W_{klua}+W_{kual}) h_{kl} \omega_{ub}\\
&=-W_{kual} h_{kl} \omega_{ub}\\
&=-W_{ubal} h_{kl} \omega_{ku}\text{~(because $W \in \Omega^2_8$ in the first (last) two indices)}\\
&=-W_{ubal} \omega_{kl} h_{ku}\text{~(because $h \in \mathcal{S}^2_{-}$)}\\
&=W_{ublk} \omega_{al} h_{ku}\\
&=-W_{klub} h_{ku} \omega_{al}\\
&=- (\mathring{W}h)_{lb} \omega_{al}\\
&= - (\omega (\mathring{W}h))_{ab},
\end{align*}
as claimed. The other case is similar, along with recalling that $W$ is Ricci-traceless. Finally, since $\mathring{W}$ and $\mathring{R}$ differ by a constant on $\mathcal{S}^2_0,$ $\mathring{R}$ also preserves that splitting.\\
These fact mean that we can consider $\hat{W}$ ($\mathring{W}$ resp.) as a self-adjoint operator only on $\Omega^2_{8}$ ($\mathcal{S}^2_{-}$ and $\mathcal{S}^2_{+0}$ resp.) which will provide better estimates when we apply the Bochner-Weitzenböck techniques.
\end{rmk}
\subsection{Harmonic forms}\label{harmshuidfh}
In this section we derive some useful properties about the harmonic forms. We will use the fact that harmonic $2$-forms lie in $\Omega^2_8$ and harmonic $3$-forms lie in $\Omega^3_{12}.$ See~\cite[Theorem 3.8]{LF}.

\begin{defn}
For $h \in \mathcal{S}^2$, let $\tilde{h} \in \mathcal{T}^2$ be defined as $$\tilde{h}_{kc}\coloneqq(\nabla_i h_{jk}) \psi^{+}_{ijc}. \qedhere$$
\end{defn}

\begin{prop}\label{randomidentities}
Let $h \in \mathcal{S}^2_-$. Then:
\begin{itemize}
\item $(\nabla_a h_{ki})\omega_{ak} =  - (\divv h)_k \omega_{ki}.$
\item $(\nabla_u h_{ik}) \psi^{-}_{uib} \omega_{ka} =  \tilde{h}_{ab}+4 (h\omega)_{ab}.$
\end{itemize}
\begin{proof}
Since, $h \in \mathcal{S}^2_-$, we have $$h_{ik}\omega_{ka} + \omega_{ik} h_{ka} = 0.$$ Differentiate it to get:
\begin{align}
0 &= (\nabla_u h_{ik})\omega_{ka} + h_{ik} (\nabla_u \omega_{ka}) +  (\nabla_u\omega_{ik}) h_{ka} +\omega_{ik} (\nabla_u h_{ka})\nonumber\\
&= (\nabla_u h_{ik})\omega_{ka} - h_{ik} \psi^{+}_{uka} - \psi^{+}_{uik} h_{ka} +\omega_{ik} (\nabla_u h_{ka}).\label{eq:randomidentities1}
\end{align}
Contract~\eqref{eq:randomidentities1} on $a,u$ to get:
\begin{align*}
0 &= (\nabla_a h_{ik})\omega_{ka} +\omega_{ik} (\nabla_a h_{ka})\\
&= - (\nabla_a h_{ki})\omega_{ak} +\omega_{ik} (\divv h)_k\\
\end{align*}
which gives the desired $$ (\nabla_a h_{ki})\omega_{ak} =  - (\divv h)_k \omega_{ki}.$$
For the second identity, contract both sides of~\eqref{eq:randomidentities1} with $\psi^{-}_{uib}$ to get:
$$0= (\nabla_u h_{ik})\omega_{ka} \psi^{-}_{uib} - h_{ik} \psi^{+}_{uka}\psi^{-}_{uib} - \psi^{+}_{uik} h_{ka}\psi^{-}_{uib} +\omega_{ik} (\nabla_u h_{ka})\psi^{-}_{uib}$$
The first term is what we need to solve for. So let us simplify the others separately:
\begin{align*}
h_{ik} \psi^{+}_{uka}\psi^{-}_{uib} &= h_{ik} \psi^{+}_{kau} \psi^{-}_{ibu}\\
&= h_{ik} (\delta_{ki}\omega_{ab}+\delta_{ab}\omega_{ki}-\delta_{kb}\omega_{ai}-\delta_{ai}\omega_{kb})\\
&= 0 + 0 - h_{ib}\omega_{ai} - h_{ak}\omega_{kb}\\
&= -(\omega h + h \omega)_{ab}\\
&= 0.
\end{align*}
Similarly, we have:
$$\psi^{+}_{uik} h_{ka}\psi^{-}_{uib} = h_{ka} \psi^{+}_{kui}\psi^{-}_{bui} = h_{ka} 4 \omega_{kb} = 4 (h\omega)_{ab},$$
and
$$\omega_{ik} (\nabla_u h_{ka})\psi^{-}_{uib} = (\nabla_u h_{ka})\psi^{-}_{ubi}\omega_{ki} = (\nabla_u h_{ka}) \psi^{+}_{ubk} =  -(\nabla_u h_{ka}) \psi^{+}_{ukb} = -\tilde{h}_{ab}.$$
Hence, $$ (\nabla_u h_{ik})\omega_{ka} \psi^{-}_{uib} = \tilde{h}_{ab}+4 (h\omega)_{ab}.\qedhere$$
\end{proof}
\end{prop}

\begin{prop}\label{ushdug}
Let $h \in \mathcal{S}^2$, so that $\tilde{h} \in \mathcal{T}^2 = \mathcal{S}^2 \oplus \Omega^2$. Then $\tilde{h}_{skew} \in \Omega^2_8.$
\begin{proof}
By Remark~\ref{mapP}, it is enough to show that  $\mathcal{P} \tilde{h} = - \tilde{h}_{skew}$. So, we compute:\\
\begin{align*}
(\mathcal{P} \tilde{h} )_{ij} &= \frac{1}{2} \tilde{h}_{ab} (\star \omega)_{ijab}\\
&= \frac{1}{2} (\nabla_u h_{va}) \psi^{+}_{uvb} (\star \omega)_{ijab}\\
&= \frac{1}{2} (\nabla_u h_{va}) (-\delta_{ui} \psi^{+}_{vja}-\delta_{uj}\psi^{+}_{iva}-\delta_{ua}\psi^{+}_{ijv}+\delta_{iv} \psi^{+}_{uja}+\delta_{jv} \psi^{+}_{iua}+\delta_{av} \psi^{+}_{iju} - \omega_{uv} \psi^{-}_{ija})\\
&= \frac{1}{2} (0+0-  (\nabla_a h_{va}) \psi^{+}_{ijv} + (\nabla_u h_{ia}) \psi^{+}_{uja} + (\nabla_u h_{ja}) \psi^{+}_{iua} + 0 - (\nabla_u h_{va})\omega_{uv} \psi^{-}_{ija})\\
&= \frac{1}{2} (-((\divv h) \intprod \psi^{+})_{ij} -  (\nabla_u h_{ai}) \psi^{+}_{uaj} + (\nabla_u h_{aj}) \psi^{+}_{uai} + (\divv h)_k \omega_{ka} \psi^{-}_{ija})\text{ (by Proposition~\ref{randomidentities})}\\
&= \frac{1}{2} (-((\divv h) \intprod \psi^{+})_{ij} - \tilde{h}_{ij} + \tilde{h}_{ji} + (\divv h)_k \psi^{+}_{kij})\\
&= \frac{1}{2} (-((\divv h) \intprod \psi^{+})_{ij} - 2(\tilde{h}_{skew})_{ij} + ((\divv h) \intprod \psi^{+})_{ij})\\
&= - (\tilde{h}_{skew})_{ij}.
\end{align*}
as claimed, concluding the proof.
\end{proof}
\end{prop}

\begin{prop}\label{harmonicthm} Let $M$ be compact nearly K\"{a}hler. Let $\beta \in \Omega^3_{12} \supseteq \mathcal{H}^3$. Hence, $\beta = h \diamond \psi^{+}$ for some unique $h\in \mathcal{S}^2_{-}.$ Then $\beta$ is harmonic iff $\divv h = 0, \tilde{h} = 2\omega h = -2 h \omega \in \mathcal{S}^2$.
\begin{proof}
Note that in fact, since $h$ is symmetric, $\omega$ is skew, and that they anticommute, we have $\omega h \in \mathcal{S}^2$. So the last condition is equivalent to $\tilde{h}_{symm} = 2 \omega h$ and $\tilde{h}_{skew} = 0.$\\
Since $M$ is compact, $\beta$ is harmonic if and only if $d^* \beta = 0$ and $d\beta = 0.$ Let us look at each of these conditions separately. So, we have:
\begin{align*}
0 &= -(d^* \beta)_{kl}\\
&= \nabla_j \beta_{jkl}\\
&= \nabla_j (h_{jp} \psi^{+}_{pkl} + h_{kp} \psi^{+}_{jpl}+h_{lp} \psi^{+}_{jkp})\\
&= (\nabla_j h_{jp}) \psi^{+}_{pkl} +(\nabla_j h_{kp}) \psi^{+}_{jpl}+ (\nabla_j h_{lp} )\psi^{+}_{jkp} + h_{jp} (\nabla_j  \psi^{+}_{pkl}) + h_{kp} (\nabla_j  \psi^{+}_{jpl})+h_{lp} (\nabla_j  \psi^{+}_{jkp})\\
&= (\divv h \intprod \psi^{+})_{kl}+ (\nabla_j h_{pk}) \psi^{+}_{jpl}- (\nabla_j h_{pl} )\psi^{+}_{jpk} + h_{jp} (\delta_{jp} \omega_{kl}+\delta_{jk} \omega_{lp}+\delta_{jl} \omega_{pk})+ 4 h_{kp} \omega_{pl} + 4 h_{lp} \omega_{kp}\text{ (by~\eqref{defnk3})}\\
&= (\divv h \intprod \psi^{+})_{kl} + \tilde{h}_{kl}-\tilde{h}_{lk}+(0+h_{kp} \omega_{lp} + h_{lp} \omega_{pk}) + 4 (h_{kp} \omega_{pl} + h_{lp} \omega_{kp})\\
&= (\divv h \intprod \psi^{+})_{kl} + 2 (\tilde{h}_{skew})_{kl}+ 3(h\omega + \omega h)_{kl}\\
&= (\divv h \intprod \psi^{+})_{kl} + 2 (\tilde{h}_{skew})_{kl},
\end{align*}
where we have used that $h \in \mathcal{S}^2_{-} \subseteq \mathcal{S}^2_{0}$. Recall that by Proposition~\ref{ushdug}, $\tilde{h}_{skew} \in \Omega^2_8$, hence, looking at the types we get:
\[d^* \beta = 0 \text{\tab if and only if\tab} \myarray{\divv h = 0, \\ \tilde{h}_{skew} = 0.}\]
Next, by Corollary~\ref{cor4diam}, we know that $d\beta = 0$ iff $\widehat{d\beta}= 0.$ We have:
\begin{align}
\widehat{d\beta}_{ia} &= (d\beta)_{ijkl} (\star \omega)_{ajkl}\nonumber\\
&= (\nabla_i \beta_{jkl}-\nabla_j \beta_{ikl}+\nabla_k \beta_{ijl}-\nabla_l \beta_{ijk}) (\star \omega)_{ajkl}\nonumber\\
&= (\nabla_i \beta_{jkl})(\star \omega)_{ajkl} -3 (\nabla_j \beta_{ikl}) (\star \omega)_{ajkl}.\label{eq:harm1}
\end{align}
We will compute each term of~\eqref{eq:harm1} separately. First we have:
\begin{align*}
(\nabla_i \beta_{jkl})(\star \omega)_{ajkl} &= \nabla_i (h_{jp} \psi^{+}_{pkl}+h_{kp} \psi^{+}_{jpl}+h_{lp} \psi^{+}_{jkp})(\star \omega)_{ajkl}\\
&=3 \nabla_i (h_{jp} \psi^{+}_{pkl}) (\star \omega)_{ajkl}\\
&=3 (\nabla_i h_{jp}) \psi^{+}_{pkl} (\star \omega)_{ajkl}+3 h_{jp} (\nabla_i \psi^{+}_{pkl}) (\star \omega)_{ajkl}\\
&=3 (\nabla_i h_{jp}) 2\psi^{+}_{paj}+3 h_{jp} (\delta_{ip} \omega_{kl} + \delta_{ik} \omega_{lp}+\delta_{il} \omega_{pk}) (\star \omega)_{ajkl}\\
&=0+3 h_{ij} \omega_{kl} (\star \omega)_{ajkl} + 6 h_{jp}\delta_{ik} \omega_{lp} (\star \omega)_{ajkl}\\
&=12 h_{ij} \omega_{aj}  - 6 h_{jp} \omega_{pl} (\star \omega)_{ajil}\\
&=12 h_{ij} \omega_{aj}  - 6 h_{jp} (\delta_{ap} \omega_{ji}+\delta_{jp} \omega_{ia}+\delta_{ip} \omega_{aj})\\
&=-12 (h \omega)_{ia}  - 6 h_{ja} \omega_{ji}- 6\trr(h) \omega_{ia}- 6 h_{ji} \omega_{aj}\\
&=-12 (h \omega)_{ia}  + 6  (\omega h+h \omega)_{ia}\\
&=-12 (h \omega)_{ia}.
\end{align*}
For the second term of~\eqref{eq:harm1}, we have:
\begin{align}
-3 (\nabla_j \beta_{ikl}) (\star \omega)_{ajkl} &= -3 \nabla_j (h_{ip} \psi^{+}_{pkl}+h_{kp} \psi^{+}_{ipl}+h_{lp} \psi^{+}_{ikp})(\star \omega)_{ajkl}\nonumber\\
&= -3 \nabla_j (h_{ip} \psi^{+}_{pkl})(\star \omega)_{ajkl} -6\nabla_j (h_{kp} \psi^{+}_{ipl})(\star \omega)_{ajkl}.\label{eq:harm2}
\end{align}
Here again, we compute both terms of~\eqref{eq:harm2} separately. First we have:
\begin{align*}
-3 \nabla_j (h_{ip} \psi^{+}_{pkl})(\star \omega)_{ajkl} &= -3 (\nabla_j h_{ip}) \psi^{+}_{pkl}(\star \omega)_{ajkl} - 3 h_{ip} (\nabla_j  \psi^{+}_{pkl})(\star \omega)_{ajkl}\\
&= -3 (\nabla_j h_{ip}) 2 \psi^{+}_{paj} - 3 h_{ip} (\delta_{jp} \omega_{kl} +\delta_{jk} \omega_{lp} +\delta_{jl} \omega_{pk})(\star \omega)_{ajkl}\\
&= -6 (\nabla_j h_{pi}) \psi^{+}_{jpa} - 3 h_{ip} \delta_{jp} \omega_{kl}(\star \omega)_{ajkl}\\
&= -6 \tilde{h}_{ia} - 3 h_{ij} 4 \omega_{aj}\\
&= -6 \tilde{h}_{ia} + 12 (h\omega)_{ia}.
\end{align*}
For the second term of~\eqref{eq:harm2} we use Proposition~\ref{randomidentities} to get:
\begin{align*}
-6\nabla_j (h_{kp} \psi^{+}_{ipl})(\star \omega)_{ajkl} &= -6(\nabla_j h_{kp}) \psi^{+}_{ipl}(\star \omega)_{ajkl} -6h_{kp}(\nabla_j \psi^{+}_{ipl}) (\star \omega)_{ajkl}\\
&= 6(\nabla_j h_{kp}) (\psi^{-}_{ipa}\omega_{jk}+\psi^{-}_{ipj}\omega_{ka}+\psi^{-}_{ipk}\omega_{aj}) - 6h_{kp}(\delta_{ji}\omega_{pl}+\delta_{jp}\omega_{li}+\delta_{jl}\omega_{ip}) (\star \omega)_{ajkl}\\
&= -6(\divv h)_s \omega_{sp} \psi^{-}_{ipa}+6 (\nabla_j h_{kp}) \psi^{-}_{ipj}\omega_{ka}+0 - 6h_{kp}\omega_{pl} (\star \omega)_{aikl}-6h_{kp}\omega_{li} (\star \omega)_{apkl}+0\\
&= 6(\divv h)_s  \psi^{-}_{iap} \omega_{sp}-6 (\nabla_j h_{pk}) \psi^{-}_{jpi}\omega_{ka}- 6h_{kp}(\delta_{ap}\omega_{ik}+\delta_{ip}\omega_{ka}+\delta_{kp}\omega_{ai})+0\\
&= 6(\divv h)_s  \psi^{+}_{ias}- 6(\tilde{h}_{ai} + 4(h\omega)_{ai})-6h_{ka} \omega_{ik} -6h_{ki}\omega_{ka} +0\\
&= 6(\divv h \intprod \psi^{+})_{ia} - 6\tilde{h}_{ai} - 24(h\omega)_{ia}-6(\omega h + h\omega)_{ia}\\
&= 6(\divv h \intprod \psi^{+})_{ia} - 6\tilde{h}_{ai} - 24(h\omega)_{ia}.
\end{align*}
Thus, combining the last two results we simplify~\eqref{eq:harm2} to get:
\begin{align*}
-3 (\nabla_j \beta_{ikl}) (\star \omega)_{ajkl} &= (-6 \tilde{h}_{ia} + 12 (h\omega)_{ia}) + (6(\divv h \intprod \psi^{+})_{ia} - 6\tilde{h}_{ai} -24(h\omega)_{ia})\\
&= 6(\divv h \intprod \psi^{+})_{ia} -12 (\tilde{h}_{symm})_{ia} -12 (h \omega)_{ia}.
\end{align*}
And so, returning to~\eqref{eq:harm1}, we have:
\begin{align*}
\widehat{d\beta}_{ia} &= -12(h\omega)_{ia} +  6(\divv h \intprod \psi^{+})_{ia} -12 (\tilde{h}_{symm})_{ia}-12 (h \omega)_{ia}\\
&= 6(\divv h \intprod \psi^{+})_{ia} -12 (\tilde{h}_{symm})_{ia} - 24 (h\omega)_{ia},
\end{align*}
which implies:
\[d \beta = 0 \text{\tab if and only if\tab} \myarray{\divv h = 0, \\ \tilde{h}_{symm} = -2h\omega.}\]
Hence, we get that $\beta$ is harmonic iff $\divv h = 0$ and $\tilde{h} = -2h\omega = 2 \omega h$.\qedhere
\end{proof}
\end{prop}
\begin{prop}\label{dfshugdsffds} Let $M$ be compact nearly K\"{a}hler. Let $\beta \in \Omega^2_{8} \supseteq \mathcal{H}^2$. Hence, $\beta = h \diamond \omega$ for some unique $h\in \mathcal{S}^2_{+0}.$ Then $\beta$ is harmonic iff $\divv h = 0, \tilde{h} = -3 h \omega \in \Omega^2_8.$
\begin{proof}
First, note that $$\beta_{ij} = (h \diamond \omega)_{ij} = h_{ip} \omega_{pj} + h_{jp} \omega_{ip} = (h\omega)_{ij} +  (\omega h)_{ij} = 2 (h\omega)_{ij}.$$
As in the proof of the previous theorem, $\beta$ is harmonic if and only if $d^{*} \beta = 0$ and $d \beta = 0$. Looking at each of the conditions separately, we get:
\begin{align*}
0 &= - (d^{*} \beta)_k = \nabla_p \beta_{pk} = 2 \nabla_p (h_{pu} \omega_{uk} )\\
&= 2 (\divv h)_u \omega_{uk} + 2h_{pu} \nabla_p \omega_{uk}\\
&= 2 (\divv h)_u \omega_{uk} - 2 h_{pu} \psi^+_{puk}\\
&= 2 (\divv h)_u \omega_{uk}.
\end{align*}
Since $\omega$ is non-degenerate, we get that:
\[d \beta = 0 \text{\tab if and only if\tab}\divv h = 0.\]
Next, by Corollary~\ref{diamcorol3}, we have that $d\beta = 0$ iff $\widehat{d\beta} = 0.$ We have:
\begin{align}
0 &= \widehat{d\beta}_{ia} = (d \beta)_{ijk} \psi^+_{ajk}\nonumber\\
&= (\nabla_i \beta_{jk}-\nabla_j \beta_{ik}+\nabla_k \beta_{ij}) \psi^+_{ajk}\nonumber\\
&= (\nabla_i \beta_{jk}) \psi^+_{ajk} - 2 (\nabla_j \beta_{ik})\psi^+_{ajk}.\label{ahn1}
\end{align}
We will compute each term of~\eqref{ahn1} separately. First we have:
\begin{align}
 (\nabla_i \beta_{jk}) \psi^+_{ajk} &= 2 \nabla_i (h_{ju} \omega_{uk})\psi^+_{ajk}\nonumber\\
&= 2 (\nabla_i h_{ju}) \omega_{uk} \psi^+_{ajk} + 2 h_{ju} (\nabla_i \omega_{uk}) \psi^+_{ajk}\nonumber\\
&= 2 (\nabla_i h_{ju}) (- \psi^-_{uaj}) + 2 h_{ju} (- \psi^+_{iuk}) \psi^+_{ajk}\nonumber\\
&= 0 - 2 h_{ju} (\delta_{ia} \delta_{uj} - \delta_{ij} \delta_{ua}-\omega_{ia} \omega_{uj} + \omega_{ij} \omega_{ua}) \text{\tab (by~\eqref{eq:pi3})}\nonumber\\
&=- 2 \delta_{ia} \trr{h} + 2 h_{ia} + 0 - 2 (\omega h \omega)_{ia}\nonumber\\
&= 2 h_{ia}  - 2 (h \omega^2)_{ia}\text{\tab (because $h \in \mathcal{S}^2_{+0}$)}\nonumber\\
&= 4 h_{ia}.\label{ahn2}
\end{align}
For the second term of~\eqref{ahn1} we have:
\begin{align}
- 2 (\nabla_j \beta_{ik})\psi^+_{ajk} &= -4 \nabla_j (h_{iu} \omega_{uk}) \psi^+_{ajk}\nonumber\\
&= -4 (\nabla_j h_{iu}) (-\psi^-_{uaj})-4h_{iu} (-\psi^+_{juk}) \psi^+_{ajk}\nonumber\\
&= 4 (\nabla_j h_{iu}) \psi^-_{uaj}-4h_{iu} (4 \delta_{ua})\text{\tab (by~\eqref{eq:pi4})}\nonumber\\
&= 4 (\nabla_j h_{iu}) \psi^-_{uaj}-16h_{ia}.\label{ahn3}
\end{align}
Conbining~\eqref{ahn2} and~\eqref{ahn3}, we get that:
\begin{align*}
0 &= \widehat{d\beta}_{ia} = 4 (\nabla_j h_{iu}) \psi^-_{uaj} - 12 h_{ia}.
\end{align*}
So, $d\beta = 0$ iff $ \widehat{d\beta} = 0$ iff $(\nabla_j h_{iu}) \psi^-_{uaj} = 3 h_{ia}$ iff $(\nabla_j h_{iu}) \psi^-_{uaj} \omega_{at} = 3 h_{ia} \omega_{at}$.\\
Since $(\nabla_j h_{iu}) \psi^-_{uaj} \omega_{at} = (\nabla_j h_{iu}) \psi^+_{ujt} = - \tilde{h}_{it}$, we get that:
\[d \beta = 0 \text{\tab if and only if\tab}\tilde{h} =  -3h\omega.\]
Hence, we conclude that $\beta$ is harmonic iff $\divv h = 0$ and $\tilde{h} = -3h\omega$.\\
It is easy to see that $h\omega \in \Omega^2$ for $h \in \mathcal{S}^2_{+0},$ so by Proposition~\ref{ushdug}, in this case we indeed have $\tilde{h} = -3h\omega \in \Omega^2_8.$
\end{proof}
\end{prop}
\subsection{Weitzenböck formulas}\label{weifojahsd}
The following formulas can be found in~\cite{CWMYK}, however we include the proofs, and when deriving sufficient conditions for vanishing of $b_2$ and $b_3$, we use slightly different forms of these formulas.
\subsubsection{$2$-forms}
We apply Corollary~\ref{weiz2simpl} to the nearly K\"{a}hler setting to get:
\begin{equation}\label{eq:weiz2forms2}
\Delta \beta =\nabla^{*} \nabla \beta + 8 \beta +  \hat{W}\beta, \text{ for any $\beta \in \Omega^2.$}
\end{equation}
\begin{prop}\label{thmnkffff1}
Let  $\beta = h \diamond \omega \in \Omega^2_8$ for some $h \in \mathcal{S}^2_{+0}$. Assume $\beta$ is harmonic. Then:
$$\nabla^* \nabla h + 6h + 2 \mathring{W} h = 0.$$
\begin{proof}
Using~\eqref{eq:weiz2forms2}, it is enough to show that $\nabla^* \nabla \beta = (\nabla^* \nabla h - 2h) \diamond \omega$ and $\hat{W} \beta = 2 (\mathring{W} h) \diamond \omega.$ So, we proceed with the first claim:
\begin{align*}
(\nabla^* \nabla \beta)_{ab} &= - \nabla_s \nabla_s \beta_{ab}\\
=& - \nabla_s \nabla_s (h \diamond \omega)_{ab}\\
=& - \nabla_s \nabla_s (h_{ap} \omega_{pb} + h_{bp} \omega_{ap})\\
=& -(( \nabla_s \nabla_s h) \diamond \omega)_{ab} - 4 (\nabla_s h_{ap}) (\nabla_s \omega_{pb})\\
&- 2 h_{ap} \nabla_s \nabla_s \omega_{pb} \text{ (because $h \in \mathcal{S}^2_{+0}$ and hence $h_{ap} \omega_{pb} = h_{bp} \omega_{ap}$)}\\
=& -(( \nabla_s \nabla_s h) \diamond \omega)_{ab} + 4 (\nabla_s h_{ap}) \psi^+_{spb} + 8 h_{ap} \omega_{pb} \text{ (by~\eqref{defnk2} and~\eqref{defnk3}, $\nabla_s \nabla_s \omega_{pb} = - 4 \omega_{pb}$)}\\
=& (( \nabla^* \nabla h) \diamond \omega)_{ab} + 4 \tilde{h}_{ab} + 8(h\omega)_{ab}\\
=& (( \nabla^* \nabla h) \diamond \omega)_{ab} - 12(h\omega)_{ab} + 8(h\omega)_{ab}\text{ (by Proposition~\ref{dfshugdsffds})}\\
=& (( \nabla^* \nabla h) \diamond \omega)_{ab} - 4(h\omega)_{ab}\\
=& (( \nabla^* \nabla h - 2h) \diamond \omega)_{ab}\text{ (becase $h \diamond \omega = 2h \omega$ for $h\in \mathcal{S}^2_{+0}$)}.\\
\end{align*}
For the second claim, we know that since $\beta \in \Omega^2_8$, then $\hat{W} \beta \in \Omega^2_8$. Hence by Proposition~\ref{diam2}, $\hat{W} \beta = f \diamond \omega$, for some $f \in \mathcal{S}^2_{+0}.$ The same proposition also tells us that $f = \frac{1}{2} (\hat{W} \beta)_{ik} \omega_{ak}$. Computing, we have:
\begin{align*}
f_{ia} &=  \frac{1}{2} (\hat{W} \beta)_{ik} \omega_{ak} = \frac{1}{2} W_{ikuv} \beta_{uv} \omega_{ak} = \frac{1}{2} W_{ikuv} (h \diamond \omega)_{uv} \omega_{ak} = W_{ikuv} h_{up} \omega_{pv} \omega_{ak}\\
&= -(W_{kuiv} + W_{uikv})h_{up} \omega_{pv} \omega_{ak}\text{\tab(by the Bianchi identity)}\\
&= (W_{ivuk} + W_{uivk}) \omega_{ka} h_{up} \omega_{vp}\\
&= W_{ivuk} \omega_{ka} h_{up} \omega_{vp} + W_{uivk} \omega_{ka} h_{up} \omega_{vp}\\
&= W_{ivak} \omega_{ku} h_{up} \omega_{vp} + W_{uiak} \omega_{kv} h_{up} \omega_{vp}\text{~(by Lemma~\ref{28omegacommute} and $W \in \Omega^2_8$ in first (last) two indices)}\\
&= W_{ivak} h_{ku} \omega_{up} \omega_{vp} - W_{uiak} h_{up} \delta_{kp}\text{~(we use that $h \in \mathcal{S}^2_{+0}$ and that $\omega^2 = -\operatorname{Id}$)}\\
&= W_{ivak} h_{ku} \delta_{uv} - W_{uiak} h_{uk}\\
&= W_{ivak} h_{kv} - W_{uiak} h_{uk}\\
&= W_{vika} h_{vk} + W_{uika} h_{uk}\\
&= 2(\mathring{W} h)_{ia},
\end{align*}
as claimed. Hence, the proof is complete.
\end{proof}
\end{prop}
\begin{thm}\label{thmnk1}
Let $M$ be a  compact nearly K\"{a}hler $6$-manifold. If $\mathcal{S}^2(\Omega^2_8) \ni \hat{W}  \geq -8$, or equivalenty $\mathcal{S}^2(\mathcal{S}^2_{+0}) \ni \mathring{W} \geq -4$, then $b_2 = 0.$
\begin{proof}
Let $\beta \in \Omega^2$ be harmonic. Then $\beta \in \Omega^2_8,$ as mentioned in the start of Section~\ref{harmshuidfh}. Substituting it in~\eqref{eq:weiz2forms2}, and using the assumption that $\hat{W} \geq -8$, we get that $\beta = 0,$ as there are no parallel non-zero $2$-forms.\\
Using the fact that $\hat{W} \beta = 2 (\mathring{W} h) \diamond \omega,$ where $\beta = h \diamond \omega \in \Omega^2_8$, for $h \in \mathcal{S}^2_{+0}$ we get the other equivalent condition.\\
Note that using Proposition~\ref{thmnkffff1} in order to get a similar result would have been worse, as we would have been able to only conclude that if $\mathcal{S}^2(\mathcal{S}^2_{+0}) \ni \mathring{W} \geq -3$ then $b_2 = 0.$ This is because  $\nabla^* \nabla \beta = (\nabla^* \nabla h - 2h) \diamond \omega$, so we can see that even though the left hand side is obviously non-negative, we cannot conclude that from the right hand side.
\end{proof}
\end{thm}
\begin{thm}\label{cornkk1}
Let $M$ be a  compact nearly K\"{a}hler $6$-manifold. Let $\delta \leq \bar{R} \leq \Delta$ with $-(\Delta + \delta) -\frac{7}{3} (\Delta - \delta) \geq -10$ \textbf{or} $(\Delta + \delta) -3 (\Delta - \delta) \geq -6$ . Then $b_2=0$.
\begin{proof} 
If the conditions above hold, then by Corollary~\ref{kjlakjkkj} we have that $\hat{R} \geq -10$. So, we use~\eqref{curidnkah} to get that $\hat{W} \geq -8$ and hence $b_2 = 0$ by Theorem~\ref{thmnk1}.
\end{proof}
\end{thm}
\subsubsection{$3$-forms}
We apply Corollary~\ref{weiz3simpl} to the nearly K\"{a}hler setting to get:
\begin{equation}\label{eq:weiz3forms2}
(\Delta \beta)_{abc} =(\nabla^{*} \nabla \beta)_{abc} + 9 \beta_{abc} +  W_{abpu} \beta_{puc}+ W_{acpu} \beta_{pbu}+ W_{bcpu} \beta_{apu}, \text{ for any $\beta \in \Omega^3.$}
\end{equation}

\begin{prop}\label{mainprop2}
Let $\beta \in \Omega^3_{12}$, so $\beta = h \diamond \psi^{+}$ for some unique $h \in \mathcal{S}^2_{-}$. Assume $\beta$ is harmonic. Then:
$$\nabla^* \nabla h + 8h + 2 \mathring{W} h = 0.$$
\begin{proof}
Substitute a harmonic $\beta$ into~\eqref{eq:weiz3forms2} to get the vanishing of the left hand side. Now, the goal is to rewrite the RHS as $A \diamond \psi^{+}$ for some $A \in \mathbb{R}g \oplus \mathbb{R}\omega \oplus \mathcal{S}^2_{-} \oplus \Omega^2_6.$ Then we can conclude that $A = 0$. By Proposition~\ref{harmonicthm}, since $\beta$ is harmonic, $\divv h =0$ and $\tilde{h} = 2\omega h.$ Keeping this in mind, we will simplify each term of the RHS of~\eqref{eq:weiz3forms2} one by one. We start with $\nabla^{*} \nabla \beta:$
\begin{align}\label{eq:3}
\begin{split}
(\nabla^{*} \nabla \beta)_{abc} =& -\nabla_s \nabla_s (h_{ap} \psi^{+}_{pbc} + h_{bp} \psi^{+}_{apc} + h_{cp} \psi^{+}_{abp})\\
=& ((\nabla^{*} \nabla h) \diamond \psi^{+})_{abc} - 2 ((\nabla_s h_{ap})(\nabla_s \psi^{+}_{pbc}) + (\nabla_s h_{bp}) (\nabla_s \psi^{+}_{apc}) + (\nabla_s h_{cp}) (\nabla_s \psi^{+}_{abp}))  \\
&- (h_{ap} (\nabla_s \nabla_s \psi^{+}_{pbc})+h_{bp} (\nabla_s \nabla_s \psi^{+}_{apc})+h_{cp} (\nabla_s \nabla_s \psi^{+}_{abp})).
\end{split}
\end{align}
First, note that:
\begin{align*}
\nabla_s \nabla_s \psi^{+}_{ijk} &= \nabla_s (\delta_{si} \omega_{jk}+\delta_{sj} \omega_{ki}+\delta_{sk} \omega_{ij})\\
&= \delta_{si} (-\psi^{+}_{sjk})+\delta_{sj} (-\psi^{+}_{ski})+\delta_{sk} (-\psi^{+}_{sij})\\
&= -3 \psi^{+}_{ijk}.
\end{align*}
Hence, the third term in~\eqref{eq:3} is equal to:
$$- (h_{ap} (\nabla_s \nabla_s \psi^{+}_{pbc})+h_{bp} (\nabla_s \nabla_s \psi^{+}_{apc})+h_{cp} (\nabla_s \nabla_s \psi^{+}_{abp})) = 3 (h_{ap} \psi^{+}_{pbc} + h_{bp} \psi^{+}_{apc} + h_{cp} \psi^{+}_{abp}) = (3 h \diamond \psi^{+})_{abc}.$$
In order to calculate the second term of~\eqref{eq:3}, we define the $3$-form $\sigma$ by:
$$\sigma_{abc} \coloneqq (\nabla_s h_{ap})(\nabla_s \psi^{+}_{pbc}) + (\nabla_s h_{bp}) (\nabla_s \psi^{+}_{apc}) + (\nabla_s h_{cp}) (\nabla_s \psi^{+}_{abp}). $$
We claim that $\sigma = 2h \diamond \psi^{+}$. In order to get this, we first calculate $\hat{\sigma}$ and then use Corollary~\ref{diamcorol3}. So,
\begin{align*}
\hat{\sigma}_{at} =& \sigma_{abc} \psi^{+}_{tbc}\\
=& ((\nabla_s h_{ap})(\nabla_s \psi^{+}_{pbc}) + (\nabla_s h_{bp}) (\nabla_s \psi^{+}_{apc}) + (\nabla_s h_{cp}) (\nabla_s \psi^{+}_{abp}))\psi^{+}_{tbc}\\
=& (\nabla_s h_{ap})(\nabla_s \psi^{+}_{pbc})\psi^{+}_{tbc}+2 (\nabla_s h_{bp}) (\nabla_s \psi^{+}_{apc})\psi^{+}_{tbc}\\
=& (\nabla_s h_{ap})(\delta_{sp}\omega_{bc}+\delta_{sb}\omega_{cp}+\delta_{sc}\omega_{pb})\psi^{+}_{tbc}+2 (\nabla_s h_{bp}) (\delta_{sa}\omega_{pc}+\delta_{sp}\omega_{ca}+\delta_{sc}\omega_{ap})\psi^{+}_{tbc}\\
=& 0+  (\nabla_b h_{ap}) \psi^{+}_{btc}\omega_{pc} +  (\nabla_c h_{ap}) \psi^{+}_{ctb}\omega_{pb}+2(\nabla_a h_{bp})\psi^{+}_{tbc}\omega_{pc}\\
&+2(\nabla_p h_{bp})\psi^{+}_{btc}\omega_{ac}-2(\nabla_c h_{bp})\psi^{+}_{cbt}\omega_{ap}\text{ (as $\divv h = 0$)}\\
=&- (\nabla_b h_{ap}) \psi^{-}_{btp}- (\nabla_c h_{ap}) \psi^{-}_{ctp}-2(\nabla_a h_{bp})\psi^{-}_{tbp}-0-2\tilde{h}_{pt}\omega_{ap}\\
=&- 2(\nabla_b h_{ap}) \psi^{-}_{btp}-0-2\tilde{h}_{pt}\omega_{ap}\\
=& 2(\nabla_b h_{pa}) \psi^{-}_{bpt}-2(\omega \tilde{h})_{at}\\
=& -2(\nabla_b h_{pa}) \psi^{+}_{bpu} \omega_{tu}-2(\omega \tilde{h})_{at}\\
=& -2\tilde{h}_{au} \omega_{tu}-2(\omega \tilde{h})_{at}\\
=& 2(\tilde{h} \omega)_{at}-2(\omega \tilde{h})_{at}\\
=& 4(\omega h \omega)_{at}-4(\omega^2 h)_{at}\text{ (because $\tilde{h} = 2 \omega h$)}\\
=& -4(h \omega^2)_{at}+4h_{at}\\
=& 8h_{at}.\\
\end{align*}
Hence, $\hat{\sigma} = 8h \in \mathcal{S}^2_0$. Thus, by Proposition~\ref{diamcorol3}, $\sigma = \frac{1}{4}\hat{\sigma} \diamond \psi^{+} = 2h \diamond \psi^{+}$, as claimed.\\
Thus, returning to~\eqref{eq:3}, we get:
\begin{align}
\nabla^{*} \nabla \beta =& (\nabla^{*} \nabla h) \diamond \psi^{+}-2 \sigma + 3 h \diamond \psi^{+}\nonumber\\
=& (\nabla^{*} \nabla h -4h+3h) \diamond \psi^{+}\nonumber\\
=& (\nabla^{*} \nabla h -h) \diamond \psi^{+}.\label{eq:dsfs}
\end{align}
Next, we proceed to the terms in~\eqref{eq:weiz3forms2} with the Weyl tensors. Recall that $W$ is in $\Omega^2_8$ with respect to the first two or the last two indices. Hence, $W_{abij} \psi^{+}_{abc} = 0.$ So, we have:
\begin{align*}
W_{abpu} \beta_{puc} &= W_{abpu} (h_{ps} \psi^{+}_{suc}+h_{us} \psi^{+}_{psc}+h_{cs} \psi^{+}_{pus})\\
&= 2  h_{ps} W_{abpu}\psi^{+}_{suc}.
\end{align*}
Similarly, we have:
\begin{align*}
W_{acpu} \beta_{pbu} &= 2  h_{ps} W_{acpu}\psi^{+}_{sbu},\\
W_{bcpu} \beta_{apu} &= 2  h_{ps} W_{bcpu}\psi^{+}_{asu}.\\
\end{align*}
Thus, the Weyl terms in~\eqref{eq:weiz3forms2} are equal to: 
$$W_{abpu} \beta_{puc}+ W_{acpu} \beta_{pbu}+ W_{bcpu} \beta_{apu} = 2h_{ps} (W_{abpu}\psi^{+}_{suc} + W_{acpu}\psi^{+}_{sbu}+W_{bcpu}\psi^{+}_{asu}).$$
Now, we define the $3$-form $\gamma$ via $$\gamma_{abc} \coloneq h_{ps} (W_{abpu}\psi^{+}_{suc} + W_{acpu}\psi^{+}_{sbu}+W_{bcpu}\psi^{+}_{asu}).$$
We claim that $\gamma = (\mathring{W} h) \diamond \psi^{+}.$ Again, to get this, we first need to calculate $\hat{\gamma}$. We have:
\begin{align}
\hat{\gamma}_{at} &= \gamma_{abc} \psi^{+}_{tbc}\nonumber\\
&= h_{ps} (W_{abpu}\psi^{+}_{suc} + W_{acpu}\psi^{+}_{sbu}+W_{bcpu}\psi^{+}_{asu}) \psi^{+}_{tbc}\nonumber\\
&= 2 h_{ps} W_{abpu}\psi^{+}_{suc}\psi^{+}_{tbc} + 0\nonumber\\
&= 2 h_{ps} W_{abpu}(\delta_{st}\delta_{ub}- \delta_{sb}\delta_{ut}+\omega_{ut}\omega_{sb}+\omega_{bu}\omega_{st})\nonumber\\
&= 0 - 2 h_{pb}W_{abpt}+2 h_{ps}W_{abpu}\omega_{ut}\omega_{sb} + 2h_{ps} W_{abpu} \omega_{bu}\omega_{st}.\label{eq:weizterms}
\end{align}
We will simplify the last two terms of~\eqref{eq:weizterms} separately. Recall Lemma~\ref{28omegacommute} which implies that $W$ and $\omega$ commute. Also, we have that $h$ and $\omega$ anticommute. Hence, for the third term we have:
\begin{align*}
2 h_{ps}W_{abpu}\omega_{ut}\omega_{sb} &= 2 W_{abpu} \omega_{ut}  h_{ps} \omega_{sb}\\
&= -2 \omega_{pu} W_{abut}  \omega_{ps} h_{sb}\\
&= -2 \delta_{us} W_{abut} h_{sb}\\
&= -2 W_{abst} h_{sb}\\
&= 2 W_{bast} h_{bs}\\
&= 2 (\mathring{W} h)_{at}.
\end{align*}
For the fourth term of~\eqref{eq:weizterms} we have:
\begin{align*}
2h_{ps} W_{abpu} \omega_{bu}\omega_{st} &= -2 W_{abpu}\omega_{ub} h_{ps}\omega_{st}\\
&= -2 \omega_{pu}W_{abub} h_{ps}\omega_{st}\\
&= 0,
\end{align*}
because $W_{abub} = 0$.
Thus, returning to~\eqref{eq:weizterms}, we get:
\begin{align*}
\hat{\gamma}_{at} &= - 2 h_{pb}W_{abpt} + 2 (\mathring{W} h)_{at} + 0\\
&= 2 W_{abtp} h_{bp} + 2 (\mathring{W} h)_{at} + 0\\
&= 4 (\mathring{W} h)_{at}.
\end{align*}
Hence, $\hat{\gamma} = 4 \mathring{W} h \in \mathcal{S}^2_0$. Thus, by Proposition~\ref{diamcorol3}, $\gamma = \frac{1}{4}\hat{\gamma} \diamond \psi^{+} = (\mathring{W} h) \diamond \psi^{+}$, as claimed.\\
This finishes the proof of the proposition, as substituting all the results into~\eqref{eq:weiz3forms2}, we get:
$$0 =((\nabla^{*} \nabla h -h) + 9h + 2\gamma)\diamond \psi^{+}= (\nabla^{*} \nabla h  + 8h + 2\mathring{W} h)\diamond \psi^{+}.\qedhere$$
\end{proof}
\end{prop}
\begin{thm}\label{thmnk2}
Let $M$ be a  compact nearly K\"{a}hler $6$-manifold. If $\mathcal{S}^2(\mathcal{S}^2_-) \ni \mathring{W}  \geq -\frac{9}{2}$ \textbf{or} $\mathcal{S}^2(\Omega^2_8) \ni \hat{W} \geq -3$, then $b_3 = 0.$
\begin{proof}
The first statement follows from the fact that using~\eqref{eq:dsfs} we can rewrite Proposition~\ref{mainprop2} as: if $\beta = h\diamond \psi^+$ is harmonic for some $h \in \mathcal{S}^2_-$, then:
$$0 = \nabla^* \nabla \beta + (9h + 2 \mathring{W} h) \diamond \psi^+.$$
Hence, assuming $\mathcal{S}^2(\mathcal{S}^2_-) \ni \mathring{W}  \geq -\frac{9}{2}$ and using the fact that there are no nonzero parallel $h \in \mathcal{S}^2_0$, we get $b_3 = 0.$\\
Note that using Proposition~\ref{mainprop2} in order to get a similar result would have been worse, as we would have been able to only conclude that if $\mathcal{S}^2(\mathcal{S}^2_{-}) \ni \mathring{W} \geq -4$ then $b_3 = 0.$ This is because  $\nabla^* \nabla \beta = (\nabla^* \nabla h - h) \diamond \psi^+$, so we can see that even though the left hand side is obviosuly non-negative, we cannot conclude that from the right hand side.\\
Next, $\mathcal{S}^2(\Omega^2_8) \ni \hat{W} \geq -3$ implies $b_3 = 0$ because of~\eqref{eq:weiz3forms2}.\\
Note that the condition $\mathcal{S}^2(\Omega^2_8) \ni \hat{W} \geq -3$ is weaker than the condition $\mathcal{S}^2(\mathcal{S}^2_-) \ni \mathring{W}  \geq -\frac{9}{2}$. This is because in the proof of Proposition~\ref{mainprop2} we show that $W_{abpu} \beta_{puc}+ W_{acpu} \beta_{pbu}+ W_{bcpu} \beta_{apu} = (2(\mathring{W} h) \diamond \psi^+)_{abc}$, for $\beta = h \diamond \psi^+ \in \Omega^3_{12}$, where $h \in \mathcal{S}^2_-.$ That means if we assume that $\mathcal{S}^2(\Omega^2_8) \ni \hat{W} \geq c$, then $\mathring{W} \geq \frac{3c}{2}$, where $c \in \mathbb{R},$ but not vice versa.
\end{proof}
\end{thm}
\begin{thm}\label{thmnk3}
Let $M$ be a  compact nearly K\"{a}hler $6$-manifold. Let $\delta \leq \bar{R} \leq \Delta$ with $\delta \geq \frac{1}{4}$ \textbf{or} $\Delta \leq \frac{17}{8}$. Then $b_3=0$.
\begin{proof}
Recall that the Einstein constant $k = 5.$ Then by Theorem~\ref{nkmaincor}, on $\mathcal{S}^2_0$, $\mathring{R} \geq -5 + 6\delta$ and $\mathring{R} \geq 5 - 4\Delta$. Hence, by~\eqref{curidnkah}, $\mathring{W} \geq -6 + 6\delta$ and $\mathring{W} \geq 4 - 4\Delta$.\\
In order for $b_3 = 0$, by Theorem~\ref{thmnk2} we want $\mathring{W}  \geq -\frac{9}{2}$. We have $-6 + 6\delta \geq -\frac{9}{2}$  iff $\delta \geq \frac{1}{4}$; and $4 -4 \Delta \geq -\frac{9}{2}$  iff $\Delta \leq \frac{17}{8}$. Hence, the result follows. Recall, that a priori, by Remark~\ref{uselessfact}, we have that $\delta \leq 1 \leq \Delta$.\\ Also, note that we do not use Corollary~\ref{kjlakjkkj} along with the statement that $\mathcal{S}^2(\Omega^2_8) \ni \hat{W} \geq -3$ implies that $b_3 = 0.$ This is because the sufficient conditions in terms of the bounds on the sectional curvature we would have obtained imply that $\Delta \leq \frac{17}{8}$ or $\delta \geq \frac{1}{4}$.
\end{proof}
\end{thm}

\section{Examples}

We only consider normal homogeneous spaces $G/H$ (see~\cite{JD}.) Denote the Lie algebras of $G$ and $H$ by $\mathfrak{g}$ and $\mathfrak{h}$ respectively. Let $\mathfrak{m}$ be the orthogonal complement of $\mathfrak{h}$ in $\mathfrak{g}$.\\
Having a bi-invariant metric on $G$ induces a metric on $G/H$ which gives us a Riemannian submersion $\pi:G \rightarrow G/H.$ The usual decomposition into vertical and horizontal subspaces corresponds to the decomposition $\mathfrak{g} = \mathfrak{m} \oplus \mathfrak{h}$. Hence, using the formula (3.30) and $\text{Corollary 3.19} $ from~\cite{JD} gives us that  for $X,Y,Z,W \in \mathfrak{m}$ we have:
\begin{align*}
R(X,Y,Z,W) =& \frac{1}{4}(\langle [X,W], [Y,Z]\rangle- \langle [X,Z], [Y,W]\rangle)+ \frac{1}{4}(\langle [X,W]_\mathfrak{h}, [Y,Z]_\mathfrak{h}\rangle- \langle [X,Z]_\mathfrak{h}, [Y,W]_\mathfrak{h}\rangle)\\
&- \frac{1}{2}\langle [Z,W]_{\mathfrak{h}}, [X,Y]_{\mathfrak{h}}\rangle.
\end{align*}
Letting $X=W, Y=Z$ yields:
\begin{equation}\label{eqseccur}
R(X,Y,Y,X) = \frac{1}{4}\|[X,Y]_{\mathfrak{m}}\|^2+\|[X,Y]_{\mathfrak{h}}\|^2.
\end{equation}
The first formula will allow us to calculate sharp bounds for $\mathring{R}$, $\hat{R}$ and the second one bounds for $\bar{R}$, which we use to check the theorems.

Before going to specific examples, we briefly outline the process of how we get the bounds for $\hat{R}$ and $\mathring{R}.$\\
Consider $\hat{R}$ first. Note that this is a self-adjoint operator, hence it is bounded by the smallest and the largest eigenvalues. So, if we take any local orthonormal frame $f_{\alpha}$ of $\Omega^2$, find all the entries of the matrix $\hat{R}_{\alpha \beta}$ corresponding to this linear operator, we can find its eigenvalues. \\
We already have that $e_{i} \wedge e_j$ for $i < j$ is an orthonormal frame for $\Omega^2$. Let $f_{\alpha}=e_i \wedge e_j,f_{\beta}= e_u \wedge e_v$ be any two such basis elements. Then from the proof of Theorem~\ref{hatthm1}, we have $$\bar{R}_{\alpha \beta} = (\hat{R} f_{\alpha}, f_{\beta}) = \frac{1}{2}(\hat{R} f_{\alpha})_{kl} (f_{\beta})_{kl} =  \frac{1}{2}(\hat{R} (e_i \wedge e_j))_{kl} (e_u \wedge e_v)_{kl} =(\hat{R} (e_i \wedge e_j))_{uv} = 2 R_{ijuv}.$$ So, we use Maple to find all the values $R_{ijuv}$ and thus the matrix $\hat{R}_{\alpha \beta}$. As mentioned before, its largest and smallest eigenvalues are the sharp bounds we are looking for.\\

Next, we want to also find the bounds for $\mathring{R}$ on $\mathcal{S}^2$. As before, it is enough to find the eigenvalues corresponding to this linear self-adjoint operator.\\
Let $M$ be of dimension $n$. Then $\mathcal{S}^2$ has dimension $\frac{n(n+1)}{2}$.\\
Let $\{e_1,\dots, e_n\}$ be an orthonormal frame, with the dual frame $\{e^1, \dots, e^n\}.$ 
Now, let: $$f_{ij} \coloneqq e_i \otimes e_j + e_j \otimes e_i \text{, for $i<j$},$$
$$f_{ii} \coloneqq e_i \otimes e_i \text{, for $i = 1, \dots, n$}.$$
Note that these $f_{ij}, f_{ii}$ form a frame for $\mathcal{S}^2$. Let they be denoted just as $f_{\alpha}$. We will still specify if the $f_{\alpha}$ we take is one of $f_{ij}$, for $i<j$, or one of $f_{ii}.$ So, now, we want to find the matrix representation of $\mathring{R}$ in terms of the basis of $f_{\alpha}$'s. Note that this frame is not orthonormal, but we do not need it to be, since the eigenvalues of the matrix will turn out to be all the same.\\
First, note that for $h \in \mathcal{S}^2$ we have:
$$h =\sum_{i,j=1}^{n} h_{ij} e_i \otimes e_j = \sum_{i<j}h_{ij} (e_i \otimes e_j + e_j \otimes e_i) + \sum_{i=1}^{n} h_{ii} e_i \otimes e_i = \sum_{i<j} h_{ij} f_{ij} + \sum_{i=1}^{n}h_{ii}f_{ii}.$$

That means that the $f_{\beta}$ component of $h$, which we will denote by $h^{\beta}$ is equal to $h_{ij}$, for $f_{\beta} = f_{ij}, i\leq j.$
Next, we need to find how $\mathring{R}$ acts on these basis elements $f_{\alpha}$. We claim that:
\begin{align}
(\mathring{R} f_{ij})_{ab} &= R_{iajb}+R_{jaib}, \text{ for $i<j$},\label{aaa1}\\
(\mathring{R} f_{ii})_{ab} &= R_{iaib}.\label{aaa2}
\end{align}
We calculate:
\begin{align*}
(\mathring{R} f_{ij})_{ab} =& \sum_{k,l} R_{kalb} (f_{ij})_{kl}\\
=& \sum_{k,l} R_{kalb} (e_i \otimes e_j + e_j \otimes e_i)_{kl}\\
=& \sum_{k,l} R_{kalb} (\delta_{ik}\delta_{jl}+\delta_{jk}\delta_{il})\\
=& R_{iajb} + R_{jaib}.
\end{align*}
Also,
\begin{align*}
(\mathring{R}f_{ii})_{ab} =& (\mathring{R} (e_i \otimes e_i ))_{ab}= R_{iaib}.
\end{align*}
as claimed.

From~\eqref{aaa1} and ~\eqref{aaa2} it is easy to get that the $f_{\beta}$ component of $\mathring{R} f_{\alpha}$, which we denote by $\mathring{R}_{\alpha \beta}$ is equal to:
\begin{itemize}
\item if $f_{\alpha} = f_{ij}, f_\beta = f_{st}, i<j$ and $s<t$, then $\mathring{R}_{\alpha \beta} = R_{isjt} +R_{jsit},$
\item if $f_{\alpha} = f_{ij}, f_\beta = f_{ss}, i<j$, then $\mathring{R}_{\alpha \beta} = 2R_{isjs},$
\item if $f_{\alpha} = f_{ii}, f_\beta = f_{st}, s<t$, then $\mathring{R}_{\alpha \beta} = R_{isit},$
\item if $f_{\alpha} = f_{ii}, f_\beta = f_{ss}$, then $\mathring{R}_{\alpha \beta} = R_{isis}.$
\end{itemize}
Note that we actually need bounds for $\hat{R}$ or $\mathring{R}$ on specific subspaces of $\Omega^2$ or $\mathcal{S}^2$ respectively, but this will be easy to get as we know what these operators do on the complements of the subspaces we are looking for.\\
Also, for both examples we identify $\mathfrak{su}(2)$ with $\mathbb{R}^3$ as follows: $\frac{a_1}{\sqrt{2}} I + \frac{a_2}{\sqrt{2}} J+\frac{a_3}{\sqrt{2}} K \longleftrightarrow (a_1,a_2,a_3)$, where $I = \begin{pmatrix}  i & 0\\ 0 & -i\end{pmatrix}, J=\begin{pmatrix}  0 & -1\\ 1 & 0\end{pmatrix}, K =  \begin{pmatrix}  0 & i\\ i & 0\end{pmatrix} $ is the standard basis for $\mathfrak{su}(2)$. This takes the inner product $\trr(a^{*} b)$ on $\mathfrak{su}(2)$ to the usual one in $\mathbb{R}^3$. For $a,b,c,d \in \mathfrak{su}(2)$, it is straightforward to verify that:
\begin{equation}\label{ex22}
\begin{split}
\langle[a,b],[c,d]\rangle &= 2(\langle a,c\rangle\langle b,d\rangle-\langle a,d\rangle\langle b,c\rangle),\\
|[a,b]|^2 &= 2 |a|^2 |b|^2 - 2\langle a,b\rangle^2.
\end{split}
\end{equation}

\subsection{$\frac{\su(3)\times \su(2)}{\uu(1)\times \su(2)}$}
We describe some of the aspects of the nearly $G_2$ structure on this $G/H$. See~\cite{AS} for more information.
By $\su(2)_d$ we denote the following embedding of $\su(2)$ into $\su(3)\times \su(2)$:
$$\su(2)_d = \{(\begin{pmatrix}
A & 0\\
0 & 0
\end{pmatrix},A), A \in \su(2)\}.$$
Also, by $\uu(1)$ we mean the following embedding into subgroup of $\su(3) \times \{\operatorname{I}\} \subseteq \su(3)\times \su(2)$:
$$\uu(1) = \{(\begin{pmatrix}
e^{it} & 0 & 0\\
0 & e^{it} & 0\\
0 & 0 & e^{-2it}
\end{pmatrix}, \ii), t \in \mathbb{R}\}.$$
Then $\frac{\su(3)\times \su(2)}{\uu(1)\times \su(2)}$ is a normal homogeneous space with the metric $B =-\frac{1}{24}(6 \trr(uv))+4 \trr(wz))$ (this is a multiple of the Killing form), for $(u,w), (v,z) \in \mathfrak{g} = \mathfrak{su}(3) \oplus \mathfrak{su}(2).$ With such a choice of a metric one obtains a nearly $G_2$ structure with $\tau_0 = -\frac{12}{\sqrt{5}}$ and hence the Einstein constant $k = \frac{54}{5}$ with $R = \frac{378}{5}.$ Then we have the following orthogonal decomposition:
$$\mathfrak{g} = \mathfrak{h} \oplus \mathfrak{m},$$
with 
$$\mathfrak{h} = \mathfrak{u}(1) \oplus \mathfrak{su}(2)_d \text{ and } \mathfrak{m} = \mathfrak{su}(2)_o \oplus \mathfrak{m}{'}$$
where:
$$\mathfrak{u}(1) = \operatorname{span}\{(\begin{pmatrix}
i & 0 & 0\\
0 & i & 0\\
0 & 0 & -2i
\end{pmatrix},0)\}, \mathfrak{su}(2)_d = \{(\begin{pmatrix}
a & 0\\
0 & 0
\end{pmatrix},a), a \in \mathfrak{su}(2)\},$$ $$\mathfrak{su}(2)_o = \{(\begin{pmatrix}
2a & 0\\
0 & 0
\end{pmatrix},-3a), a \in \mathfrak{su}(2)\}, \mathfrak{m}{'} = \{(\begin{pmatrix}
0 & z\\
-\bar{z}^{T} & 0
\end{pmatrix},0), z \in \mathbb{C}^2\}.$$
We define the following quantities:
\begin{align*}
f_1(a)& \coloneqq (\begin{pmatrix}
a & 0\\
0 & 0
\end{pmatrix},a) \in \mathfrak{su}(2)_d \subseteq \mathfrak{h}, \text{ for } a \in \mathfrak{su}(2),\\
f_2(a)& \coloneqq (\begin{pmatrix}
2a & 0\\
0 & 0
\end{pmatrix},-3a)\in \mathfrak{su}(2)_o \subseteq \mathfrak{m}, \text{ for } a \in \mathfrak{su}(2),\\
g_1(r) & \coloneqq (r \begin{pmatrix}
i & 0 & 0\\
0 & i & 0\\
0 & 0 & -2i
\end{pmatrix},0) \in \mathfrak{u}(1) \subseteq \mathfrak{h}, \text{ for } r \in \mathbb{R},\\
g_2(z)& \coloneqq (\begin{pmatrix}
0 & z\\
-\bar{z}^{T} & 0
\end{pmatrix},0)\in \mathfrak{m}{'} \subseteq \mathfrak{m}, \text{ for } z \in \mathbb{C}^2,\\
|z|^2&\coloneqq|z_1|^2+|z_2|^2 = \bar{z}^{T} z, \text{ for } z = \begin{pmatrix}
z_1 \\
z_2
\end{pmatrix} \in \mathbb{C}^2.
\end{align*}
Note that all $f_1,f_2,g_1,g_2$ are linear. Next, we compute their norms with respect to the metric $B$, where the norm squared is denoted by $\|\cdot \|^2 = B(\cdot, \cdot).$ So:
\begin{align}
\|f_1(a)\|^2 &= -\frac{1}{24} (6 \trr(\begin{pmatrix}
a & 0\\
0 & 0
\end{pmatrix}^2) + 4 \trr(a^2))= -\frac{1}{24} (6\trr(a^2) + 4 \trr(a^2))= -\frac{5}{12} \trr(a^2)= \frac{5}{12} |a|^2.\nonumber\\
\|f_2(a)\|^2 &= -\frac{1}{24} (6 \trr(\begin{pmatrix}
2a & 0\\
0 & 0
\end{pmatrix}^2) + 4 \trr((-3a)^2))= -\frac{1}{24} (24\trr(a^2) + 36 \trr(a^2))= -\frac{5}{2} \trr(a^2)= \frac{5}{2} |a|^2.\label{something}\\
\|g_1(r)\|^2 &=-\frac{1}{24} (6 r^2 \trr(\begin{pmatrix}
i & 0 & 0\\
0 & i & 0\\
0 & 0 & -2i
\end{pmatrix}^2))= -\frac{1}{24} (6 r^2 (-6))= \frac{3}{2} r^2.\nonumber\\
\|g_2(z)\|^2 &= -\frac{1}{24} (6 \trr(\begin{pmatrix}
0 & z\\
-\bar{z}^{T} & 0
\end{pmatrix}^2))=-\frac{1}{4} \trr(\begin{pmatrix}
-z\bar{z}^{T} & 0\\
0 & -\bar{z}^{T}z
\end{pmatrix})= \frac{1}{4} (\trr(z\bar{z}^{T}) + \bar{z}^{T}z)= \frac{1}{4} 2\bar{z}^{T}z= \frac{1}{2} |z|^2.\nonumber
\end{align}
Now, we want to find bounds on $\bar{R}$, so take $X = f_2(a) + g_2(z), Y = f_2(b) + g_2(w) \in \mathfrak{m}$ for some $a,b \in \mathfrak{su}(2)$ and $z,w \in \mathbb{C}^2$, with $\|X\|^2 = \|Y\|^2 = 1, B(X,Y) = 0.$ So, we have:
\begin{align*}
1 &= \|X\|^2 = \|f_2(a)\|^2 + \|g_2(z)\|^2 =  \frac{5}{2} |a|^2 + \frac{1}{2} |z|^2,\\
1 &= \|Y\|^2 = \|f_2(b)\|^2 + \|g_2(w)\|^2 =  \frac{5}{2} |b|^2 + \frac{1}{2} |w|^2,\\
0 &= B(X,Y) = B( f_2(a) + g_2(z),  f_2(b) + g_2(w))\\
&= B((\begin{pmatrix}
2a & z\\
-\bar{z}^{T} & 0
\end{pmatrix}, -3a),(\begin{pmatrix}
2b & w\\
-\bar{w}^{T} & 0
\end{pmatrix}, -3b))\\
&= -\frac{1}{24} (6 \trr(\begin{pmatrix}
2a & z\\
-\bar{z}^{T} & 0
\end{pmatrix} \begin{pmatrix}
2b & w\\
-\bar{w}^{T} & 0
\end{pmatrix}) + 4 \trr((-3a)(-3b)))\\
&= -\frac{1}{24} (6 \trr(\begin{pmatrix}
4ab - z \bar{w}^{T} & 2aw\\
-2\bar{z}^{T}b & -\bar{z}^{T}w
\end{pmatrix}) + 36 \trr(ab))\\
&= -\frac{1}{24} (24 \trr(ab) - 6\trr(z \bar{w}^{T}) -6\bar{z}^{T}w + 36 \trr(ab))\\
&= -\frac{1}{24} (60 \trr(ab) - 6\bar{w}^{T}z -6\bar{z}^{T}w),\\
&\text{and thus:}\\
& \bar{w}^{T}z + \bar{z}^{T}w = 10 \trr(ab).
\end{align*}
Next, we need to calculate $[X,Y]$. We have:
\begin{align*}
[X,Y] &= [f_2(a) + g_2(z),  f_2(b) + g_2(w)]\\
&= [f_2(a), f_2(b)] +  [g_2(z),f_2(b)] + [f_2(a),g_2(w)] +  [ g_2(z),g_2(w)].
\end{align*}
We will calculate each term separately:
\begin{align*}
[f_2(a),  f_2(b)] &= [(\begin{pmatrix}
2a & 0\\
0 & 0
\end{pmatrix},-3a), (\begin{pmatrix}
2b & 0\\
0 & 0
\end{pmatrix},-3b)]\\
&=  (\begin{pmatrix}
4[a,b] & 0\\
0 & 0
\end{pmatrix},9[a,b])\\
&= 6 f_1([a,b]) - f_2([a,b]).
\end{align*}
Next:
\begin{align*}
[g_2(z),  f_2(b)] &= [(\begin{pmatrix}
0 & z\\
-\bar{z}^{T} & 0
\end{pmatrix},0), (\begin{pmatrix}
2b & 0\\
0 & 0
\end{pmatrix},-3b)]\\
&= ([\begin{pmatrix}
0 & z\\
-\bar{z}^{T} & 0
\end{pmatrix},\begin{pmatrix}
2b & 0\\
0 & 0
\end{pmatrix}],0)\\
&= (\begin{pmatrix}
0 & 0\\
-2\bar{z}^{T}b & 0
\end{pmatrix} - \begin{pmatrix}
0 & 2bz\\
0 & 0
\end{pmatrix},0)\\
&= (\begin{pmatrix}
0 & -2bz\\
-2\bar{z}^{T}b & 0
\end{pmatrix} ,0)\\
&=-2 g_2(bz).
\end{align*}
Similarly:
\begin{align*}
[f_2(a), g_2(w)] = - [g_2(w),  f_2(a)] = 2 g_2(aw).
\end{align*}
Finally:
\begin{align*}
[g_2(z),g_2(w)] &= ([\begin{pmatrix}
0 & z\\
-\bar{z}^{T} & 0
\end{pmatrix},\begin{pmatrix}
0 & w\\
-\bar{w}^{T} & 0
\end{pmatrix}],0)\\
&=(\begin{pmatrix}
-z\bar{w}^{T} & 0\\
0 & -\bar{z}^{T} w
\end{pmatrix}-\begin{pmatrix}
-w\bar{z}^{T} & 0\\
0 & -\bar{w}^{T} z
\end{pmatrix},0)\\
&=(\begin{pmatrix}
-z\bar{w}^{T} +  w\bar{z}^{T}& 0\\
0 & -\bar{z}^{T} w + \bar{w}^{T} z
\end{pmatrix},0)\\
&= (\frac{-\bar{z}^{T} w + \bar{w}^{T} z}{-2i}\begin{pmatrix}
i & 0 & 0\\
0 & i & 0\\
0 & 0 & -2i
\end{pmatrix},0) + (\begin{pmatrix}
-z\bar{w}^{T} + w \bar{z}^T + \frac{-\bar{z}^{T} w + \bar{w}^{T} z}{2} \operatorname{I} & 0\\
0 & 0
\end{pmatrix},0)\\
&\text{(Let $A \coloneqq -z\bar{w}^{T} + w \bar{z}^T + \frac{-\bar{z}^{T} w + \bar{w}^{T} z}{2} \operatorname{I} \in \mathfrak{su}(2)$)}\\
&= g_1 (\frac{-\bar{z}^{T} w + \bar{w}^{T} z}{-2i})+ (\begin{pmatrix}
A & 0\\
0 & 0
\end{pmatrix},0)\\
&= g_1 (\frac{-\bar{z}^{T} w + \bar{w}^{T} z}{-2i}) + \frac{3}{5}f_1(A) +  \frac{1}{5}f_2(A).
\end{align*}
Hence we conclude that:
$$[X,Y] = [X,Y]_{\mathfrak{m}} +  [X,Y]_{\mathfrak{h}},$$
where
$$[X,Y]_{\mathfrak{m}} = f_2(-[a,b]+ \frac{1}{5}A) + g_2(2(aw-bz)),$$
$$[X,Y]_{\mathfrak{h}} = f_1(6[a,b]+ \frac{3}{5}A) + g_1 (\frac{-\bar{z}^{T} w + \bar{w}^{T} z}{-2i}).$$
Applying the formula~\eqref{eqseccur} for the sectional curvature, along with~\eqref{something}, we get:
\begin{align*}
\bar{R}(X \wedge Y) =& \frac{1}{4}\|[X,Y]_{\mathfrak{m}}\|^2 + \|[X,Y]_{\mathfrak{h}}\|^2\\
=& \frac{1}{4} (\| f_2(-[a,b]+ \frac{1}{5}A)\|^2 + \| g_2(2(aw-bz))\|^2) + \| f_1(6[a,b]+ \frac{3}{5}A)\|^2 + \| g_1 (\frac{-\bar{z}^{T} w + \bar{w}^{T} z}{-2i})\|^2\\
=& \frac{1}{4} (\frac{5}{2} |-[a,b]+ \frac{1}{5}A|^2 + \frac{1}{2} |2(aw-bz)|^2) + \frac{5}{12} |6[a,b]+ \frac{3}{5}A|^2 + \frac{3}{2} \Big(\frac{-\bar{z}^{T} w + \bar{w}^{T} z}{-2i}\Big)^2\\
=& \frac{5}{8}(|[a,b]|^2 + \frac{1}{25}|A|^2 + \frac{2}{5} \trr([a,b]A))+\frac{1}{2}|aw-bz|^2\\
&+ \frac{5}{12}(36|[a,b]|^2 + \frac{9}{25}|A|^2 - \frac{36}{5} \trr([a,b]A)) - \frac{3}{8} (-\bar{z}^{T} w + \bar{w}^{T} z)^2\\
=& \frac{125}{8}|[a,b]|^2 + \frac{7}{40} |A|^2-\frac{11}{4}\trr([a,b]A)+\frac{1}{2}|aw-bz|^2- \frac{3}{8} (-\bar{z}^{T} w + \bar{w}^{T} z)^2.
\end{align*}

It is straightforward to check that for $a \in \mathfrak{su}(2), w \in \mathbb{C}^2$ we have: 
\begin{equation}\label{asdaaf}
|aw|^2 = \frac{1}{2} |a|^2 |w|^2.
\end{equation}
Polarizing, we also get:
$$\langle aw, az \rangle = \frac{1}{2} |a|^2 \langle w,z \rangle,$$
$$\langle az, bz \rangle + \langle bz, az \rangle= |z|^2 \langle a,b \rangle,$$
$$\langle az, bw \rangle + \langle bw, az \rangle + \langle bz, aw \rangle + \langle aw, bz \rangle=  \langle a,b \rangle(\langle z,w \rangle + \langle w,z \rangle).$$

For simplicity, define:
\begin{align*}
\alpha & \coloneqq \langle z, w \rangle \in \mathbb{C},\\
\sigma & \coloneqq -\langle aw, bz \rangle \in \mathbb{C},\\
\varphi & \coloneqq -\langle az, bw \rangle \in \mathbb{C}.
\end{align*}
Then some calculation yields that:
\begin{align*}
\frac{125}{8}|[a,b]|^2 &= \frac{125}{4} |a|^2 |b|^2 - \frac{125}{4}\langle a, b \rangle^2,\\
\frac{7}{40} |A|^2 &= \frac{7}{20} |z|^2|w|^2 - \frac{7}{80}\alpha^2 - \frac{7}{80}\bar{\alpha}^2 - \frac{7}{40}\alpha \bar{\alpha},\\
-\frac{11}{4}\trr([a,b]A) &= \frac{11}{4}(\sigma + \bar{\sigma} - \varphi - \bar{\varphi}),\\
\frac{1}{2}|aw-bz|^2 &= \frac{1}{4}|a|^2|w|^2+\frac{1}{4}|b|^2|z|^2 - \sigma - \bar{\sigma},\\
- \frac{3}{8} (-\bar{z}^{T} w + \bar{w}^{T} z)^2 &=- \frac{3}{8}\alpha^2 - \frac{3}{8}\bar{\alpha}^2 + \frac{3}{4}\alpha \bar{\alpha}.
\end{align*}
Recall that we assumed:
\begin{align*}
1 &= \frac{5}{2} |a|^2 + \frac{1}{2} |z|^2,\\
1 &= \frac{5}{2} |b|^2 + \frac{1}{2} |w|^2,\\
& \alpha + \bar{\alpha} = -10 \langle a, b \rangle .
\end{align*}

Isolating $\langle a, b \rangle, |a|^2, |b|^2$ and substituting these resullts into expressions found earlier, we get that:
\begin{align*}
\bar{R}(X \wedge Y) =&~ 5 - \frac{12}{5} |z|^2 - \frac{12}{5} |w|^2 + \frac{3}{2} |z|^2 |w|^2 - \frac{31}{40} \alpha^2 - \frac{31}{40}\bar{\alpha}^2 - \frac{1}{20}\alpha \bar{\alpha}\\
& + \frac{7}{4} \sigma + \frac{7}{4} \bar{\sigma} - \frac{11}{4} \varphi - \frac{11}{4} \bar{\varphi} .
\end{align*}

Note that each $\sigma, \bar{\sigma}, \varphi, \bar{\varphi}$ in absolute value is $\leq \frac{1}{2} |a||b||z||w|$, by Cauchy-Schwarz and~\eqref{asdaaf}. Hence:
\begin{align*}
\bar{R}(X \wedge Y) &\leq  5 - \frac{12}{5} |z|^2 - \frac{12}{5} |w|^2 + \frac{3}{2} |z|^2 |w|^2 + \frac{3}{2} |\alpha|^2 + \frac{9}{2} |a||b||z||w|\\
&\leq  5 - \frac{12}{5} |z|^2 - \frac{12}{5} |w|^2 + 3 |z|^2 |w|^2 + \frac{9}{2} |a||b||z||w|.
\end{align*}
One can check that on $1 = \frac{5}{2} |a|^2 + \frac{1}{2} |z|^2, 1 = \frac{5}{2} |b|^2 + \frac{1}{2} |w|^2$, $$\bar{R}(X \wedge Y) \leq \frac{37}{5}.$$

Numerical evidence suggests that $\frac{1}{5} \leq \bar{R}(X \wedge Y)$, however the author was unable to verify this. Nevertheless, we have $0 \leq \bar{R}(X \wedge Y)$ and we can show that both values $\frac{1}{5}$ and $\frac{37}{5}$ can be achieved:
$$a=b=0, z = \sqrt{2} \begin{pmatrix} 1 \\ 0 \end{pmatrix}, w = \sqrt{2} \begin{pmatrix} i \\ 0 \end{pmatrix}$$
gives orthonormal $X,Y$ with $\bar{R}(X \wedge Y) = \frac{37}{5}$, and
$$a=-\frac{\sqrt{5}}{5} I, b=0, z = 0, w = \sqrt{2} \begin{pmatrix} 0 \\ i \end{pmatrix}$$
gives orthonormal $X,Y$ with $\bar{R}(X \wedge Y) = \frac{1}{5}$.

A computation on Maple reveals that eigenvalues of $\hat{R}$ on $\Omega^2$ are $-\frac{114}{5}_1, -\frac{66}{5}_3, -\frac{18}{5}_7,\frac{6}{5}_{10}$ where by the subscript we denote its multiplicity. Note that this makes sense, because for $\beta \in \Omega^2_7$, from Remark~\ref{g2description}, $\beta = X \intprod \varphi$ and hence a quick calculation gives that $\hat{R} \beta = X \intprod (\hat{R} \varphi) = X \intprod (-\frac{\tau_0^2}{8} \varphi) = -\frac{\tau_0^2}{8} \beta = -\frac{18}{5} \beta,$  so we get seven eigenvalues $-\frac{18}{5}$. Hence, we conlude that on $\Omega^2_{14}$, $-\frac{114}{5} \leq \hat{R} \leq \frac{6}{5}.$\\
Similarly, Maple shows that eigenvalues of $\mathring{R}$ on $\mathcal{S}^2$ are $-\frac{54}{5}_1,-\frac{47}{5}_1,-\frac{23}{5}_7, \frac{13}{5}_8, 5_5, \frac{37}{5}_6.$ Again, this makes sense, as we know that $\mathring{R}g = -\frac{3\tau^2_0}{8} g = -\frac{54}{5} g$. Hence, on $\mathcal{S}^2_0$,$-\frac{47}{5} \leq \mathring{R} \leq \frac{37}{5}.$\\
So, we summarize and check the theorems. We have:
\begin{align*}
\frac{1}{5} &\leq \bar{R} \leq \frac{37}{5},\\
-\frac{114}{5} &\leq \hat{R} \leq \frac{6}{5}\text{ on $\Omega^2_{14}$},\\
-\frac{47}{5} &\leq \mathring{R} \leq \frac{37}{5}\text{ on $\mathcal{S}^2_0$}.
\end{align*}
Corollary~\ref{improvedaa} gives us that $-\frac{122}{5} \leq \hat{R} \leq \frac{46}{5}$ on $\Omega^2_{14}$ which is correct.\\
Corollary~\ref{nkmaincor} gives us that $-\frac{47}{5} \leq \mathring{R} \leq \frac{49}{5}$ on $\mathcal{S}^2_0$, which is also correct with the first inequality being sharp.\\
For the main Theorems, we know in this case that $b_2 = 1$. So it must be false that $\hat{W} \geq -18$ on $\Omega^2_{14}$, by Theorem~\ref{thmm2}. By~\eqref{hatlemmag2}, we have that $\hat{W} \geq -19.2$ on $\Omega^2_{14}$, with the eigenvalue value $-19.2$ achieved. Hence, we get no contradiction.\\
As for Theorem~\ref{thmm1} we cannot predict whether $\mathring{W} \geq -\frac{54}{5}$ on $\mathcal{S}^2_0$ or $\hat{W} \geq -\frac{36}{5}$ on $\Omega^2_{14}$ must hold or not, because $b_3 =0.$ However, these inequalities do not hold as we have $\mathring{W} \geq -\frac{56}{5}$ on $\mathcal{S}^2_0$ and $\hat{W} \geq -\frac{96}{5}$ on $\Omega^2_{14}$ with these lower bounds attained. This shows that, in general, these sufficient conditions are not necessary.
\subsection{$\frac{\su(2) \times \su(2) \times \su(2)}{\su(2)}$}\label{dsusdhf}

First, we describe the nearly K\"{a}hler structure on this $G / H.$\\
The $\su(2)$ in the denominator is embedded diagonally in the numerator, meaning it is: $$\{(A,A,A) \in \su(2) \times \su(2) \times \su(2):A \in \su(2)\}.$$ Recall that we decompose $\mathfrak{g} = \mathfrak{h} \oplus \mathfrak{m}$. Then $\mathfrak{m} = \{(a,b,-(a+b): a,b \in \mathfrak{su}(2))\}.$\\ Equipping $G$ with the metric $B((a,b,c),(u,v,w)) = \frac{1}{3} (\trr(a^*u)+\trr(b^*v)+\trr(c^*w))$, for $(a,b,c),(u,v,w) \in \mathfrak{g}$ makes $G/H$ into a normal homogeneous space with scalar curvature $30.$\\
The nearly K\"{a}hler structure is obtained from the almost complex structure, which is defined as follows: $$J((a,b,c)) = \frac{2}{\sqrt{3}} (b,c,a) +  \frac{1}{\sqrt{3}} (a,b,c),$$ for $(a,b,c) \in \mathfrak{m}.$ It follows from $c = - a - b$ that $J^2 = -\operatorname{I}.$\\
Define the following quantities:
\begin{align*}
f(a) &\coloneqq (a,a,a) \in \mathfrak{h} \subseteq \mathfrak{g}, \text{ for } a\in \mathfrak{su}(2)\\
g(b,c) &\coloneqq (b,c,-(b+c)) \in \mathfrak{m} \subseteq \mathfrak{g}, \text{ for } b,c\in \mathfrak{su}(2).
\end{align*}
Then for $(a,b,c) \in \mathfrak{g}$, we have that
\begin{equation}\label{ex21}
\begin{split}
(a,b,c)_{ \mathfrak{h}} &= f(\frac{a+b+c}{3}),\\
(a,b,c)_{ \mathfrak{m}} &= g(\frac{2a-b-c}{3},\frac{-a+2b-c}{3}).
\end{split}
\end{equation}
Note that also $|f(a)|^2 = |a|^2,$ and $|g(b,c)|^2 = \frac{1}{3}(|b|^2+|c|^2+|b+c|^2)$.\\
We want to calculate the bounds on $\bar{R}$. Clearly from the formula~\eqref{eqseccur} for $\bar{R}$, we see that $0 \leq \bar{R}$. We claim that $\bar{R} \leq \frac{9}{4}$. Take $X,Y \in \mathfrak{m}$ with $\|X\|^2 = 1 = \|Y\|^2, B(X,Y) = 0.$ Let $X = g(b,c), Y = g(d,e),$ for $b,c,d,e \in \mathfrak{su}(2).$ Then:
\begin{align*}
[X,Y] &= [(b,c,-(b+c)),(d,e,-(d+e))]\\
&=([b,d],[c,e],[b,d]+[c,e]+[b,e]+[c,d]).
\end{align*}
Let $A \coloneqq [b,d], B \coloneqq [c,e], C \coloneqq [b,e] + [c,d]$ so that $[X,Y] = (A,B,A+B+C).$ By~\eqref{ex21}:
\begin{align*}
[X,Y]_{\mathfrak{h}} &= f(\frac{1}{3}(2A+2B+C)),\\
[X,Y]_{\mathfrak{m}} &= g(\frac{1}{3}(A-2B-C),\frac{1}{3}(-2A+B-C)).
\end{align*}
Hence, equation~\eqref{eqseccur} gives:
\begin{align*}
\bar{R}(X \wedge Y) =& \frac{1}{4}|[X,Y]_{\mathfrak{m}}|^2+|[X,Y]_{\mathfrak{h}}|^2\\
=&  \frac{1}{4} \cdot \frac{1}{3} (|\frac{1}{3}(A-2B-C)|^2+|\frac{1}{3}(-2A+B-C)|^2 + |\frac{1}{3}(-A-B-2C)|^2)+|\frac{1}{3}(2A+2B+C)|^2\\
=&  \frac{1}{108} |A-2B-C|^2+\frac{1}{108} |-2A+B-C|^2+\frac{1}{108} |-A-B-2C|^2+\frac{1}{9} |2A+2B+C|^2\\
=& \frac{1}{108}(|A|^2+4|B|^2+|C|^2-4\langle A,B\rangle-2\langle A,C\rangle+4\langle B,C\rangle)\\
&+ \frac{1}{108}(4|A|^2+|B|^2+|C|^2-4\langle A,B\rangle+4\langle A,C\rangle-2\langle B,C\rangle)\\
&+ \frac{1}{108}(|A|^2+|B|^2+4|C|^2+2\langle A,B\rangle+4\langle A,C\rangle+4\langle B,C\rangle)\\
&+ \frac{1}{9}(4|A|^2+4|B|^2+|C|^2+8\langle A,B\rangle+4\langle A,C\rangle+4\langle B,C\rangle)\\
=& \frac{1}{2}|A|^2+\frac{1}{2}|B|^2+\frac{1}{6}|C|^2+\frac{5}{6}\langle A,B\rangle+\frac{1}{2}\langle A,C\rangle+\frac{1}{2}\langle B,C\rangle.
\end{align*}
Using~\eqref{ex22} we now get:
\begin{align*}
|A|^2 &= |[b,d]|^2 = 2 |b|^2 |d|^2 - 2\langle b,d\rangle^2,\\
|B|^2 &= |[c,e]|^2 = 2 |c|^2 |e|^2 - 2\langle c,e\rangle^2,\\
|C|^2 &= |[b,e]+[c,d]|^2 =  |[b,e]|^2+ |[c,d]|^2 + 2\langle[b,e],[c,d]\rangle\\
&= 2 |b|^2 |e|^2 - 2\langle b,e\rangle^2+2 |c|^2 |d|^2 - 2\langle c,d\rangle^2 +4(\langle b,c\rangle\langle e,d\rangle-\langle b,d\rangle\langle c,e\rangle),\\
\langle A,B\rangle &= \langle [b,d], [c,e]\rangle = 2(\langle b,c\rangle\langle d,e\rangle-\langle b,e\rangle\langle c,d \rangle),\\
\langle A,C\rangle &= \langle [b,d],[b,e]+[c,d]\rangle = \langle [b,d], [b,e]\rangle + \langle [b,d], [c,d]\rangle\\
&= 2 (|b|^2 \langle d,e\rangle-\langle b,e\rangle \langle b,d\rangle+ \langle b,c\rangle|d|^2-\langle b,d\rangle \langle c,d\rangle),\\
\langle B,C\rangle &= \langle [c,e],[b,e]+[c,d]\rangle = \langle [c,e], [b,e]\rangle + \langle [c,e], [c,d]\rangle\\
&= 2 (\langle c,b\rangle|e|^2-\langle c,e\rangle \langle b,e\rangle+ |c|^2\langle d,e\rangle-\langle c,d\rangle \langle c,e\rangle).
\end{align*}
Substituting the above into $\bar{R}(X\wedge Y)$, we get:
\begin{equation}\label{ex23}
\begin{split}
\bar{R}(X\wedge Y) =& |b|^2 |d|^2 - \langle b,d\rangle^2 + |c|^2 |e|^2 - \langle c,e\rangle^2\\
&+\frac{1}{3} |b|^2 |e|^2 - \frac{1}{3}\langle b,e\rangle^2+\frac{1}{3} |c|^2 |d|^2 - \frac{1}{3}\langle c,d\rangle^2 +\frac{2}{3}\langle b,c\rangle\langle e,d\rangle-\frac{2}{3}\langle b,d\rangle\langle c,e\rangle + \frac{5}{3}\langle b,c\rangle\langle d,e\rangle- \frac{5}{3}\langle b,e\rangle\langle c,d \rangle\\
&+|b|^2 \langle d,e\rangle-\langle b,e\rangle \langle b,d\rangle+ |d|^2 \langle b,c\rangle-\langle b,d\rangle \langle c,d\rangle + |e|^2\langle c,b\rangle-\langle c,e\rangle \langle b,e\rangle+ |c|^2\langle d,e\rangle-\langle c,d\rangle \langle c,e\rangle.
\end{split}
\end{equation}
Recall that we assumed $\|X\|^2 = 1 = \|Y\|^2, B(X,Y) = 0$. Hence $1 = \|X\|^2 \|Y\|^2-  B(X,Y)^2.$ (Note this is just saying $\|X \wedge Y\|^2 = 1$. We could have assumed just this, however, the first assumption makes the argument easier.) We have:
\begin{align}
\|X\|^2 &= \frac{1}{3}(|b|^2+|c|^2 + |b+c|^2) = \frac{2}{3}(|b|^2+|c|^2 + \langle b,c\rangle).\nonumber\\
\|Y\|^2 &= \frac{1}{3}(|d|^2+|e|^2 + |d+e|^2) = \frac{2}{3}(|d|^2+|e|^2 + \langle d,e\rangle).\nonumber\\
B(X,Y) &= \frac{1}{3}(\langle b,d\rangle+\langle c,e\rangle+\langle b+c,d+e\rangle)\nonumber\\
&= \frac{1}{3}(2\langle b,d\rangle+2\langle c,e\rangle+\langle b,e\rangle+\langle c,d\rangle).\label{ex24}
\end{align}
Hence,
\begin{align*}
1 =& \|X\|^2 \|Y\|^2-B(X,Y)^2\\
=&  \frac{4}{9}(|b|^2+|c|^2 + \langle b,c\rangle)(|d|^2+|e|^2 + \langle d,e\rangle) - \frac{1}{9}(2\langle b,d\rangle+2\langle c,e\rangle+\langle b,e\rangle+\langle c,d\rangle)^2.
\end{align*}
or equivalently,
\begin{align*}
\frac{9}{4} =& |b|^2 |d|^2+|b|^2 |e|^2 + |b|^2 \langle d,e\rangle + |c|^2 |d|^2+|c|^2 |e|^2 + |c|^2 \langle d,e\rangle +|d|^2\langle b,c\rangle+ |e|^2\langle b,c\rangle + \langle b,c\rangle \langle d,e\rangle\\
&-\langle b,d\rangle^2-\langle c,e\rangle^2-\frac{1}{4}\langle b,e\rangle^2 -\frac{1}{4}\langle c,d\rangle^2- 2\langle b,d\rangle \langle c,e\rangle-\langle b,d\rangle \langle b,e\rangle - \langle b,d\rangle \langle c,d\rangle\\
&-\langle c,e\rangle \langle b,e\rangle - \langle c,e\rangle \langle c,d\rangle-\frac{1}{2}\langle b,e\rangle\langle c,d\rangle,
\end{align*}
which can be rearranged to get:
\begin{align*}
&|b|^2 |d|^2 + |b|^2 \langle d,e\rangle +|c|^2 |e|^2 + |c|^2 \langle d,e\rangle +|d|^2\langle b,c\rangle+ |e|^2\langle b,c\rangle -\langle b,d\rangle^2-\langle c,e\rangle^2\\ -&\langle b,d\rangle \langle b,e\rangle - \langle b,d\rangle \langle c,d\rangle -\langle c,e\rangle \langle b,e\rangle - \langle c,e\rangle \langle c,d\rangle\\
&=\frac{9}{4}-|b|^2 |e|^2 - |c|^2 |d|^2 - \langle b,c\rangle \langle d,e\rangle + \frac{1}{4}\langle b,e\rangle^2+\frac{1}{4}\langle c,d\rangle^2 + 2\langle b,d\rangle \langle c,e\rangle+\frac{1}{2}\langle b,e\rangle\langle c,d\rangle.
\end{align*}
Substituting this into~\eqref{ex23}, we get:
\begin{align*}
\bar{R}(X,Y) &= \frac{9}{4}-\frac{2}{3}|b|^2 |e|^2 -\frac{2}{3} |c|^2 |d|^2 +\frac{4}{3} \langle b,c\rangle \langle d,e\rangle - \frac{1}{12}\langle b,e\rangle^2-\frac{1}{12}\langle c,d\rangle^2 + \frac{4}{3}\langle b,d\rangle \langle c,e\rangle-\frac{7}{6}\langle b,e\rangle\langle c,d\rangle\\
&= \frac{9}{4}-\frac{1}{12}\Bigg(8\Big(|b|^2 |e|^2 + |c|^2 |d|^2 - 2 \langle b,c\rangle \langle d,e\rangle - 2 \langle b,d\rangle \langle c,e\rangle+ 2\langle b,e\rangle\langle c,d\rangle\Big) + (\langle b,e\rangle-\langle c,d\rangle)^2\Bigg).
\end{align*}
We claimed that $\bar{R} \leq \frac{9}{4}$, hence it is enough to show that $$|b|^2 |e|^2 + |c|^2 |d|^2 - 2 \langle b,c\rangle \langle d,e\rangle - 2 \langle b,d\rangle \langle c,e\rangle+ 2\langle b,e\rangle\langle c,d\rangle \geq 0.$$
We note that this expression is equal to:
\begin{align}
&(|b|^2 |e|^2 - \langle b,e\rangle^2 + |c|^2 |d|^2 - \langle c,d\rangle^2 - 2 \langle b,c\rangle \langle d,e\rangle + 2 \langle b,d\rangle\langle c,e\rangle) + (\langle b,e\rangle^2 +2\langle b,e\rangle \langle c,d\rangle + \langle c,d\rangle^2) - 4\langle b,d\rangle\langle c,e\rangle \nonumber\\
&= (\frac{1}{2}|[b,e]|^2 +  \frac{1}{2}|[c,d]|^2 - \langle [b,e],[c,d]\rangle) + (\langle b,e\rangle+\langle c,d\rangle)^2 - 4\langle b,d\rangle\langle c,e\rangle\nonumber\\
&= \frac{1}{2}|[b,e]-[c,d]|^2+(\langle b,e\rangle+\langle c,d\rangle)^2 - 4\langle b,d\rangle\langle c,e\rangle.\label{ex25}
\end{align}
We assumed that $B(X,Y)=0$, so from~\eqref{ex24}, we get that $\langle b,e\rangle+\langle c,d\rangle = -2 (\langle b,d\rangle+\langle c,e\rangle)$. Thus, continuing with~\eqref{ex25}, we get:
$$\frac{1}{2}|[b,e]-[c,d]|^2+4(\langle b,d\rangle+\langle c,e\rangle)^2 - 4\langle b,d\rangle\langle c,e\rangle,$$
which is always non-negative because for any real $x,y$, we have $4(x+y)^2-4xy = 4(x^2+xy+y^2) \geq 0.$\\
Finally, we need to show that the bounds $0 \leq \bar{R} \leq \frac{9}{4}$ are sharp. To do this, we take an explicit orthonormal basis for $\mathfrak{m}$:
$$e_1 = g(\frac{\sqrt{3}}{2}I,0), \tab e_2 = g(\frac{\sqrt{3}}{2}J,0),\tab e_3 = g(\frac{\sqrt{3}}{2}K,0),$$$$e_4 = g(\frac{1}{2}I,-I),\tab e_5 = g(\frac{1}{2}J,-J),\tab e_6 = g(\frac{1}{2}K,-K).$$
For this basis we also have: $Je_{i} = e_{i+3}, 1\leq i \leq 3.$\\
Then one easily calculates that $\bar{R}(e_1 \wedge e_4) = 0, \bar{R}(e_1 \wedge e_2) = \frac{9}{4}$, as we claimed.\\
A computation on Maple reveals that eigenvalues of $\hat{R}$ on $\Omega^2$ are $-7_3, -2_7, 1_5$, where by the subscript we denote its multiplicity. Note that this makes sense, as we know from the discussion following the proof of Proposition~\ref{curiden} that $\hat{R}$ is $-2\operatorname{Id}$ on $\Omega^2_1$ and $\Omega^2_6$, so we get seven eigenvalues $-2$. Hence, we conlude that on $\Omega^2_8$, $-7 \leq \hat{R} \leq 1.$\\
Similarly, Maple shows that eigenvalues of $\mathring{R}$ on $\mathcal{S}^2$ are $-5_1, -4_2,-\frac{3}{2}_3, 2_{10},\frac{5}{2}_5.$ Again, this makes sense, as we know that $\mathring{R}g = -5 g$, because the Einstein constant is $5$. Hence, on $\mathcal{S}^2_0$,$-4 \leq \mathring{R} \leq \frac{5}{2}.$
Furthermore, using that $\hat{W} \beta = 2(\mathring Wh) \diamond \omega$ for $\beta = h \diamond \omega \in \Omega^2_8$, where $h \in \mathcal{S}^2_{+-}$, it can be easily shown that in fact the eigenvalues $-\frac{3}{2}_3, \frac{5}{2}_5$ occur on $\mathcal{S}^2_{+0}$ and $-4_2,2_{10}$ occur on $\mathcal{S}^2_{-}.$
So, we summarize and check the theorems. We have:
\begin{align*}
0 &\leq \bar{R} \leq \frac{9}{4},\\
-7 &\leq \hat{R} \leq 1\text{ on $\Omega^2_8$},\\
-4 &\leq \mathring{R} \leq \frac{5}{2}\text{ on $\mathcal{S}^2_0$},\\
-\frac{3}{2} &\leq \mathring{R} \leq \frac{5}{2}\text{ on $\mathcal{S}^2_{+0}$},\\
-4 &\leq \mathring{R} \leq 2\text{ on $\mathcal{S}^2_{-}$}.
\end{align*}
Corollary~\ref{kjlakjkkj} gives us that $-\frac{15}{2} \leq \hat{R} \leq 3$ on $\Omega^2_8$ which is correct.\\
Corollary~\ref{nkmaincor} gives us that $-4 \leq \mathring{R} \leq 5$ on $\mathcal{S}^2_0$, which is correct with the first inequality being sharp.\\ 
Corollary~\ref{sdgfhgffgsdccc} gives us that $-\frac{7}{4} \leq \mathring{R} \leq \frac{7}{2}$ on $\mathcal{S}^2_{+0}$, which is also correct.\\
For the main theorems, we know in this case that $b_3 = 2$. So it must be false that $\mathring{W} \geq -4$ on $\mathcal{S}^2_0$ , by Theorem~\ref{thmnk2}. By~\eqref{curidnkah}, we have that $-5 \leq \mathring{W}$ on $\mathcal{S}^2_0$, with the eigenvalue value $-5$ achieved. Hence, we get no contradiction. Similarly, $\hat{R}$ achieves $-7$, so $\hat{W}$ achieves $-5$, hence we indeed have that $\hat{W} \geq -3$ is false.\\
As for Theorem~\ref{thmnk1} we cannot predict whether $\hat{W} \geq -8$ on $\Omega^2_8$ (or equivalenty $\mathcal{S}^2(\mathcal{S}^2_{+0}) \ni \mathring{W} \geq -4$) must hold or not, because $b_2 =0.$ However, we can actually deduce the vanishing of $b_2$, since by~\eqref{curidnkah}, we can get that $-5 \leq \hat{W}$ on $\Omega^2_8$, so the assumption of the theorem is satisfied. Finally, note that it is even possible to deduce that $b_2=0$ from Theorem~\ref{cornkk1}, since we get that $-(\Delta + \delta) -\frac{7}{3} (\Delta - \delta) \geq -10$ and $(\Delta + \delta) -3 (\Delta - \delta) \geq -6$ both hold.

\end{document}